\documentclass[preprint,12pt]{elsarticle}
\usepackage{todonotes}
\usepackage{float}
\usepackage{amsmath}
\usepackage{amssymb}
\usepackage{amsthm}
\usepackage{mathtools}
\usepackage{pifont}
\usepackage{graphicx}
\usepackage{xcolor}
\usepackage{fullpage}
\usepackage[percent]{overpic}

\newcommand{\nodes}{{\mathcal A}}
\newcommand{\vphi}{\varphi}
\graphicspath{{./Figures/}}
\usepackage{algpseudocode}
\usepackage{algorithm}
\usepackage[percent]{overpic}

\newcommand{\pad}[2]{\frac{\partial #1}{\partial #2}}
\newcommand{\RR}{\mathbb{R}}
\newcommand\norm[1]{\lVert#1\rVert}
\newcommand\abs[1]{\lvert#1\rvert}
\newcommand\normh[1]{\norm{#1}_{h,\alpha}}

\renewcommand{\vec}[1]{\mathbf{#1}}

\newtheorem{theorem}{Theorem}[section]
\newtheorem{corollary}{Corollary}[theorem]
\newtheorem{lemma}[theorem]{Lemma}
\newtheorem{proposition}[section]{Proposition}

\usepackage{amssymb}

\begin{document}

\begin{frontmatter}



\title{A nodal ghost method based on variational formulation and  regular square grid for elliptic problems on arbitrary domains in two space dimensions}


\author[1,2]{Clarissa Astuto}
\author[2,3]{Daniele Boffi}
\author[1]{Giovanni Russo}
\author[4]{Umberto Zerbinati}

\affiliation[1]{organization={Department of Mathematics and Computer Science, University of Catania},
            city={Catania},
            postcode={95125},
            country={Italy}}

\affiliation[2]{organization={King Abdullah University of Science and Technology (KAUST)},
            city={Thuwal},
            postcode={23955},
            country={Saudi Arabia}}

\affiliation[3]{organization={Dipartimento di Matematica ``F. Casorati'', University of Pavia},
            city={Pavia},
            postcode={27100},
            country={Italy}}

\affiliation[4]{organization={Mathematical Institute, University of Oxford},
        city={Oxford},
        postcode={OX2 6GG},
        country={England}}

\begin{abstract}
This paper focuses on the numerical solution of elliptic partial differential equations (PDEs) with Dirichlet and mixed boundary conditions, specifically addressing the challenges arising from irregular domains. Both finite element method (FEM) and finite difference method (FDM), face difficulties in dealing with arbitrary domains. The paper introduces a novel nodal symmetric ghost {method based on a variational formulation}, which combines the advantages of FEM and FDM. The method employs bilinear finite elements on a structured mesh and provides a detailed implementation description. A rigorous a priori convergence rate analysis is also presented. The convergence rates are validated with many numerical experiments, in both one and two space dimensions.

\end{abstract}




\end{frontmatter}



\section{Introduction}
Elliptic partial differential equations (PDEs) are used to describe a large variety of physical phenomena with applications ranging from modelling the diffusion of a certain chemical species in a solvent to determining the electric field corresponding to a given charge density.
In the interest of keeping our presentation as streamlined as possible, we will only focus our attention on the simplest elliptic PDE, i.e.,\ the Poisson equation with mixed boundary conditions,
\begin{equation}
	\label{eq:poissonMixed}
	-\Delta u = f \; {\rm in }\;\Omega, \; u=g_D\;{\rm on }\;\Gamma_D \subset \partial\Omega \; {\rm and } \;\pad{u}{n}  = g_N \; {\rm on }\;\Gamma_N\subset\partial\Omega.
\end{equation}
The previous equation is defined on a domain $\Omega\subset R\subseteq\RR^d$,
where $R$ is a rectangular region, and $d=1,2$ is the number of space
dimensions.

When solving \eqref{eq:poissonMixed} in two space dimensions, a significant difficulty comes from the arbitrary shape of $\Omega$. Indeed, while for simple geometries, such as a square or a circle, it is possible to find an analytical solution in terms of trigonometric and Bessel functions, for more complex geometries this is not possible and one has to resort to numerical methods.

Among several schemes developed to solve \eqref{eq:poissonMixed} on arbitrary domains, two methods of great importance are the finite element method (FEM) and the finite difference method (FDM).

The body-fitted FEM is a very general method that relies on a tessellation of the domain $\Omega$, the construction of ``good'' tessellation is far from being a trivial task, \cite{Edelsbrunner2000, Edelsbrunner2001, Henshaw1996, Schoberl1997}.
Since the generation of a mesh for $\Omega$ is such a delicate task, using descriptions that are closely related to the geometric data acquisition can dramatically reduce the cost of preprocessing an acquired geometry into representations suitable for the FEM.

For example, the generalization of standard finite elements, known as isogeometric analysis, uses a choice of basis function compatible with the description of the geometry adopted by most CAD software, \cite{Hughes2005, Sangalli2014}.
Another extremely powerful approaches, known as the CutFEM, are finite element discretizations where the boundary of $\Omega$ is represented on a background mesh, which can be of low quality, using a level set function. In the CutFEM the background mesh is also used to represent the numerical solution of \eqref{eq:poissonMixed}, \cite{Burman2015, Lehrenfeld2016, Lehrenfeld2018, Burman2022}.
Since CutFEM does not require the use of a mesh fitted to the body over which we are solving the Poisson equation, it belongs to the realm of the ``unfitted'' finite element methods.

A hybrid approach between fitted and unfitted finite elements methods is the Fictitious Domain method, first presented in \cite{Glowinski}. In this method, a body-fitted mesh is only constructed for boundary $\partial\Omega$ and a background mesh of the computational domain is used to represent the numerical solution.
In the Fictitious Domain method, the boundary conditions are imposed via distributed Lagrange multipliers, similar to \cite{GastaldiBoffi}  and \cite{GastaldiBoffiHelati}.
A comparative study between CutFEM and the Lagrange multiplier method previously discussed can be found in \cite{Helatai}.

The CutFEM is just one of the methods that makes use of level sets to describe the geometry of the domain. Level set functions in general allow us to describe the geometry of the domain in a very flexible way, \cite{Osher2002, Sussman1994, Russo2000}.

Not only the FEM can be used in combination with level set techniques to describe the geometry of the domain, but there is also a vast literature of unfitted finite difference schemes \cite{Gibou2002,Gibou2007,Shortley1938,Coco2013,Coco2018,semplice2016adaptive}, just to mention a few.
A common feature of most unfitted FDMs previously mentioned is the addition of ghost points, i.e. points that are not part of the domain $\Omega$ but are used to impose the boundary conditions and retrieve the optimal order of convergence. We can find many applications in \cite{COCO2020109623,astuto2023multiscale,ASTUTO2023111880,astuto2023time,astuto2024high}.

While the unfitted FDM, compared to unfitted FEM, might be easier to implement and understand for a novice user, these finite difference schemes have the disadvantage of preserving less structural properties of the continuous problem, such as the symmetry and the positive definiteness of the stiffness matrix.
In the present paper, we aim to introduce ideas from the unfitted FEM and FDM literature and introduce the reader to a nodal ghost finite element method for elliptic problems on irregular domains with mixed boundary conditions. An explicit comparison is shown in \cite{astuto2024comparison}.

The choice of bilinear finite elements on a structured mesh yields a similar stencil to the one typical of FDM, allowing a much more natural and straightforward implementation for a novice user.

To achieve this goal, the paper is structured as follows: we begin deriving the variational formulation of \eqref{eq:poissonMixed} in the one-dimensional case and describe our nodal symmetric ghost point method for this simple case. We then proceed to extend our nodal symmetric ghost point method to the two-dimensional case. 
In Section~\ref{sec:comp}, we provide a detailed description of the implementation aspects of the scheme we propose. In particular, we discuss the exact quadrature of the boundary integral, which is one of the main novelties of our work. In Section ~\ref{sec:analysis},  we prove a priori convergence rate for the proposed scheme. Furthermore, to the best of our knowledge, we give for the first time a rigorous treatment of the snapping back to grid perturbation. ~Lastly, we validate the method and the a priori estimates proved via extensive numerical experiments in both one and two space dimensions.

We hope that our treatment of the proposed nodal symmetric ghost point method, which starts from the most basic concepts and builds up to a more complicated scenario, will facilitate the novice reader's understating of our scheme.

\section{The continuous one-dimensional problem}
\label{sec:one_dim}
In this section, we will introduce the reader to the variational formulation of the Poisson equation with Dirichlet and mixed boundary conditions in one dimension.
While we know these ideas are well known to the reader familiar with the finite element literature, we decided to present them anyway to make the paper as self-contained as possible for the novice reader.

Let us start our discussion from the Poisson equation with homogeneous Dirichlet conditions at the boundary defined
in a domain $\Omega = [a,b]\subset R = [0,1]$: 
\begin{equation}
	-u''(x) = f(x), \quad x\in[a,b], \quad u(a) = 0, \> u(b) = 0.
	\label{eq1D0}
\end{equation}
Multiplying the above equation by a test function ${v\in V =  H^1_0([a,b])}$, provides the following variational
formulation of the problem: 
\begin{equation}
        \int_a^b \left[u'(x) v'(x) - f(x)v(x)\right] \, dx  = u'(b)v(b)-u'(a)v(a), \quad \forall v\in H^1_0([a,b])
        \label{eq:weak1}
\end{equation}
where $H^1_0([a,b])$ is the standard Sobolev space with homogeneous trace at the boundary \cite{Brezis}.
Thanks to this condition, the right-hand side of the equation vanishes,
leading to the formulation: 
\begin{equation}
        \int_a^b \left[u'(x) v'(x) - f(x)v(x)\right] \, dx  = 0,  \quad \forall v\in H^1_0([a,b])
	\label{eq:D0}
\end{equation}

The variational formulation can be obtained by setting to zero the first variation of the functional 
\begin{equation}
	J_{D0}[u] = \int_a^b \left[\frac12 (u'(x))^2\, dx - f(x) u(x)\right]\, dx,
	\label{JD0[u]}
\end{equation}
as can be immediately verified by performing the first variation of the functional by replacing 
$u$ by $u+\varepsilon v$, with $u$ and $v$ in $H^1_0([a,b])$, and imposing
\begin{equation}
	\left. \frac{d}{d\varepsilon} J_{D0}[u+\varepsilon v]\right|_{\varepsilon = 0} = 0 \quad \forall v\in
H^1_0([a,b]).
\end{equation}
In particular, we notice that since $J_{D0}$ is a convex functional, looking for functions for which the first variation of $J_{D0}$ vanishes is equivalent to minimizing the functional.

A different way to approximately impose zero Dirichlet conditions on the boundary is to penalise any
discrepancy from zero. This is obtained by adding a penalisation term to the functional:
\begin{equation}
	J_{D0\lambda}[u] = \int_a^b \left[\frac12 (u'(x))^2\, dx - f(x) u(x)\right]\, dx 
+ \frac{\lambda}2[u(a)^2+u(b)^2].
	\label{JD0L[u]}
\end{equation}
{where $\lambda$ is a positive penalty parameter \cite{Burman2012} that will be discussed in Section \ref{sec:analysis}, and it is chosen proportional to a power of the space step.}

Notice that the energy functional $J_{D0\lambda}$ corresponds to \eqref{eq1D0} where the Dirichlet boundary conditions are replaced by Robin boundary conditions.
It is a well-established fact that as the penalisation grows to infinity it will force zero Dirichlet conditions as
$\lambda\to\infty$ \cite{Courant,HilbertCourant}.
When integrating by parts, we obtain the following variational formulation of the problem: 
\begin{align}
            \label{eq:AubinBabuska}
            \int_a^b \left[u'(x)v'(x) - f(x)v(x)\right] \, dx  & + [\lambda u(a)+u'(a)]v(a) \\
    						& +  [\lambda u(b)-u'(b)]v(b) = 0 , \quad \forall v\in H^1([a,b])\nonumber
\end{align}
where now we do not need to impose that $v$ is zero on the boundary.
Notice that the variational formulation \eqref{eq:AubinBabuska} is the starting point for the penalty finite element proposed in \cite{Aubin, ConfPenaltyBabuska} and analyzed in \cite{Li, PenaltyElliott}.

If we want to impose non-homogeneous Dirichlet conditions, say we impose $u(a) = u_a$, $u(b) = u_b$, we have to consider the functional 
\begin{equation}
	J_{D\lambda}[u] = \int_a^b \left[\frac12 (u'(x))^2 - f(x) u(x)\right]\, dx + 
	\frac{\lambda}2[(u(a)-u_a)^2+(u(b)-u_b)^2]
	\label{JDL[u]}
\end{equation}
which leads to the following variational formulation:
\begin{equation}
\begin{split}
        \int_a^b \left[u'(x)v'(x) \right. & 
			\left. - f(x)v(x)\right] \, dx   + [\lambda (u(a)-u_a) + u'(a)] v(a) \\
				& +  [\lambda (u(b)-u_b) - u'(b)] v(b) = 0 , \quad \forall v\in H^1([a,b])
	\label{eq:weakDL}
\end{split}
\end{equation}

The variational formulation \eqref{eq:weakDL} is non-symmetric, because of the term $\displaystyle \left.(u'v)\right|_a^b$.
We will now address this issue by introducing an alternative formulation, which is equivalent at a continuous level and is obtained as follows. 
We add and subtract the quantity $\displaystyle \left.(uv')\right|_a^b$:
\begin{equation}
        \int_a^b u'(x) v'(x) \, dx - \left. (u'v) \right|_a^b - \left. (uv') \right|_a^b = \int_a^b f(x)v(x) \, dx - u_b v'(b) + u_a v'(a) 
        \label{eq:weak_Nitsche}
\end{equation}
where we made use of the Dirichlet boundary conditions in $a$ and $b$ on the right-hand side.

The coercivity of the operator may be improved by adding a penalisation term, therefore for continuous formulation can be written as 
\begin{align}
        \int_a^b u'(x) v'(x) \, dx - \left. (u'v) \right|_a^b - \left. (uv') \right|_a^b + &\lambda(u(a)-u_a)v(a)
       + \lambda(u(b)-u_b)v(b) \nonumber \\  = \int_a^b f(x)v(x) \, dx - u_b v'(b) + u_a v'(a) 
        \label{eq:weak3}
\end{align}
and bringing to the right-hand side the known term, we have
\begin{align}
        &\int_a^b u'(x) v'(x) \, dx - \left. (u'v) \right|_a^b - \left. (uv') \right|_a^b + \lambda u(a)v(a) 
       + \lambda u(b)v(b) \nonumber \\  =& \int_a^b f(x)v(x) \, dx - u_b v'(b) + u_a v'(a) +\lambda u_av(a) + \lambda u_b v(b)
        \label{eq:weak4}
\end{align}
It is trivial to obtain the analogous variational formulation for the Neumann problem, in fact, if we consider
\begin{equation}
    -u''(x) = f(x), \quad x\in[a,b]\qquad u'(a) = g_a, \quad u'(b) = g_b \label{eq1N},
\end{equation}
we can multiply by a test function $v\in H^1[a,b]$, integrate over $[a,b]$ and integrate by parts to see that the solution of \eqref{eq1N}$ is $ $u\in H^1[a,b]$ such that
\begin{equation}
	\int_a^b u'(x)v'(x)\, dx = \int_a^bf(x)v(x)\,dx + g_bv(b) - g_bav(a)  \quad \forall v\in H^1([a,b]).
	\label{eq:weak3N}
\end{equation}
The only aspects we need to be careful when considering a variational formulation for the Neumann problem is that a compatibility condition has to be verified by $f\in L^2([a,b])$ and  we need to choose a unique solution in the kernel of the gradient, imposing a null-average condition, i.e.
\begin{equation}
\nonumber
    \int_a^b f \,dx = g_b-g_a,\qquad \int_a^b u \,dx= 0.
\end{equation}
Let us  now consider the following mixed problem: 
\begin{equation} 
	-u''(x) = f(x), \quad x\in[a,b]
	\qquad
	u(a) = 0, \quad u'(b) = g_b \label{eq1DM_d}
\end{equation}


We shall look for solutions of
\eqref{eq1DM_d} in $H^1([a,b])$, respecting the boundary conditions.
Multiplying by a test function
\begin{equation} \nonumber
    H^1_{a}([a,b]) := \left\{v\in H^1([a,b]) \;:\; v(a) = 0 \right\},
\end{equation}
and integrating by parts in $[a,b]$ one
obtains:
\begin{equation} \nonumber
    \int_a^b u'(x)v'(x)\, dx - \left. (u'v)\right|_a^b = \int_a^bf(x)v(x)\,dx.
\end{equation}
Now, making use of the boundary condition in $b$ and the fact that for any $v\in H^1_a([a,b])$ $v(a)=0$, we obtain
the following variational formulation of problem \eqref{eq1DM_d}:
\begin{equation}
	\int_a^b u'(x)v'(x)\, dx = \int_a^bf(x)v(x)\,dx + g_bv(b) \quad \forall v\in H^1_a([a,b]).
	\label{eq:weak3b}
\end{equation}

This equation can be obtained by looking for the function $u\in H^1_a([a,b])$ that minimizes the following functional: 

\begin{equation}
	J[u] = \frac12\int_a^b [(u'(x))^2\, dx - f(x) u(x)]\, dx + g_b u(b)
	\label{JM0[u]}
\end{equation}
Indeed, by equating to zero the first variation of the functional, i.e. imposing 
\begin{equation} \nonumber
	\left. \frac{d}{d\varepsilon} J[u+\varepsilon v]\right|_{\varepsilon = 0} = 0
\end{equation}
$\forall v\in H^1_a([a,b])$, gives exactly \eqref{eq:weak3b}. 
As in the case of the Dirichlet problem, we may impose an approximate Dirichlet condition on the left hand boundary by penalising the functional as follows:
\begin{equation*}
	J[u] = \int_a^b \left[\frac12 (u'(x))^2\, dx - f(x) u(x)\right]\, dx + g_b u(b) +
\frac{\lambda}{2}(u(a)-u_a)^2 , \quad u\in H^1([a,b])
\end{equation*}
which leads to the following variational formulation of the problem in $H^1([a,b])$: 
\begin{equation}
	\int_a^b u'(x)v'(x)\, dx + 
	[\lambda(u(a)-u_a)+u'(a)]v(a)= 
	\int_a^bf(x)v(x)\,dx + g_b v(b).
	\label{eq:weakML}
\end{equation}

As before, we look for a term that symmetries \eqref{eq:weakML}. For this reason, we add and subtract the term $u(a)v'(a)$, making use of the Dirichlet boundary condition on $a$
\begin{align}
       & \int_a^b u'(x) v'(x) \, dx  + [\lambda(u(a)-u_a)+u'(a)]v(a) + u(a)v'(a)  \\ \nonumber &= \int_a^b f(x)v(x) \, dx + u_av'(a) + g_bv(b)
\end{align}
and bringing to the right-hand side the known terms, we have
\begin{align}
        \label{eq:weak_Nitsche_mixed}
       & \int_a^b u'(x) v'(x) \, dx  + \lambda u(a) v(a)+u'(a)v(a) + u(a)v'(a)  \\ \nonumber &= \int_a^b f(x)v(x) \, dx + \lambda u_a v(a) + u_av'(a) + g_bv(b) .
\end{align}
The technique we used to symmetrise the variational problem is known as the Nitsche technique \cite{Nitsche}.

\section{The discrete one-dimensional problem}
\label{sec:one_dim_discr}
The region $R = [0,1]$ is discretized by a uniform mesh of grid size $h=1/N$. 
Without loss of generality, we assume that $a\in[0,h]$, and $b\in [1-h,1]$. If $a=0$ and $b=1$, standard finite element discretization in the interval $[0,1]$ can be used.

\begin{figure}[H]
\begin{center}
\includegraphics[width=0.7\textwidth]{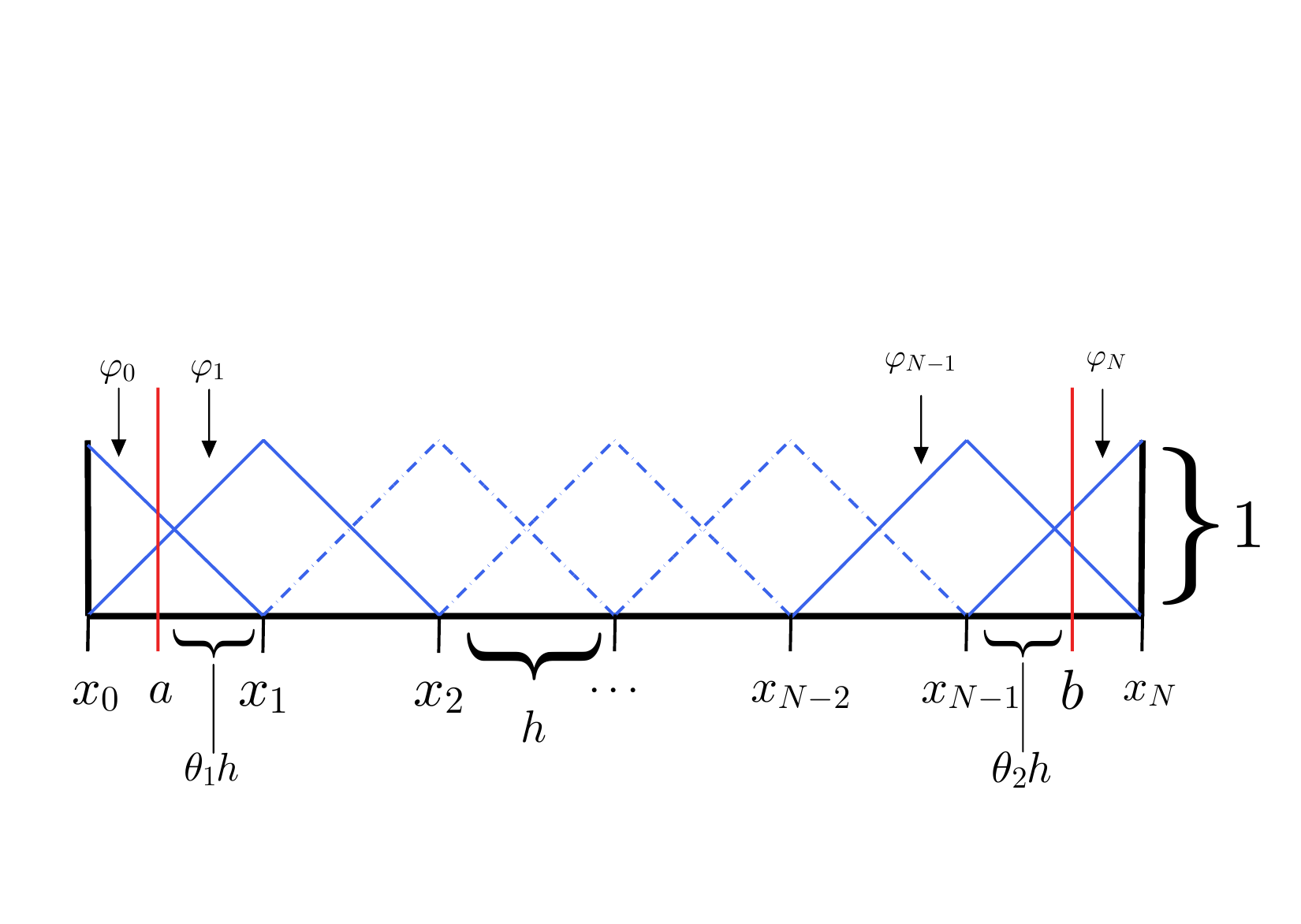}
\end{center}
\caption{\textit{Discretization setup in 1D: the region $R = [0,1]\supset [a,b] = \Omega $, and the basis function $\varphi_j, j=0,\ldots,N$, defined in \eqref{eq_basis_f}.}}
\label{fig:1D_setup}
\end{figure}

We denote by $\theta_1 = 1-a/h$ and $\theta_2 = { 1 - (1-b)/h}$, so that $\theta_1, \theta_2 \in
[0,1]$.
We denote by $I_j = [j-1,j]h$, with $j=1,\ldots,N$, the $j$-th interval of the
subdivision. 
We consider the space $V_h$ of dimension $N+1$ as the space of the splines $S_1$
on the mesh $\{I_j\}_{i=1}^N$.  
As basis for the space $V_h$ we choose the functions $\varphi_i(x)$ with $i=0,\ldots,N,$
defined as 
%
\begin{equation}
\label{eq_basis_f}
	\varphi_i(x) = \max\{1-|x-x_i|/h,0\}, \quad x\in [0,1], \quad i=0,\ldots,N.
\end{equation}
Fig.~\ref{fig:1D_setup} provides a schematic representation of the discretization setup. 

A generic function $u_h \in V_h$ can be expressed using the basis $\{\varphi_i\}_{i=1}^N$ as follows:
\begin{equation}
	u_h(x) = \sum_{j=1}^N u_j \varphi_{j}(x)
	\label{eq:u_h}
\end{equation}

We now make use of the conformity of the discrete space, i.e., $V_h\subset H^1_0([a,b])$, to obtain a discrete variational problem, by considering \eqref{eq:weakDL} with both test and trial functions in $V_h$:
\begin{equation*}
\begin{split}
       & \int_a^b u_h'(x) \varphi_i'(x) \, dx -  \int_a^b f(x)\varphi_i(x)\, dx   + \lambda
(u_h(a)-u_a)\varphi_i(a) + u_h'(a) \varphi_i(a)  \\
				& + u_h(a) \varphi_i'(a) +  \lambda (u_h(b)-u_b)\varphi_i(b) - u_h'(b)\varphi_i(b) - u_h(b)\varphi_i'(b) = 0 , \quad i=0,\ldots,N
\end{split}
\end{equation*}

By replacing $u_h$ with its representation \eqref{eq:u_h}, we obtain the following system of
equations for the unknowns $u_i$, $i=1,\ldots,N$:
\begin{align} \nonumber
   & \sum_{j = 0}^N u_j\underbrace{(\varphi_j',\varphi_i')_{L^2([a,b])}}_{\rm S} + \sum_{j = 0}^N u_j\underbrace{\left[\varphi_j'(a)\varphi_i(a) + \varphi_j(a)\varphi_i'(a) - \varphi_j'(b)\varphi_i(b) - \varphi_j(b)\varphi_i'(b) \right]}_{\rm S_T} \\ 
   & + \lambda \sum_{j = 0}^N u_j\left(\underbrace{\varphi_j(a)\varphi_i(a)}_{\rm P_a} +  \underbrace{\varphi_j(b)\varphi_i(b)}_{\rm P_b} \right) \\ \newline 
    & = \sum_{j = 0}^N f_j\underbrace{(\varphi_j,\varphi_i)_{L^2([a,b])}}_{\rm M} \, dx + \sum_{j = 0}^N {u_a} \underbrace{\varphi_j(a)\varphi_i'(a)}_{\rm D_a} + \sum_{j = 0}^N {u_b} \underbrace{\varphi_j(b)\varphi_i'(b)}_{\rm D_b} + \\ 
    & + \lambda \left(\sum_{j = 0}^N{u_a} \underbrace{\varphi_j(a) \varphi_i(a)}_{\rm P_a} + \sum_{j = 0}^N{u_b} \underbrace{\varphi_j(b) \varphi_i(b)}_{\rm P_b} \right), \quad i=0,\ldots,N 
\end{align}
and in matrix form it becomes
\begin{equation}
	\sum_{j=0}^N a^D_{ij} u_j = f^D_i, \quad i=0,\ldots,N
	\label{eq:system1D}
\end{equation}
where $A^D = \{a^D_{i,j}\}$ and $F^D = \{f^D_i\}$, and
\begin{equation}
\begin{array}{c}
	\displaystyle	A^D  = {\rm S + S_T +} \lambda\, \left({\rm P_a+ P_b}\right)\\
	\displaystyle F^D  = {\rm M}f + \left({\rm D_a} +  \lambda {\rm P_a}\right) u_a +  \left({\rm D_b} +  \lambda {\rm P_b}\right) u_b
\end{array}
\label{coeffD1}
\end{equation}
where the matrices ${\rm A^D, M, D_a, D_b, P_a, P_b} \in\RR^{(N+1)\times (N+1)}$.


We now turn our attention to the mixed Dirichlet-Neumann problem \eqref{eq1DM_d}. 
We consider the penalised variational formulation in \eqref{eq:weak_Nitsche_mixed} 
and look for a solution in $V_h$. 
Replacing $u$ by $u_h$ and $v$ by $\varphi_j$ in \eqref{eq:u_h}, we obtain the system 
\begin{align}
   & \int_a^b u_h'(x)\varphi_i'(x)\, dx + u_h'(a)\varphi_i(a) + u_h(a)\varphi_i'(a)  + \lambda u_h(a)\varphi_i(a) \label{eq:weakMLN} \\ \newline 
    & = \int_a^b f_h(x)\varphi_i(x) \, dx + {u_a} \varphi_i'(a) + {g_b}_h\varphi_i(b) + \lambda {u_a}  \varphi_i(a), \quad \forall i = 1,\cdots,N \nonumber
\end{align}
and if we use the expression in \eqref{eq:u_h} for the function $u_h$, it becomes
\begin{align} \nonumber
   & \sum_{j = 0}^N u_j\underbrace{(\varphi_j',\varphi_i')_{L^2([a,b])}}_{\rm S} + \sum_{j = 0}^N u_j\underbrace{\left[\varphi_j'(a)\varphi_i(a) + \varphi_j(a)\varphi_i'(a) \right]}_{\rm S_T}  + \lambda \sum_{j = 0}^N u_j\underbrace{\varphi_j(a)\varphi_i(a)}_{\rm P_a} \\ \nonumber \newline 
    & = \sum_{j = 0}^N f_j\underbrace{(\varphi_j,\varphi_i)_{L^2([a,b])}}_{\rm M} \, dx + \sum_{j = 0}^N {u_a} \underbrace{\varphi_j\varphi_i'(a)}_{\rm D_a} + \sum_{j = 0}^N{g_b}_j\underbrace{\varphi_j(b)\varphi_i(b)}_{\rm N_b} \\
    &+ \lambda \sum_{j = 0}^N{u_a} \underbrace{\varphi_j \varphi_i(a)}_{\rm P_a}, \quad i=0,\ldots,N 
	\label{eq:weakMLN_discrete}
\end{align}
which leads to the following system of equations for $u_i$:
\begin{equation}
	\sum_{j=0}^N a^M_{ij} u_j = f^M_i, \quad i=0,\ldots,N
	\label{eq:system1M}
\end{equation}
where $A^M = \{a^M_{i,j}\}$ and $F^M = \{f^M_i\}$, and
\begin{equation}
\begin{array}{c}
	\displaystyle	A^M  = {\rm S + S_T +} \lambda\, {\rm P_a}\\
	\displaystyle F^M  = {\rm M}f + \left({\rm D_a} +  \lambda {\rm P_a}\right) u_a + {\rm N_b }\, g_b
\end{array}
\label{coeffD1M}
\end{equation}
The non-zero values of the stiffness matrices $A^D$ and $A^M$, and of the source terms $F^D$ and $F^M$, are computed explicitly in \ref{Appendix:B}.

\section{The two-dimensional problem} 
In this section, we introduce the two-dimensional problem and define our numerical scheme based on finite-element discretization. We begin by imposing Dirichlet boundary conditions on the elliptic equation and then extend our discussion to include Neumann and mixed boundary conditions. Detailed derivations of the spatial discretization are provided in \ref{app:sec_2D}. 

Let us consider now the Poisson equation with Dirichlet boundary conditions in $\Omega\subseteq R\subset \RR^2$, where $R$ is a rectangular region:
\begin{equation}
	-\Delta u = f \quad {\rm in}\>  \Omega, \quad u = g_D \> {\rm on}\>  \partial \Omega
	\label{eq2D0}	
\end{equation}
{Defining the space of test functions 
\[ V = \{ v \in H^1(\Omega)\},
\]} multiplying Eq.\eqref{eq2D0} by $v\in V$, integrating in $\Omega$ and then integrating by parts, one obtains: 
\begin{equation}
	\int_{\Omega}\nabla u\cdot \nabla v \, d\Omega = \int_\Omega f v\, d\Omega + \int_{\partial
\Omega}v\pad{u}{n}\, ds.
\end{equation}

\begin{figure}[h]
    \centering
		\begin{overpic}[abs,width=0.7\textwidth,unit=1mm,scale=.25]{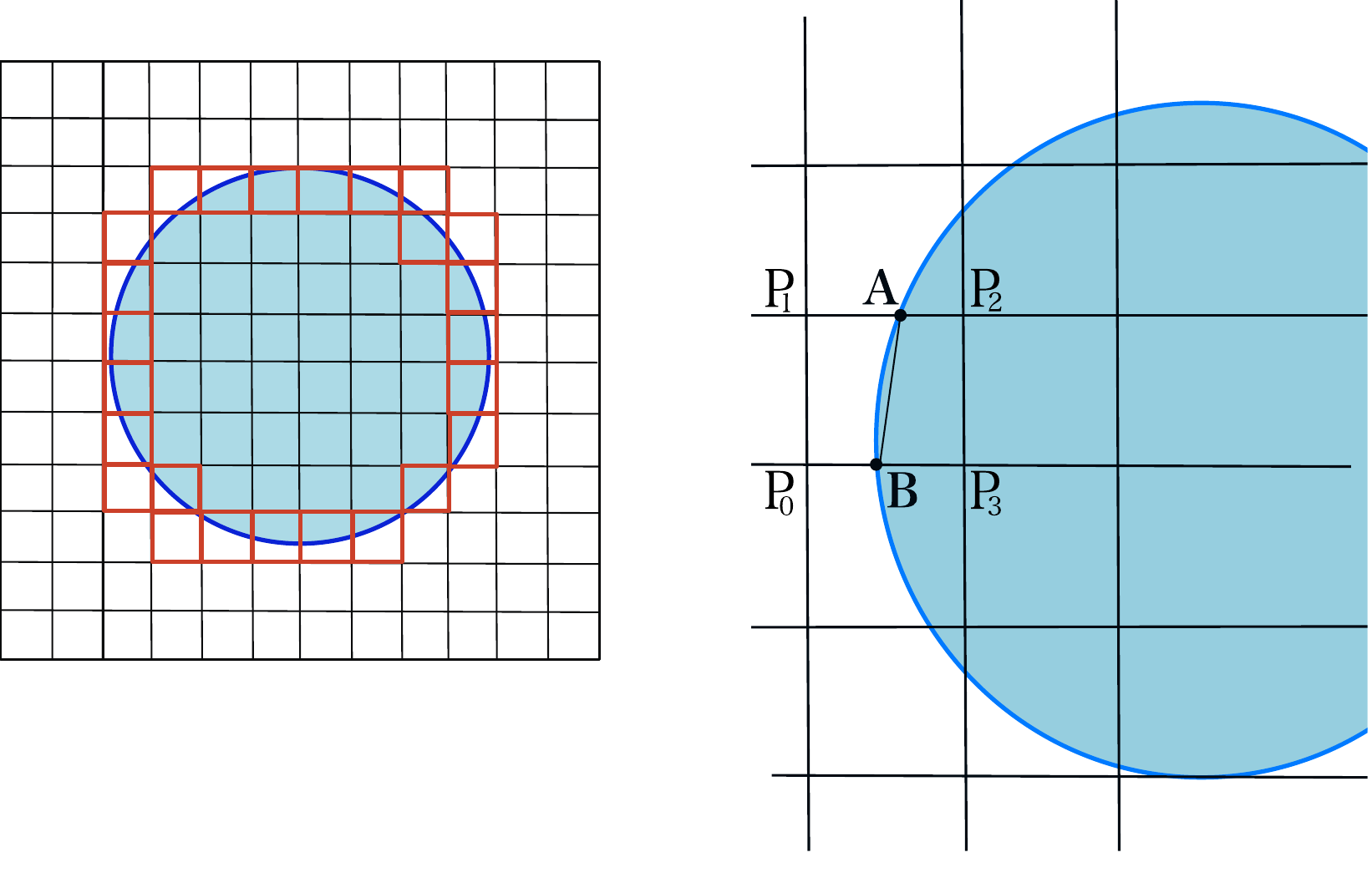}
			\put(-7,58){(a)}
                \put(-2,44){$\nearrow$}
                \put(-2,25){$\xrightarrow{\hspace*{2cm}}$}
                \put(-15,27){$\phi<0$}
                \put(-15,40){$\phi>0$}
                \put(48,58){(b)}
                \put(21.5,30.8){$\Omega$}
                \put(11,56){$R$}
                \put(65,56){$R$}
                \put(80,35){$\Omega$}
		\end{overpic}        
\caption{\textit{(a): Square regular mesh of the region R and of the domain $\Omega$. In red, we have those cells that intersect the boundary $\Gamma$. (b): Representation of the two points \textbf{A}, \textbf{B} defined in Algorithm~\ref{alg_ab}, making use of the points $\{P_0,\ldots, P_4 \}$, where $P_0$ and $P_1$ are ghost points.}} 
\label{fig:Domain2D}
\end{figure}

The discretization of the problem is obtained by discretizing the rectangular region $R\subset \RR^2$, containing the domain $\Omega$, by a regular rectangular grid (see Figure~\ref{fig:Domain2D}).
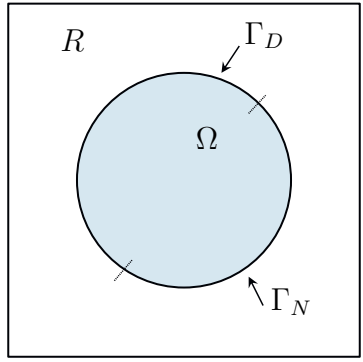
\begin{figure}
    \centering
    \begin{overpic}[abs,width=0.4\textwidth,unit=1mm,scale=.25]{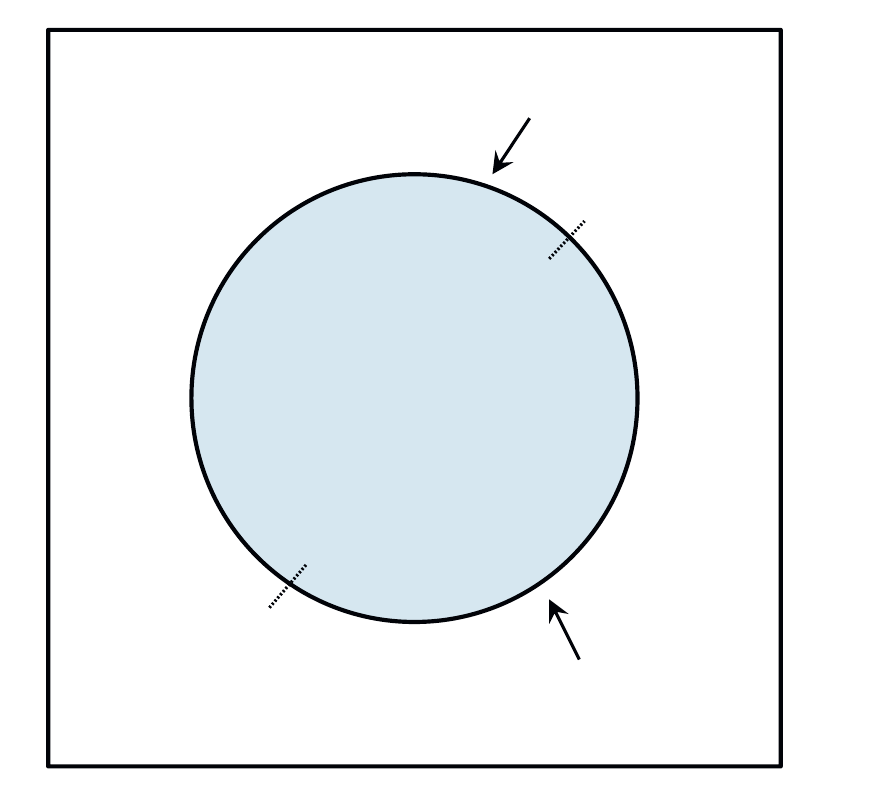}
    \put(10,43){$R$}
    \put(34.5,44){$\Gamma_D$}
    \put(38,8.5){$\Gamma_N$}
        \put(28,30){$\Omega$}
		\end{overpic}   
    \caption{\textit{Two dimensional problem domain $\Omega\subset R$: $\Gamma_D$ and $\Gamma_N$ are portions of $\Gamma = \partial \Omega$, with Dirichlet and Neumann boundary conditions, respectively. $R$ is the rectangular region where we define the square grid.}}
    \label{fig:domain_bc}
\end{figure}
Following the approach shown in  \cite{Sussman1994,Osher2002,Russo2000,book:72748}, the domain $\Omega$ is implicitly defined by a level set function $\phi(x,y)$ that is negative inside $\Omega$, positive outside and zero on the boundary $\Gamma$ :
\begin{eqnarray} \nonumber
	\Omega = \{(x,y): \phi(x,y) < 0\}, \qquad
	\Gamma = \{(x,y): \phi(x,y) = 0\}.
\end{eqnarray}
The unit normal vector $n$ in \eqref{eq2DM} can be computed as $n = \frac{\nabla \phi }{|\nabla \phi|}$.
For example, for a circular domain centered in $(x_c,y_c)$, the most convenient level-set function in terms of numerical stability is the signed distance function between $(x,y)$ and $\Gamma$, i.e.,\ $\phi=r-\sqrt{(x-x_c)^2+(y-y_c)^2}$, where $r$ is the radius of the circle.

The set of grid points will be denoted by $\mathcal N$, with $\# \mathcal N = (1+N)^2$, the active nodes (i.e.,\ internal $\mathcal{I}$ or ghost $\mathcal{G}$) by $\mathcal{I}\cup\mathcal{G} = \nodes \subset \mathcal N$, while the set of inactive points by $\mathcal O \subset \mathcal N$, with $\mathcal O\cup\mathcal A = \mathcal N$ and $\mathcal O \cap \mathcal A = \emptyset$ and the set of cells by $\mathcal C$, with $\# \mathcal C = N^2$. {Finally, we denote by $\Omega_c = R\setminus \Omega$ the outer region in $R$.}

Here, we define the set of ghost points $\mathcal{G}$, which are {grid} points that belong to $\Omega_c$, with at least an internal point as neighbor, formally defined as
\begin{equation}
\notag
	(x,y) \in \mathcal{G} \iff (x,y) \in {\mathcal N}\cap \Omega_c  \text{ and } \{(x \pm h,y),(x,y\pm h), (x \pm h,y\pm h) \} \cap \mathcal I \neq \emptyset.
\end{equation}

{Now we introduce the family of finite-elements spaces 
\begin{equation}\label{eq:Vh}
V_h = \{ v_h \in V\, : {v_h}|_K \in \mathbb Q_2(K)\quad \forall K\in \Omega_h\},
\end{equation}
where $\mathbb Q_2(K)$ denotes the space of piecewise bilinear functions on $K$,}
which are continuous in $R$.
Making use of the basis functions as test functions, and of the representation of $u_h\in V_h$, $u_h(x,y) = \sum_{i\in\mathcal{A}} u_i\varphi_i(x,y)$, we obtain the following system of equations for $u_i$:
\begin{equation}
	\sum_{j\in\nodes} a^D_{ij} u_j = f^D_i, \quad i\in\nodes
	\label{eq:system2D}
\end{equation}
where 
\begin{equation}
\begin{array}{c}
	\displaystyle	a^D_{ij}  = (\nabla\varphi_i,\nabla\varphi_j)_{L^2(\Omega_h)} + \lambda (\varphi_i,\varphi_j)_{L^2(\Gamma_{D,h})} 
	- (n\cdot\nabla\vphi_j,\vphi_i)_{L^2(\Gamma_{D,h})} \\ - (\vphi_j,n\cdot\nabla\vphi_i)_{L^2(\Gamma_{D,h})},\\[4mm]
	\displaystyle	f^D_i  = \lambda (g_D,\vphi_i)_{L^2(\Gamma_{D,h})}  - (g_D,n\cdot\nabla\vphi_i)_{L^2(\Gamma_{D,h})}  + (f,\vphi_i)_{L^2(\Omega_h)}.
\end{array}
\label{coeffD2}
\end{equation}

As discussed in the one-dimensional case, it is easy to obtain a variational formulation for Neumann boundary conditions
\begin{equation}
	-\Delta u = f \quad {\rm in}\>  \Omega, \qquad \pad{u}{n} = g_N \> {\rm on}\>  \Gamma_N,
	\label{eq2N}	
\end{equation}
where $\Gamma_N$ is the all perimeter of $\Omega$, i.e $\Gamma_N=\partial\Omega$. Notice that for purely Neumann problems there is no penalisation term, since $\Gamma_N = \Gamma$ and $\Gamma_D = \emptyset$.
For such problems, we assume the standard compatibility condition is satisfied, i.e.
\begin{equation} \nonumber
    \int_\Omega f\, d\Omega + \int_\Gamma g_N\, d\Gamma = 0.
\end{equation}
Making use of the basis functions as test functions in \eqref{eq2N}, we obtain the following system of equations for $u_i,$ $i \in \mathcal{A}$:
\begin{equation}
	\sum_{j\in\nodes} a^N_{ij} u_j = f^N_i, \quad i\in\nodes
	\label{eq:system2N}
\end{equation}
where $A^N = \{a^N_{ij} \}$ and $F^N = \{f_i \}$, and
\begin{align}
\begin{split}
\label{coeffN}
A^N = {\rm S},\\
F^N = {\rm M} f + {\rm N_{\Gamma}}\, g_N,
\end{split}
\end{align}
and the entries of the matrices are 
\begin{align}
\begin{split}
    \label{eq_coeff_aN}
	\displaystyle	a^N_{ij}  & = (\nabla\varphi_i,\nabla\varphi_j)_{L^2(\Omega_h)},\\ 
	\displaystyle	f^N_i & = (f,\vphi_i)_{L^2(\Omega_h)} + (g_N,\vphi_i)_{L^2(\Gamma_{N,h})}.
 \end{split}
\end{align}

We now consider the mixed problem in two dimensions:
\begin{equation}
	-\Delta u = f \quad {\rm in}\>  \Omega, \quad u = g_D \> {\rm on}\>  \Gamma_D, \> 
	\pad{u}{n} = g_N \> {\rm on}\>  \Gamma_N
	\label{eq2DM}	
\end{equation}
where the boundary $\partial\Omega$ is partitioned into a Dirichlet boundary $\Gamma_D$ and 
a Neumann boundary $\Gamma_N$, $\Gamma_D\cup \Gamma_N = \partial \Omega$, $\Gamma_D\cap\Gamma_N=\emptyset$.

This results in the discrete system 
\begin{equation}
	\sum_{j\in\nodes} a^M_{ij} u_j = f^M_i, \quad i\in\nodes
	\label{eq:system2M}
\end{equation}
where $A^M = \{a^M_{ij} \}$ and $F^M = \{f_i \}$, and
\begin{align}
\begin{split}
\label{coeffM2}
A^M = {\rm S - S_T} +\lambda {\rm P_{\Gamma_D}},\\ 
F^M = {\rm M} f + (\lambda{\rm P_{\Gamma_D}- D_{\Gamma_D}}) g_D + {\rm N_{\Gamma_N}}\, g_N,
\end{split}
\end{align}
and the entries of the matrices are 
\begin{equation}
    \begin{array}{c}
    \label{eq_coeff_aM}
	\displaystyle	a^M_{ij}  = (\nabla\varphi_i,\nabla\varphi_j)_{L^2(\Omega_h)} 
	+ \lambda (\varphi_i,\varphi_j)_{L^2(\Gamma_{D,h})} 
	- (n\cdot\nabla\vphi_j,\vphi_i)_{L^2(\Gamma_{D,h})} \\ - (\vphi_j,n\cdot\nabla\vphi_i)_{L^2(\Gamma_{D,h})},\\[4mm]
	\displaystyle	f^M_i  = \lambda (g_D,\varphi_i)_{L^2(\Gamma_{D,h})} - (g_D,n\cdot\nabla\vphi_i)_{L^2(\Gamma_{D,h})}
	+ (g_N,\vphi_i)_{L^2(\Gamma_{N,h})} \\ + (f,\vphi_i)_{L^2(\Omega_h)}.
\end{array}
\end{equation}
\section{Computational aspects and snapping back to grid}
\label{sec:comp}
As previously discussed we are interested in dealing with the discrete variational formulations in \eqref{eq:system2D},\eqref{eq:system2N} and \eqref{eq:system2M}, and for this reason, in this section, we will detail the implementation aspects regarding how to compute the entries of those linear systems.

\begin{figure}[h]
    \centering
    \begin{minipage}{.49\textwidth}    
		\begin{overpic}[abs,width=0.8\textwidth,unit=1mm,scale=.25]{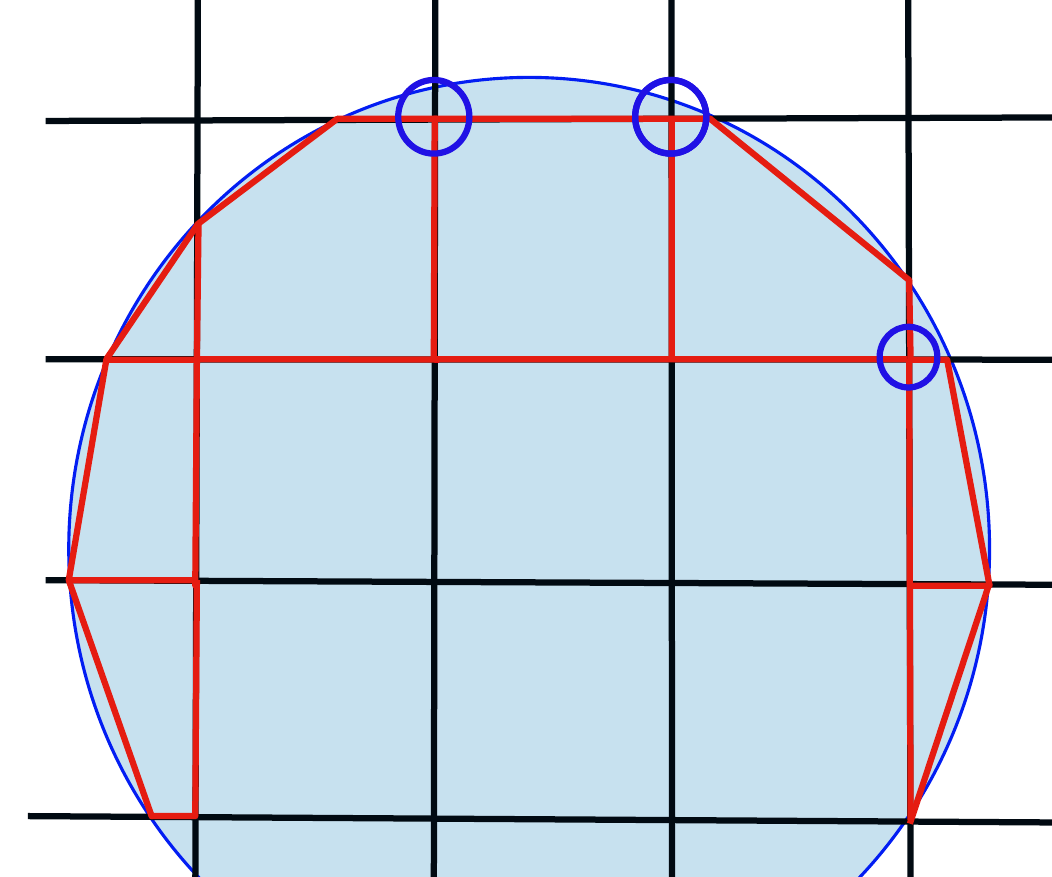}
			\put(-5,48){(a)}
    \put(27,47){convex}
   \put(48,32){non--convex}
                \put(28,20){$\Omega$}
		\end{overpic}    
    \end{minipage}
    \begin{minipage}{.49\textwidth}    
		\begin{overpic}[abs,width=\textwidth,unit=1mm,scale=.25]{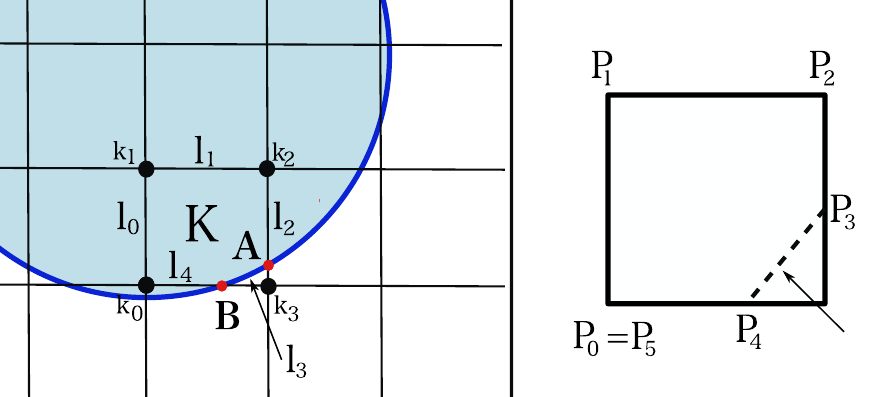}
			\put(-5,40){(b)}
                \put(25,28){$\Omega$}
                \put(52,16){$\widehat K$}
                \put(44,13){$\widehat l_0$}
                \put(53,25){$\widehat l_1$}
                \put(53,8.5){$\widehat l_4$}
                \put(65.5,17){$\widehat l_2$}
                \put(66.5,2){$\widehat l_3$}
		\end{overpic}    
    \end{minipage}
\caption{\textit{(a): Possible shapes of those polygons that intersect the boundary $\Gamma$. These elements can assume different shapes, from triangles to quadrilateral and pentagons. In red we show the set of cells $\mathcal C$  after the snapping back to grid.  (b):  Representation of the indices $\{k_0, \ldots, k_{m-1}\}$, the edges $\{l_0,\ldots,l_{m-1} \}$ of a generic element $K$, and the intersections with $\Gamma$: $\bf A$ and $\bf B$; on the right, a zoom-in that represents the reference element $\widehat K$, with edges $\{\widehat l_0,\ldots,\widehat l_{m-1} \}$, and vertices $\{P_1,\ldots, P_{m+1} \}$. }}  
\label{figure_class_cell_generic_cell}
\end{figure}

In particular, we will adopt a reference element type of approach that is standard in the finite element community. Hence, we consider a structured tessellation $\mathcal{T}_h$ (see Fig.~\ref{figure_class_cell_generic_cell} (a)) of the computational domain $R\subset \mathbb{R}^2$ and map each element of the tessellation to reference unit square $\widehat{K}$ where \eqref{coeffD2},\eqref{eq_coeff_aN} and \eqref{eq_coeff_aM} are evaluated, \cite{BrennerScott, Ciarlet}. We will denote $F_K:K\to\widehat{K}$ the map from the physical element $K\in \mathcal{T}_h$ to the reference element $\widehat{K}$.

If $i$ is an inactive point, i.e., $i\not\in \mathcal{A}$, we shall set $f^\bullet_i=0$, and the components of matrix $A^\bullet$ corresponding to such a point are $a^\bullet_{ij}=a^\bullet_{ji} = \delta_{ij}$. In this way, matrix $A^\bullet$ will be non-singular and the solution will be zero on the inactive points. We adopted the notation $f^\bullet,a^\bullet_{ij},A^{\bullet}$ to denote either $f^D,a^D_{ij},A^D$, $f^N,a^N_{ij},A^N$ or $f^M,a^M_{ij},A^M $ as in \eqref{coeffD1M},\eqref{eq_coeff_aN} and \eqref{coeffM2}.

We consider now the various integrals appearing in \eqref{coeffD2},\eqref{eq_coeff_aN} and \eqref{eq_coeff_aM}. Let us start from the integral $\left( \varphi_i , \varphi_j\right)_{L^2(K)}$. Since $\varphi_{{j}}({x}_{i},y_i)=\delta_{{i}{j}}, \, i,j \in \mathcal N$,
it is possible to reduce the summations to the vertices of the $K-$th cell (see Fig.~\ref{figure_class_cell_generic_cell} (b)), yielding 
\begin{align}
\label{eq_int_discr2}
      \left(\varphi_i, \varphi_j\right)_{L^2(K)} = \sum_{\eta,\mu = 0}^{m-1}\left(\varphi_{k_\eta},\varphi_{k_\mu} \right)_{L^2(K)},
\end{align}
where $m$ is the number of edges of the $K-$th cell, and it is $m = 4$ for the internal cells (those cells whose vertices belong to $\Omega_h$), but it can also be $m = 3$ or $m=5$ when a cell intersect the boundary (see, for example, Fig.~\ref{figure_class_cell_generic_cell} (a)). We denote by capital letters the cell indices, i.e., $K \in \mathcal{C}$, and by small letters the vertices of the cells. 
\begin{algorithm}
\caption{Computation of the intersection of the boundary with the grid}\label{alg_ab}
\begin{algorithmic}
\State $k_4 = k_0$
\For{i = 0:3}
\If{$\phi(k_i)\phi(k_{i+1})<0$} 
   \State $\theta = \phi(k_i)/(\phi(k_i)-\phi(k_{i+1})) $ 
   \State $P = \theta k_{i+1}+(1-\theta)k_i$
   \If{$\phi(k_i)<0 $}
       \State ${\bf A}:=P$
       \Else
       \State ${\bf B}:=P$
   \EndIf
\EndIf 
\EndFor
\end{algorithmic}
\end{algorithm}
Two fundamental novelties of the scheme we propose are the following:
first, the unknowns are always located on the original square grid, more precisely, on interior or ghost nodes.
Second, the computation of the integrals that define all matrix elements, characterizing systems \eqref{eq:system2D},\eqref{eq:system2N} and \eqref{eq:system2M}, are computed by exact quadrature formulas defined on the boundary of each cell, as explained below.

We observe that the product of two elements in $V_h$ appearing in \eqref{eq_int_discr2}, is an element of $\mathbb{Q}_2(K)$, i.e., the set of bi-quadratic polynomials in $K$. 

\begin{figure}
    \centering
   \includegraphics[width=0.4\textwidth]{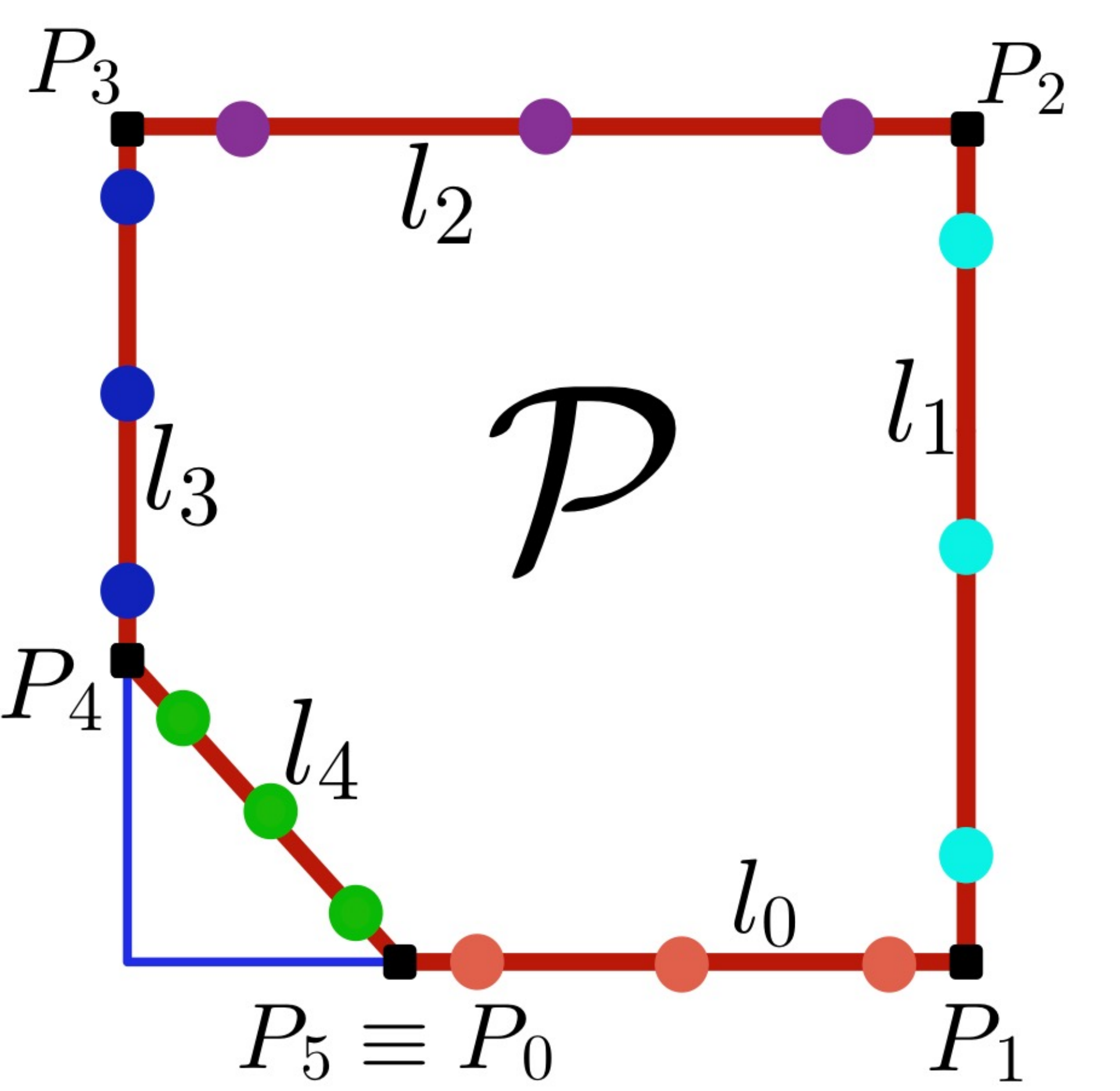}
    \caption{\textit{Scheme of the three quadrature points (circles) for each edge $l_i, \, i = 0,\cdots,4$. The squared points represent the vertices $P_i,\, i = 0,\cdots,4,$ of the polygon $\mathcal P$.}}
    \label{fig:quadr_cells}
\end{figure}

Let us consider a general integrable function $f$ defined in $\Omega$, with $F(y) = \int f(x,y) dx$ (in our strategy we consider the primitive in $x$ direction; analogue results can be obtained integrating in $y$ direction). We now define a vector $\textbf{F}= (F,0)^\top$ such that
$\nabla \textbf{F} = f$.
Thus, we have
\begin{align}
\label{eq_diverg}
    \int_K f \,dx dy = \int_K \nabla \cdot \vec F \,dx dy.
\end{align}
We can apply the Gauss theorem to \eqref{eq_diverg}, to obtain
\begin{align}
\label{eq_gauss_th}
   \int_K \nabla \cdot \vec F\, dxdy = \int_{\partial K} \vec F\cdot \vec{n}\, dl,
\end{align}
where $\vec{n} = (n_x,n_y)$ is the unit normal vector to $K$.
If $\mathcal{P}$ is the generic polygon with $m$ edges $l_r, \, r = 0,\ldots,m-1$ (see Fig.~\ref{figure_class_cell_generic_cell} (b) and Fig.~\ref{fig:quadr_cells}) 
 and express \eqref{eq_gauss_th} as
\begin{align}
\label{eq_gauss_interval}
   \int_\mathcal{P} \nabla \cdot \vec F\, d\mathcal P = \int_{\partial \mathcal{P}} \vec F\cdot \vec{n}\, d l = \sum_{r=0}^{m-1} \int_{l_r} \vec F\cdot \vec{n}\, dl = \sum_{r=0}^{m-1} \int_{l_r}  F {n}_x\, dl = \sum_{r=0}^{m-1} \int_{l_r}  F dy. 
\end{align}
To evaluate the integral over the generic edge $l_r$ we choose the three-point Gauss-Legendre quadrature rule, which is exact for polynomials in $\mathbb P_5(\mathbb R)$. Therefore, we write
\begin{equation}
     \int_{l_r}  F dy = \sum_{s = 1}^3w_sF(\widehat x_{r,s},\widehat y_{r,s}) (y_{P_{r+1}}-y_{P_r})
\end{equation}
(see Fig.~\ref{figure_class_cell_generic_cell} (b)) where $w_s \in \{5/18,4/9,5/18\}$, $\widehat x_{r,s} \in \{1-\sqrt{3/5},1,1+\sqrt{3/5} \}/2 \subset [0,1]$, with $s = 1,2,3$.

Choosing $f=\varphi_{k_\eta}\varphi_{k_\mu}\in \mathbb{Q}_2(K)\subset \mathbb P_4(K)$, and making use of \eqref{eq_gauss_interval}, we write
\begin{equation}
     \left(\varphi_i, \varphi_j\right)_{L^2(K)} = \sum_{r=0}^{m-1} \left(\sum_{s=1}^3 w_s\Phi_{i,j}(\widehat x_{r,s},\widehat y_{r,s}) (y_{P_{r+1}}-y_{P_r})\right)|l_r|.
\end{equation}
where $\Phi_{i,j} = \int \varphi_i \varphi_j\, dx$, and the formula is exact because $\Phi_{ij}, \in \mathbb P_5(K), \, \forall i,j \in \mathcal N$. Analogously, the integral for the stiff term in \eqref{coeffD2}, \eqref{eq_coeff_aM} and \eqref{eq_coeff_aM} is computed as follows
\begin{equation}
     \left(\nabla \varphi_i, \nabla \varphi_j\right)_{L^2(K)} = \sum_{r=0}^{m-1} \left(\sum_{s=1}^3 w_s\Psi_{i,j}(\widehat x_{r,s},\widehat y_{r,s}) (y_{P_{r+1}}-y_{P_r})\right)|l_r|.
\end{equation}
where $\Psi_{i,j} = \int \nabla\varphi_i\nabla\varphi_j \, dx$.

Stability analysis plays an important role in the analysis of numerical schemes when solving a partial differential equation in arbitrary domains \cite{bulenda2001stability}. 
In Fig.~\ref{figure_class_cell_generic_cell}, we show different types of stability failure, where we modify the discrete domain $\Omega_h$. When the distance between an internal point and the boundary is proportional to a power of the space step, i.e., $\propto h^\alpha$, we neglect that cell. These cases are {non} trivial because the computation of the integrals in \eqref{coeffD2} and \eqref{eq_coeff_aM} could increase the conditioning number of the computational matrices $A^D$ or $A^M$. The {snapping back to grid} technique, that we show in Algorithm~\ref{alg_snap} (see Fig.~\ref{figure_class_cell_generic_cell} (a)) is a way to avoid an ill-conditioned matrix while maintaining the order of accuracy of the method. 

Various techniques, both physical and mathematical, were employed to discover the optimal shape of the computational domain that improves the stability of the numerical schemes \cite{manuello2016step}. In Algorithm~\ref{alg_snap} we identify internal cells located in proximity to the boundary $\Gamma$. Evaluating the level set function $\phi$ at the vertices of each cell, if the value is less than $  h^\alpha$, for a chosen $ \alpha$, we disregard the respective cell. In all our numerical results we set {\tt eps} in the algorithm to Matlab machine epsilon.

In Fig.~\ref{figure_class_cell_generic_cell} (a), we present two distinct scenarios starting from a convex domain: one, where the remaining computational domain retains its convex nature, and a second scenario where convexity is lost due to the application of the snapping back to grid technique. For these two different cases, a different analysis will be conducted in Section~\ref{sec:analysis}.

\begin{algorithm}[H]
\caption{Snapping back to grid}\label{alg_snap}
\begin{algorithmic}
\For{$k \in \mathcal N$} 
\If{$\phi(k)<0 \quad \& \quad |\phi(k)|<h^{ \alpha}$} 
   \State $\phi(k) := {\tt eps} $ 
\EndIf 
\EndFor
\end{algorithmic}
\end{algorithm}

\section{A priori convergence analysis}
\label{sec:analysis}
In this section, we will focus our attention on developing a priori error estimates for the scheme described in the paper. 
{ 
We will begin the analysis of the method under the following geometrical constraint on the boundary $\partial \Omega$:
\begin{equation}
    \forall K \in \{\mathcal{T}_h\}_{h>0}\quad \widehat{C}_\Omega h^\alpha \leq \abs{\partial \Omega \cap K} \leq \widehat{C}_\Omega h^\gamma,
\end{equation}
with $\gamma < \alpha$. In principle $\gamma$ could be arbitrarily small or even negative and $\alpha$ could be arbitrarily large.
Yet we are interested in the case when $\partial\Omega\cap K$ is comparable to the edges of the cell, which translates to $\gamma = 1$.
Furthermore, we are interested in allowing $\partial\Omega\cap K$ to have smaller length than the edges of the cell, but not too small, i.e. $\alpha = 2$.
While general values of $\alpha$ and $\gamma$ are considered in the analysis, our work has as main focus the case when $\alpha \sim 2$ and $\gamma \sim 1$.
}
The analysis here presented will follow closely \cite{Lehrenfeld2018}.  
To lighten the notation we introduce the discrete bilinear form $a_{h,\alpha}:H^{\frac{3}{2}+\varepsilon}(\Omega)\times H^{\frac{3}{2}+\varepsilon}(\Omega) \to \mathbb{R}$ corresponding to \eqref{eq:weak2M_Nitsche_discrete}, i.e., 
\begin{equation}
    \label{eq:ah}
    a_{h,\alpha}(u,v) := (\nabla u,\nabla v)_{L^2(\Omega)} + C_a h^{-\alpha} (u,v)_{L^2(\Gamma_D)}-(\partial_nu,v)_{L^2(\Gamma_D)}-(u,\partial_n v)_{L^2(\Gamma_D)}.
\end{equation}

To study the well-posedness of the discrete variational problem: find $u_h\in V_h$ such that
\begin{equation}
    \label{eq:discreteVarProblem}
    a_{h,\alpha}(u_h,v_h) := (f,v_h)_{L^2(\Omega)}+(g_N,v_h)_{L^2(\Gamma_N)}-(g_D,\partial_n v_h)_{L^2(\Gamma_D)}, \qquad v_h \in V_h
\end{equation}
we will consider an auxiliary eigenvalue problem, i.e., find $(w_h,\lambda)\in V_h\times \mathbb{R}$ such that 
\begin{equation}
    b(w_h,z_h) := (\partial_n w_h,\partial z_h)_{L^2(\partial \Omega)} = \lambda a(w_h,z_h), \qquad \forall z_h\in V_h.
    \label{eq:eig}
\end{equation}
We characterize the first eigenvalue problem using the Rayleigh quotient we have,
\begin{equation}
    \lambda= \underset{z_h\in V_h}{\max} \frac{b(z_h,z_h)}{a(z_h,z_h)} = \underset{z_h\in V_h}{\max} \norm{\partial_n z_h}_{L^2(\partial \Omega)}^2\norm{\nabla z_h }^{-2}_{L^2(\partial \Omega)}.
\end{equation}
The above characterization ensures that for any $v_h\in V_h$ the following estimate holds
\begin{equation}
    \norm{\partial_n v_h}^2_{L^2(\Omega)}\leq \lambda\norm{\nabla v_h}_{L^2(\Omega)}^2.
\end{equation}
We will now bound $\lambda$ in terms of the mesh size, this will be useful to prove the coercivity of the bilinear form $a_h(\cdot,\cdot)$ for a sufficiently large value of $\alpha$.
The following lemma is an adaptation of \cite[Lemma 12.8]{ErnGuermond1}.
\begin{lemma}
Let $V_h$ be the space of bilinear functions on each element $K$ of $\mathcal{T}_h$, then
\begin{equation}
    \label{eq:discreteTraceInequality}
    \norm{\partial_n u_h}^2_{L^2(\partial\Omega\cap K)}\leq C_2 h^{\gamma-{ \alpha}-1}\norm{\nabla u_h}_{L^2(K)}^2,  \qquad \forall v_h \in V_h,
\end{equation}
provided that the following geometric assumption is verified by $\mathcal{T}_h$:
\begin{equation}
    \label{eq:GeoAssumption1}
    \widehat{C}_\Omega h^{ \alpha} \leq \abs{\partial\Omega\cap K}\leq C_\Omega h^{\gamma}, \quad { \forall K\in \mathcal{T}_h.}
\end{equation}
\end{lemma}
\begin{proof}
    We will obtain this estimate via a scaling argument, for this reason, we consider the map from the physical element $K$ to the reference element $\widehat{K}$, i.e., $F_K: K\to \widehat{K}$.
    Notice that since we are working with a structured mesh the map $F_K: K \to \widehat{K}$ is an affine map, i.e., $F_K(\vec{x})=\mathbb{A}\vec{x}+\vec{b}$.
    We then consider the pullback of $v_h$, i.e., 
    \begin{equation}
        \widehat{v}_h(\widehat{\vec{x}}):= v_h\left(F_K^{-1}(\widehat{\vec{x}})\right):\widehat{K}\to \mathbb{R}.
    \end{equation}    
    We can then compute via a scaling argument the following bound,
    \begin{align}
        \label{eq:scalingArgument}
        \norm{\partial_n u_h}_{L^2(\partial\Omega\cap K)} &= \norm{\nabla u_h \cdot n}_{L^2(\partial\Omega\cap K)}\leq \norm{\vec{n}}_{L^\infty(\partial\Omega\cap K)}\norm{\nabla u_h}_{L^2(\partial\Omega\cap K)}\\
        & \leq \norm{\mathbb{A}}^{-1}\left(\frac{\abs{\partial\Omega\cap K}}{\abs{F_K(\partial\Omega\cap K)}}\right)^{\frac{1}{2}}\norm{\nabla \widehat{u}_h}_{L^2\left(F_K(\partial\Omega\cap K)\right)}.
    \end{align}
    We now notice that since all norms in finite dimension are equivalent there exists a constant $C_1$ such that
    \begin{equation}
        \label{eq:normEquivalence}
        \norm{\nabla \widehat{u}_h}_{L^2\left(F_K(\partial\Omega\cap K)\right)}\leq C_1 \norm{\nabla \widehat{u}_h}_{L^2(\widehat{K})},
    \end{equation}
    where $C_1$ may depend on $ \alpha$ and $\gamma$ but it does not depend on $h$ because we are working on the reference element.
    Combining \eqref{eq:normEquivalence} with \eqref{eq:scalingArgument} we have
    \begin{align}
        \norm{\partial_n u_h}_{L^2(\partial\Omega\cap K)} &\leq C_1 \norm{\mathbb{A}}^{-1}\left(\frac{\abs{\partial\Omega\cap K}}{\abs{F_K(\partial\Omega\cap K)}}\right)^{\frac{1}{2}}\norm{\nabla \widehat{u}_h}_{L^2(\widehat{K})}\\
        &\leq C_1 \norm{\mathbb{A}}^{-1}\norm{\mathbb
        {A}}\left(\frac{\abs{\partial\Omega\cap K}}{\abs{F_K(\partial\Omega\cap K)}}\right)^{\frac{1}{2}}\left(\frac{\abs{\widehat{K}}}{\abs{K}}\right)^{\frac{1}{2}}\norm{\nabla u_h}_{L^2(K)}\\
        &\leq C_1\left(\frac{\abs{\widehat{K}}}{\abs{F_K(\partial\Omega\cap K)}}\right)^{\frac{1}{2}}\left(\frac{\abs{\partial \Omega\cap K}}{\abs{K}}\right)^{\frac{1}{2}}\norm{\nabla u_h}_{L^2(K)}\label{eq:metaTraceInequality}
    \end{align}
    We proceed bounding $\left(\frac{\abs{\widehat{K}}}{\abs{F_K(\partial\Omega\cap K)}}\right)^{\frac{1}{2}}\left(\frac{\abs{\partial \Omega\cap K}}{\abs{K}}\right)^{\frac{1}{2}}$ using \eqref{eq:GeoAssumption1}.
    We begin observing that $\widehat{K}$ has a unit area and that from \eqref{eq:GeoAssumption1} we know
    \begin{equation}
         \abs{F_K(\partial\Omega\cap K)} \geq \widehat{C}_\Omega h^{{ \alpha}-1}\;\Rightarrow\; \abs{F_K(\partial\Omega\cap K)}^{-1} \leq \widehat{C}_\Omega^{-1} h^{1-{ \alpha}}.
         \label{eq:GeoBound1}
    \end{equation}
    Furthermore $\abs{K}=h^2$ and from \eqref{eq:GeoAssumption1} we know that
    \begin{equation}
        \frac{\abs{\partial \Omega\cap K}}{\abs{K}}\leq C_\Omega h^{\gamma-2}. \label{eq:GeoBound2}
    \end{equation}
    Combining \eqref{eq:GeoBound1}, \eqref{eq:GeoBound2} with \eqref{eq:metaTraceInequality} we have,
    \begin{align}
        \norm{\partial_n u_h}_{L^2(\partial\Omega\cap K)} &\leq C_1\left(\frac{\abs{\widehat{K}}}{\abs{F_K(\partial\Omega\cap K)}}\right)^{\frac{1}{2}}\left(\frac{\abs{\partial \Omega\cap K}}{\abs{K}}\right)^{\frac{1}{2}}\norm{\nabla u_h}_{L^2(K)}\\
        &\leq C_1C_\Omega^{\frac{1}{2}}\widehat{C}_\Omega^{-\frac{1}{2}}h^{\frac{\gamma-{ \alpha}-1}{2}}\norm{\nabla u_h}_{L^2(K)},
    \end{align}
    lastly we call $C_2^{\frac{1}{2}}$ the constant equal to $C_1\widehat{C}_\Omega^{-\frac{1}{2}} C_\Omega^{\frac{1}{2}}$.
\end{proof}

To study the a priori convergence rate of \eqref{eq:discreteVarProblem} we need to introduce a discrete norm, i.e.,
\begin{equation}
    \normh{\cdot}^2 = \abs{\cdot}_{H^1(\Omega)}^2+h^{-\alpha}C_a \norm{\cdot}^2_{L^2(\Gamma_D)}+ C_a^{-1}h^{\alpha} \norm{\partial_n \cdot}^2_{L^2(\Gamma_D)},
\end{equation}
which is well posed for any argument in $H^{\frac{3}{2}+\varepsilon}(\Omega)$.
\begin{lemma}
    The bilinear form $a_{h,\alpha}:H^{\frac{3}{2}+\varepsilon}(\Omega)\times H^{\frac{3}{2}+\varepsilon}(\Omega) \to \mathbb{R}$, defined as in \eqref{eq:ah}, is continuous and coercive, when restricted to $V_h$, with respect to $\norm{\cdot}_{h,\alpha}$, provided we choose $ C_a = \frac{1}{2}C_2^{-1}$ and $ \gamma \geq 1$.
\end{lemma}
\begin{proof}
    We begin proving coercivity. We have
    \begin{align}
        a_h(v_h,v_h) & =  (\nabla v_h,\nabla v_h)_{L^2(\Omega)} + C_a h^{-\alpha} (v_h,v_h)_{L^2(\Gamma_D)}-2(\partial_n v_h,v_h)_{L^2(\Gamma_D)}\\ 
        & = \abs{v_h}^2_{H^1(\Omega)} + C_ah^{-\alpha}\norm{v_h}^2_{L^2(\Gamma_D)}-2(C_a^{-\frac{\alpha}{2}}h^{\frac{\alpha}{2}}\partial_n v_h,C_a^{\frac{\alpha}{2}}h^{-\frac{\alpha}{2}}v_h)_{L^2(\Gamma_D)}\\
        &\geq \abs{v_h}^2_{H^1(\Omega)} + (1-\varepsilon) C_a^\alpha h^{-\alpha}\norm{v_h}^2_{L^2(\Gamma_D)} - \varepsilon^{-1}C_a^{-\alpha}h^{\alpha}\norm{\partial_n v_h}^2_{L^2(\Gamma_D)},
    \end{align}
    where the last inequality comes from Young's inequality with $\varepsilon\in(0,1)$.
    We add and subtract $\varepsilon \abs{v_h}_{H^1(\Omega)}^2$, fix $ \gamma\geq 1$ and use \eqref{eq:discreteTraceInequality} to control from below $\varepsilon \abs{v_h}_{H^1(\Omega)}^2$, i.e.,
    \begin{align}
        a_{h,\alpha}(v_h,v_h) &\geq \abs{v_h}^2_{H^1(\Omega)} + (1-\varepsilon) C_a^\alpha h^{-\alpha}\norm{v_h}^2_{L^2(\Gamma_D)} - \varepsilon^{-1}C^{-\alpha}_a h^{\alpha}\norm{\partial_n v_h}^2_{L^2(\Gamma_D)}\\
        & \geq (1-\varepsilon) \abs{v_h}^2_{H^1(\Omega)} + (1-\varepsilon) C_a^\alpha h^{-\alpha}\norm{v_h}^2_{L^2(\Gamma_D)}\\
        &\qquad\qquad\qquad\qquad+ (2\varepsilon-\varepsilon^{-1})C_a^{-\alpha} h^{\alpha}\norm{\partial_n v_h}^2_{L^2(\Gamma_D)}\geq C_3 \norm{v_h}_{h,\alpha},
    \end{align}   
    where { we choose} $C_a =  \frac{1}{2}C_2^{-1}$, and $C_3= (1-\varepsilon^*)$, with $\varepsilon^*$ being the positive root of $1-\varepsilon = { 2}\varepsilon - \varepsilon^{-1}$. In particular, notice that $C_3 \geq \frac{1}{5}$.
    We now proceed to prove the continuity of the bilinear form for $u,v \in H^{\frac{3}{2}+\varepsilon}(\Omega)$. We have
        \begin{align}
        &a_{h,\alpha}(u,v) & \leq \abs{u}_{H^1(\Omega)}\abs{v}_{H^1(\Omega)}&+ \norm{h^{\frac{\alpha}{2}}C_a^{-\frac{\alpha}{2}}\partial_n u}_{L^2(\Gamma_D)}\norm{h^{-\frac{\alpha}{2}}C_a^{\frac{\alpha}{2}}v}_{L^2(\Gamma_D)}\\
        &&&+ \norm{h^{\frac{\alpha}{2}}C_a^{-\frac{\alpha}{2}}\partial_n v}_{L^2(\Gamma_D)}\norm{h^{-\frac{\alpha}{2}}C_a^{\frac{\alpha}{2}}u}_{L^2(\Gamma_D)}\\
        &&&+ \norm{h^{\frac{\alpha}{2}}C_a^{-\frac{\alpha}{2}}u}_{L^2(\Gamma_D)}\norm{h^{\frac{\alpha}{2}}C_a^{-\frac{\alpha}{2}}v}_{L^2(\Gamma_D)}\\
        &&&{+   \norm{h^{\frac{\alpha}{2}}C_a^{-\frac{\alpha}{2}}\partial_{n}u}_{L^2(\Gamma_D)}\norm{h^{\frac{\alpha}{2}}C_a^{-\frac{\alpha}{2}}\partial_{n}v}_{L^2(\Gamma_D)}}\\
        & & \leq \abs{u}_{H^1(\Omega)}\abs{v}_{H^1(\Omega)}&+ \left(h^{\frac{\alpha}{2}} C_a^{-\frac{\alpha}{2}}\norm{\partial_n u}_{L^2(\Gamma_D)} + h^{-\frac{\alpha}{2}} C_a^{\frac{\alpha}{2}}\norm{u}_{L^2(\Gamma_D)} \right)\\
        &&&\cdot \left(h^{\frac{\alpha}{2}} C_a^{-\frac{\alpha}{2}}\norm{\partial_n v}_{L^2(\Gamma_D)} + h^{-\frac{\alpha}{2}} C_a^{\frac{\alpha}{2}}\norm{v}_{L^2(\Gamma_D)} \right)
    \end{align}
    \begin{align}
        a_{h,\alpha}(u,v)&\leq\left [\abs{u}^2_{H^1(\Omega)}+\left(h^{\frac{\alpha}{2}} C_a^{-\frac{\alpha}{2}}\norm{\partial_n u}_{L^2(\Gamma_D)} + h^{-\frac{\alpha}{2}} C_a^{\frac{\alpha}{2}}\norm{u}_{L^2(\Gamma_D)} \right)^2\right]^{\frac{1}{2}}\\
        &\cdot\;\; \left [\abs{v}^2_{H^1(\Omega)}+\left(h^{\frac{\alpha}{2}} C_a^{-\frac{\alpha}{2}}\norm{\partial_n v}_{L^2(\Gamma_D)} + h^{-\frac{\alpha}{2}} C_a^{\frac{\alpha}{2}}\norm{v}_{L^2(\Gamma_D)} \right)^2\right]^{\frac{1}{2}}\\
        &\leq \left(\abs{u}^2_{H^1(\Omega)}+{ 2h^{\alpha} C_a^{-\alpha}\norm{\partial_n u}^2_{L^2(\Gamma_D)} + 2h^{-\alpha} C_a^{\alpha}\norm{u}^2_{L^2(\Gamma_D)}}\right)^{\frac{1}{2}}\\
        &\cdot\;\; \left(\abs{v}^2_{H^1(\Omega)}+{  2h^{\alpha} C_a^{-\alpha}\norm{\partial_n v}^2_{L^2(\Gamma_D)} + 2h^{-\alpha} C_a^{\alpha}\norm{v}^2_{L^2(\Gamma_D)}}\right)^{\frac{1}{2}}\\
        &\leq C_4 \norm{u}_{h,\alpha}\norm{v}_{h,\alpha},
    \end{align}
    hence $C_4=2$ is the continuity constant.
\end{proof}
\begin{corollary}
Using Lax--Milgram lemma we know that for any $(f,g)\in L^2(\Omega)\times H^{\frac{1}{2}}(\partial\Omega)$ there exists a unique $u_h\in V_h$ solving \eqref{eq:discreteVarProblem}.
\end{corollary}
\begin{proposition}
    \label{prop:discreteAPriori}
    The discrete solution $u_h\in V_h$ of \eqref{eq:discreteVarProblem} is an approximation of the solution $u\in H^{1+\delta}(\Omega)$, with respect to $\norm{\cdot}_{h,\alpha}$ i.e.,
    \begin{equation}
        \norm{u-u_h}_{h,\alpha} \leq \underset{v_h\in V_h}{\inf}\norm{u-v_h}_{h,\alpha} \leq C_5 { h^{\min\{\delta,1,\frac{1}{2}+\delta-\frac{\alpha}{2},\frac{\alpha}{2}\}}}\abs{u}_{H^{1+\delta}(\Omega)},
    \end{equation}
    provided \eqref{eq:GeoAssumption1} holds, $\delta\in (\frac{1}{2},\frac{3}{2}]$ and $ \gamma \geq 1$.
\end{proposition}
\begin{proof}
    Via triangle inequality, we know that for any $v_h\in V_h$,
    \begin{align}
        \norm{u-u_h}_{h,\alpha}\leq \norm{u-v_h}_{h,\alpha}+\norm{v_h-u_h}_{h,\alpha},
    \end{align}
    we can now bound the second term in the inequality, {  using the coercivity of the bilinear form $a_{h,\alpha}(\cdot,\cdot)$ with constant $\frac{1}{5}$:}
    \begin{align}
        \label{eq:GA}
        \norm{v_h-u_h}_{h,\alpha}^2\leq 5a_h(v_h-u_h,v_h-u_h)=5a_h(u-v_h,u_h-v_h),
    \end{align}
    where the last identity is a consequence of Galerkin orthogonality, i.e.,
    \begin{align}
        a_h(u_h-v_h,u_h-v_h) &= \nonumber a_h(u_h,u_h-v_h) - a_h(v_h,u_h-v_h) \\ \nonumber
        &= (f,u_h-v_h)_{L^2(\Omega)} - a_h(v_h,u_h-v_h) \\ \nonumber
        &= a_h(u,u_h-v_h)-a_h(v_h,u_h-v_h)\\
        &=a_h(u-v_h,u_h-v_h).
    \end{align}
    Since $1+\delta>3/2$, from~\eqref{eq:GA} and the continuity of $a_h$, we obtain
    \begin{align}
        \norm{v_h-u_h}_{h,\alpha}^2\leq 5a_h(u-v_h,u_h-v_h){ \leq}10\norm{v_h-u_h}_{h,\alpha}\norm{u-v_h}_{h,\alpha}.
    \end{align}
    Hence we have that $\forall v_h \in V_h$ the following inequality holds,
    \begin{align}
         \norm{u-u_h}_{h,\alpha} \leq 11\norm{u-v_h}_{h,\alpha}
    \end{align}
    and, since the previous inequality holds for all $v_h$ in $V_h$ we have that,
    \begin{align}
        \norm{u-u_h}_{h,\alpha}\leq 11\underset{v_h\in V_h}{\inf} \norm{u-v_h}_{h,\alpha}.
    \end{align}
    { 
    We choose as $v_h$ the nodal $\mathbb{Q}_1$ interpolant $\mathcal{I}_h u$ of $u$ and proceed to bound every term in $\norm{u-\mathcal{I}_h u}_{h,\alpha}$, i.e.
    \begin{equation}
        \norm{u-\mathcal{I}_h u}_{h,\alpha}^2 = \abs{u-\mathcal{I}_h u}_{H^1(\Omega)}^2 + h^{-\alpha}C_a\norm{u-\mathcal{I}_h u}_{L^2(\Gamma_D)}^2 + C_a^{-1}h^{\alpha}\norm{\partial_n(u-\mathcal{I}_h u)}_{L^2(\Gamma_D)}^2.
    \end{equation}
    From standard interpolation estimates we know that 
    \begin{equation}
        \abs{u-\mathcal{I}_h u}_{H^1(\Omega)}\leq D_1 h^{\min\{\delta,1\}}\abs{u}_{H^{1+\delta}(\Omega)}.
    \end{equation}
    Since the nodal interpolant along the polygonal boundary corresponds to the one-dimensional interpolant we have that
    \begin{equation}
        \norm{u-\mathcal{I}_h u}_{L^2(\Gamma_D)}\leq D_2 h^{\frac{1}{2}+\delta}\abs{u}_{H^{\frac{1}{2}+\delta}(\partial \Omega)}.
    \end{equation}
    We now need to bound the last remaining term, i.e.
    \begin{align}
        \norm{\partial_n(u-\mathcal{I}_h u)}_{L^2(\Gamma_D)}^2&\leq \sum_{F} \norm{\partial_n(u-\mathcal{I}_h u)}_{L^2(F)}^2 \leq \sum_{F} \norm{\nabla(u-\mathcal{I}_h u)}_{L^2(F)}^2\\
        &\leq D_3 \sum_K \norm{\nabla(u-\mathcal{I}_h u)}_{H^{\delta}(K)}^2 \leq D_3 \norm{u}_{H^{1+\delta}(\Omega)}^2,
    \end{align}
    where $F$ are the faces of $\Omega$, i.e. the polygonal domain described interesting the boundary $\Gamma_D$ with the points where the level set function evaluates to zero and $K$ are the elements assosciated to the faces $F$.
    The last bound comes from the element wise stability of the nodal interpolant with respect tot the $H^{1+\delta}$ norm.
    Combining all of the previous estimates we have that
    \begin{equation}
        \norm{u-\mathcal{I}_h u}_{h,\alpha}\leq D_4 h^{\min\{\delta,1,\frac{1}{2}+\delta-\frac{\alpha}{2},\frac{\alpha}{2}\}}\abs{u}_{H^{1+\delta}(\Omega)}.
    \end{equation}
    }
\end{proof}
We can now proceed to use the Aubin--Nitsche duality argument to give an estimate of the convergence rate of the scheme here proposed with respect to the $L^2(\Omega)$ norm.
\begin{proposition}
    Under the assumption that \eqref{eq:GeoAssumption1} holds the discrete solution $u_h\in V_h$ of \eqref{eq:discreteVarProblem} is an approximation of the solution $u\in H^{1+\delta}(\Omega)$ with respect to the $L^2(\Omega)$ norm i.e.,
    \begin{equation}
        \norm{u-u_h}_{L^2(\Omega)} \leq C_6{ h^{2\min\{\delta,1,\frac{1}{2}+\delta-\frac{\alpha}{2},\frac{\alpha}{2}\}}\abs{u}_{H^{1+\delta}(\Omega)}},\label{eq:L2}
    \end{equation}
    provided \eqref{eq:GeoAssumption1} holds, $ \delta\in (\frac{1}{2},\frac{3}{2}]$ and $ \gamma \geq 1$ { and that the adjoint problem is sufficiently regular.}
\end{proposition}
\begin{proof}
    We consider the auxiliary problem, find $w\in H^1_{\Gamma_D}(\Omega)$ such that 
    \begin{equation}
        a(w,z) = (u-u_h,z)_{L^2(\Omega)}, \qquad \forall z \in H^1_{\Gamma_D}(\Omega),
    \end{equation}
    {  and we assume $w\in H^{1+\delta}(\Omega)$}.
    We now notice that
    \begin{align}
        \norm{u-u_h}^2_{L^2(\Omega)} &= \left(\nabla w,\nabla(u-u_h)\right)_{L^2(\Omega)}-(u-u_h,\partial_n w)_{L^2(\partial\Omega)}\\
        &= \left(\nabla w,\nabla(u-u_h)\right)_{L^2(\Omega)}-(u-u_h,\partial_n w)_{L^2(\Gamma_D)}\\
        &-(\partial_n u-\partial_n u_h,w)_{L^2(\Gamma_D)}+C_ah^{-1}\left(u-u_h,w\right)_{L^2(\Gamma_D)}\\
        &=a_h(u-u_h,w).
    \end{align}
    We now notice that $a_h(u-u_h,w_h) = (f,w_h)_{L^2(\Omega)}-(f,w_h)_{L^2(\Omega)}=0$ therefore using the continuity of $a_h$ we can rewrite the previous identity as
    \begin{align}
        \norm{u-u_h}^2_{L^2(\Omega)} & = a_h(u-u_h,w-w_h)\leq \norm{u-u_h}_{h,1}\norm{w-w_h}_{h,1}.
    \end{align}
    Using Proposition~\ref{prop:discreteAPriori} we can bound the $L^2(\Omega)$ error as
    \begin{align}
        \norm{u-u_h}^2_{L^2(\Omega)} & \leq C_5^2h^\delta\abs{u}_{H^2(\Omega)}h^\delta\abs{w}_{H^2(\Omega)}\leq C_5 h^{2\delta}\norm{u-u_h}_{L^2(\Omega)}\abs{u}_{H^2(\Omega)}.
    \end{align}
    Dividing the previous inequality by $\norm{u-u_h}_{L^2(\Omega)}$ we obtain the desired result.
\end{proof}
We now focus our attention on the analysis of the snapping back to the grid perturbation discussed in Section~\ref{sec:comp}.

When adopting a snapping back to grid technique we are solving \eqref{eq:discreteVarProblem} on a domain $\widetilde{\Omega}$, which is a perturbation of the original $\Omega$.
We introduce this perturbation to get rid of cells that do not respect the geometrical assumption \eqref{eq:GeoAssumption1}.

We aim to use the theory presented in \cite{SS}, and to develop an a priori error estimate for the snapping back to grid approach viewed as a perturbation method.

Before diving into the results presented in \cite{SS} we need to discuss the eventual regularity issue that might arise from considering the perturbed domain $\widetilde{\Omega}$.
In the one-dimensional case, we know that the snapping presented in Section~\ref{sec:comp} does not affect the elliptic regularity discussed at the beginning of this section.
Instead, in the two-dimensional case the worst possible scenario comes when, as a consequence of the snapping back to grid, the perturbed domain $\widetilde{\Omega}$ presents a re-entrant corner of aperture $\frac{3}{2}\pi$.
An example of such a situation can be seen in Figure \ref{figure_class_cell_generic_cell}(a).

In this case elliptic regularity only guarantees that the $\widetilde{u}\in H^1_{\Gamma_D}(\widetilde\Omega)$ solving
\begin{equation}
    (\nabla \widetilde{u},\nabla v)_{L^2(\widetilde{\widetilde\Omega})} = (f,v)_{L^2(\widetilde\Omega)}+\prescript{}{H^{-\frac{1}{2}}(\widetilde\Gamma_N)}{\langle} g_N,v\rangle_{H^{\frac{1}{2}}(\widetilde\Gamma_N)}, \qquad \forall v \in H^1_{\Gamma_D}(\widetilde\Omega),\label{eq:perturbedVarProblem}
\end{equation}
{ lives in the Sobolev space $H^{\frac{5}{3}+\varepsilon}(\widetilde\Omega)$ no matter how regular $f$, $g_D$, and $g_N$ are.
Hence we might have a mismatch between the regularity of $\widetilde{u}$ and $u$, and this should be taken into consideration when using \eqref{eq:L2}.}

Let us state the following simplified result from \cite[Theorem 8.4]{SS} to bound the $L^2$ distance over $R$ between $\widetilde{u}$ and $u$.
\begin{theorem}
    \label{thm:SS}
    Let $\Omega \subset R$ and $\widetilde{\Omega}\subset R$ be Lipschitz domain with corners of aperture at most $\frac{3}{2}\pi$.
    We consider $u\in H^1(\Omega)$ and $\widetilde{u}\in H^{1}(\widetilde{\Omega})$ solutions of the original variational problem and \eqref{eq:perturbedVarProblem} for a data $f\in L^2(R)$ and with homogeneous Dirichlet boundary condition, i.e., $g_D=0$ on the entire boundary, i.e., $\Gamma_D = \partial\Omega$ and $\widetilde\Gamma_D = \partial \widetilde{\Omega}$ respectively.
    Under these hypotheses the following bound for the $L^2$ distance over $R$ between $\widetilde{u}$ and $u$, holds:
    \begin{equation}
        \norm{u-\widetilde{u}}_{L^2(R)}\leq C_7 \left(\sup_{\vec{x}\in \Omega\Delta \widetilde{\Omega}}d\left(\vec{x},\partial \Omega \right)\right)\norm{f}_{L^2(R)},
    \end{equation}
    where $\Omega\Delta\widetilde{\Omega}$ denotes the set difference between $\Omega$ and $\widetilde{\Omega}$ and $C_7$ is a constant depending only the diameter of $R$ and the aperture of all the re-entrant corners of $\widetilde{\Omega}$.
\end{theorem}
\begin{corollary}
    Let $\widetilde{u}_h\in V_h$ be the solution of the perturbed discrete variational problem, 
    \begin{align}
        (\nabla \widetilde{u}_h,\nabla {v}_h)_{L^2(\widetilde {\Omega})} &+ C_a h^{-\alpha} (\widetilde{u}_h,{v}_h)_{L^2(\Gamma_D)}
        \\&-(\partial_n\widetilde{u}_h,{v}_h)_{L^2(\Gamma_D)}-(\widetilde{u}_h,\partial_n {v}_h)_{L^2(\Gamma_D)}= (f,{v}_h)_{L^2(\widetilde\Omega)}, \;\; \forall {v}_h \in V_h,\nonumber
    \end{align}
    where $V_h$ is the usual space of $\mathbb{Q}_1$ functions over each element of a mesh $\mathcal{T}_h$ for $R\subset \mathbb{R}$. 
    If $\widetilde{\Omega}$ is the geometry obtained from $\Omega$ after snapping back to cell for every $K$ such that
    \begin{equation}
        \abs{\partial\Omega\cap K}\leq \widehat{C}_\Omega h^{ \alpha},
    \end{equation}
    then for $\gamma\geq 1$ there exists $\widetilde\delta \in (\frac{1}{2},\frac{5}{2}]$, depending only on the regularity of $\partial\widetilde{\Omega}$, such that
    \begin{equation}
         \norm{u-\widetilde{u}_h}_{L^2(R)}\leq C_8 h^{\min\{\widetilde\delta,1,\frac{1}{2}+\delta-\frac{\alpha}{2},\frac{\alpha}{2}\}}\norm{f}_{L^2(R)}.
    \end{equation}
\end{corollary}
\begin{proof}
    The proof of this result is almost immediate.
    { 
    It is sufficient to observe that in Algorithm~\ref{alg_snap} we are snapping back to the grid the cells that do not respect the geometrical assumption \eqref{eq:GeoAssumption1}, hence the distance between the newly obtained domain $\widetilde{\Omega}$ and the original domain $\Omega$ is bounded by the diameter of the cells that have been snapped back, implying that
    \begin{equation}
        \sup_{\vec{x}\in \Omega\Delta \widetilde{\Omega}}d\left(\vec{x},\partial \Omega \right)\leq \sqrt{2}\widehat{C}_\Omega h^{\alpha}.
    \end{equation}
    }
    We then observe that as discussed previously $\widetilde{\Omega}$ can only have re-entrant corners of aperture $\frac{3}{2}\pi$ and therefore we can apply the previous theorem to obtain,
    \begin{equation}
        \norm{u-\widetilde{u}}_{L^2(R)}\leq \sqrt{2}\widehat{C}_\Omega  C_7 h^{ \alpha}\norm{f}_{L^2(R)}.
    \end{equation}
    Since we know that $\widetilde{u}\in H^{1+\widetilde\delta}$, with $ \widetilde\delta\in(\frac{1}{2},\frac{5}{2}]$, and that \eqref{eq:GeoAssumption1} is verified by the pair $(\mathcal{T}_h,\widetilde\Omega)$ if $ \gamma\geq 1$ we can apply \eqref{eq:L2} to obtain the following bound,
    \begin{equation}
        \norm{\widetilde{u}-\widetilde{u}_h}_{L^2(\widetilde{\Omega})}\leq C_6 h^{ 2\min\{\widetilde\delta,1,\frac{1}{2}+\widetilde\delta-\frac{\alpha}{2},\frac{\alpha}{2}\}}\norm{f}_{L^2(\widetilde{\Omega})}.
    \end{equation}
    Using the fact that $\norm{\cdot}_{L^2(\widetilde{\Omega})}\leq \norm{\cdot}_{L^2(R)}$ and the triangular inequality we get that,
    \begin{equation}
        \norm{u-\widetilde{u}_h}_{L^2(R)}\leq C_8 h^{ 2\min\{\widetilde\delta,1,\frac{1}{2}+\widetilde\delta-\frac{\alpha}{2},\frac{\alpha}{2}\}}\norm{f}_{L^2(R)},
    \end{equation}
    where $u$ and $\widehat u$ are extended to zero outside of $\Omega$ and $\widehat\Omega$ and $C_8 = (C_6+\sqrt{2}\widehat{C}_\Omega  C_7)$.
\end{proof}
\section{Numerical tests}
\label{sec:num}
In this section, we provide numerical tests showing second-order accuracy in one and two-dimensions. We test different smooth geometries, and for each one we show the accuracy of the solution and of its divergence, which are computed as follows:
\begin{align}
\label{eq:error}
    {\rm error} & = \frac{\sqrt{\int_\mathcal{A} |u_h - u_{\rm exa}|^2 dx dy}}{\sqrt{\int_\mathcal{A} |u_{\rm exa}|^2 dx dy}} \\ \label{eq:Derror}
    \nabla{\rm error} & = \frac{\sqrt{\int_\mathcal{A} |\nabla \cdot u_h - \nabla \cdot u_{\rm exa}|^2 dx dy}}{\sqrt{\int_\mathcal{A} |\nabla \cdot u_{\rm exa}|^2 dx dy}},
\end{align}
where $u_h = \sum_j u_h^j\varphi_j$ and $\nabla \cdot u_h = \sum_j u_h^j\nabla \cdot \varphi_j$. The integral over $\mathcal{A}$ is computed using the midpoint quadrature rule. The number of nodes in the quadrature rule must be chosen such that it is incommensurate with the number of points in the spatial discretization. In the accuracy analysis, the penalisation parameter plays an important role. We find a compromise between accuracy and a well-conditioned numerical scheme.

\subsection{Numerical tests in 1D}
Here, we perform numerical tests to check that the method achieves the expected accuracy
independently of the parameters $\theta_1$ and $\theta_2$ (see Fig.~\ref{fig:1D_setup}).

\begin{figure}
    \centering
 \begin{overpic}[abs,width=0.5\textwidth,unit=1mm,scale=.25]{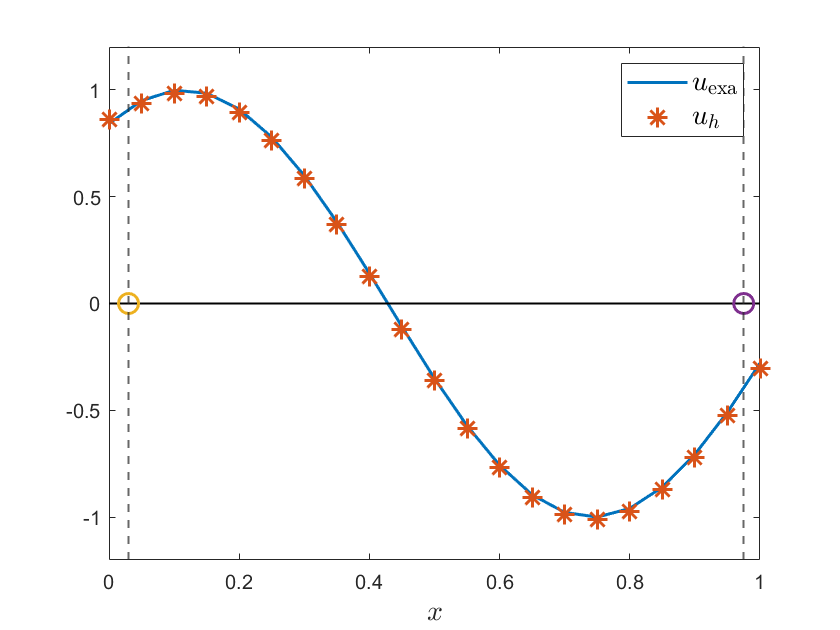}
\put(-2,48){(a)}\end{overpic} \\  
\begin{minipage}[b]{.31\textwidth}
    \centering
 \begin{overpic}[abs,width=1.1\textwidth,unit=1mm,scale=.25]{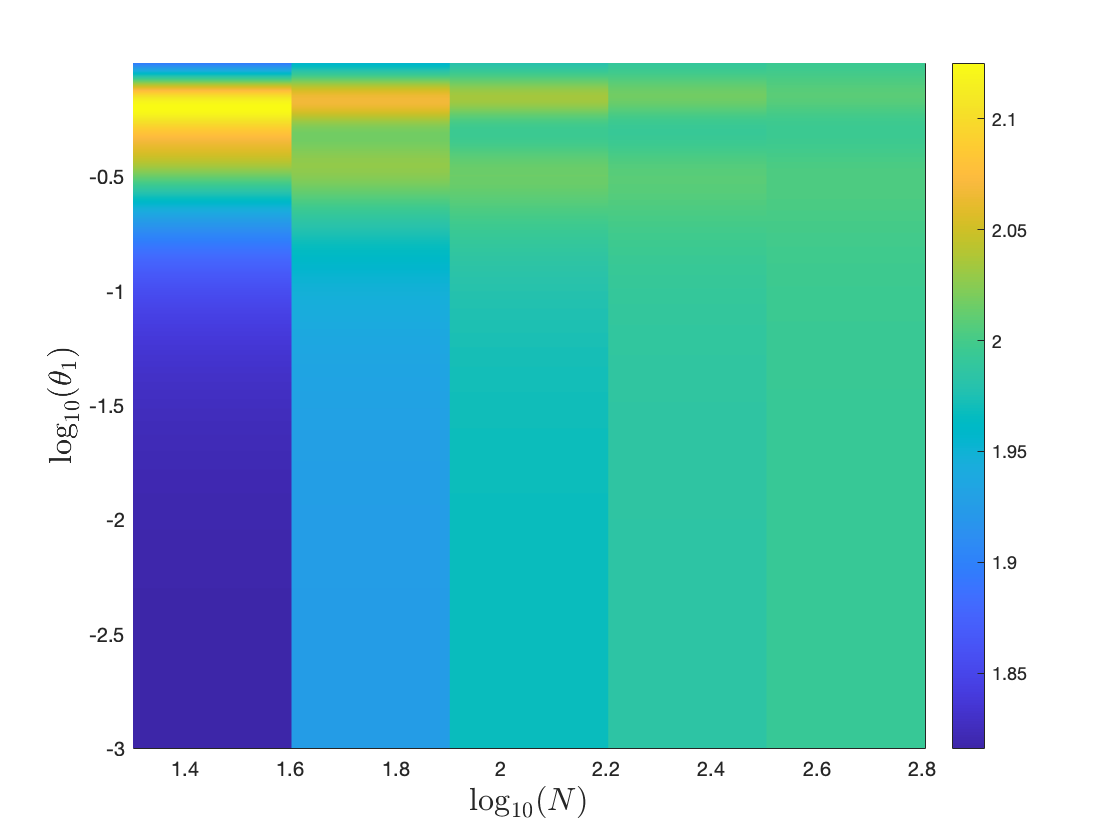}
\put(-0,34){(b)}
\end{overpic}     
\end{minipage}
\begin{minipage}[b]{.31\textwidth}
    \centering
 \begin{overpic}[abs,width=1.1\textwidth,unit=1mm,scale=.25]{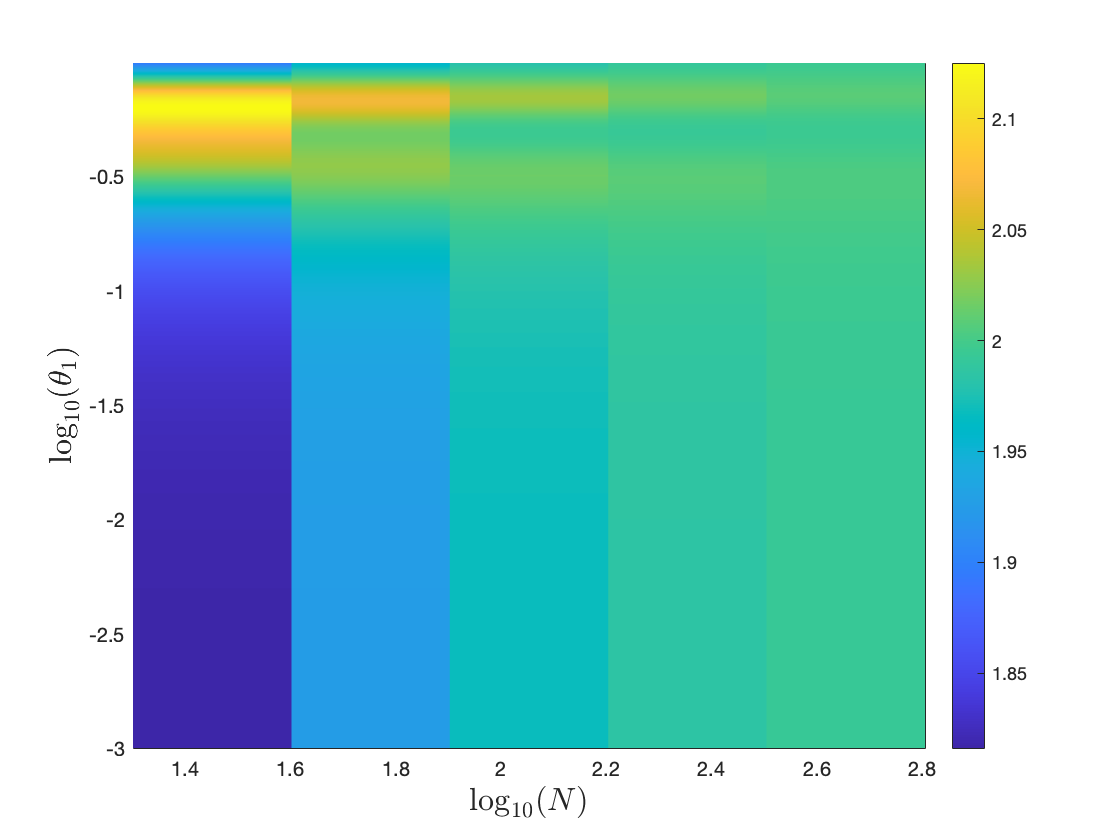}
\put(-0,34){(c)}
\end{overpic}     
\end{minipage}
\begin{minipage}[b]{.31\textwidth}
    \centering
 \begin{overpic}[abs,width=1.1\textwidth,unit=1mm,scale=.25]{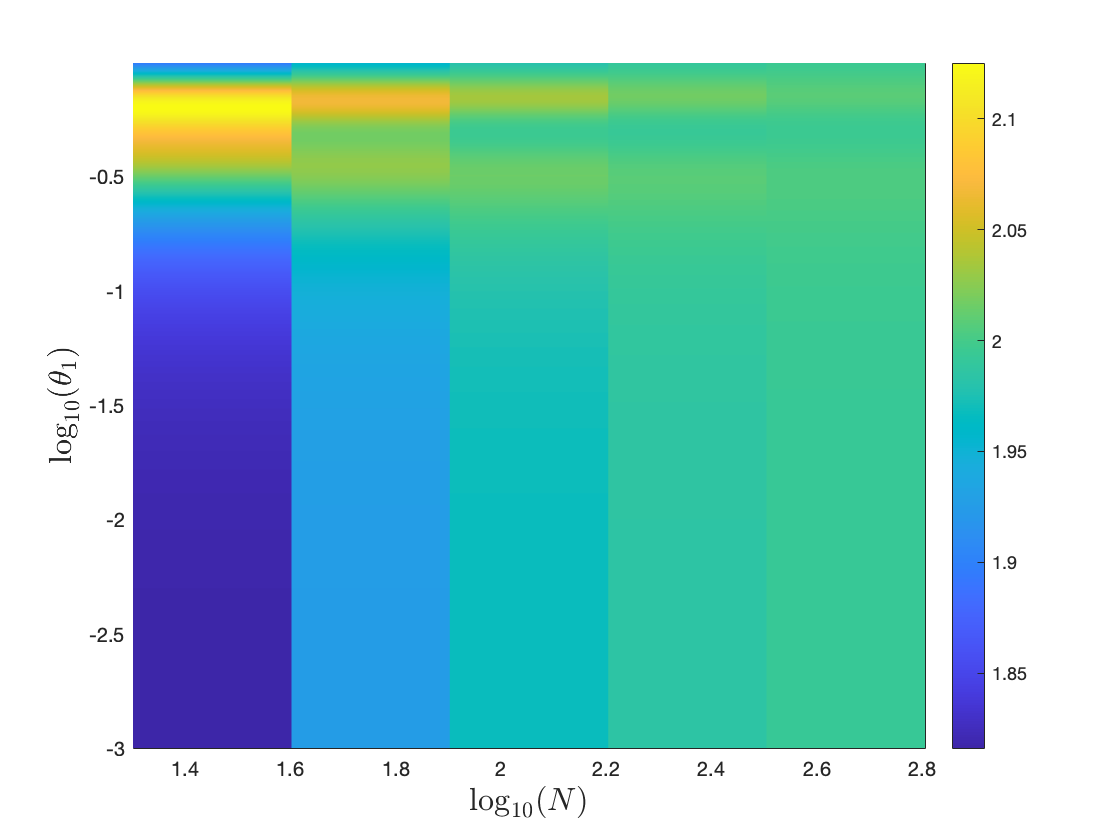}
\put(-0,34){(d)}
\end{overpic}     
\end{minipage}
\begin{minipage}[b]{.31\textwidth}
    \centering
 \begin{overpic}[abs,width=1.1\textwidth,unit=1mm,scale=.25]{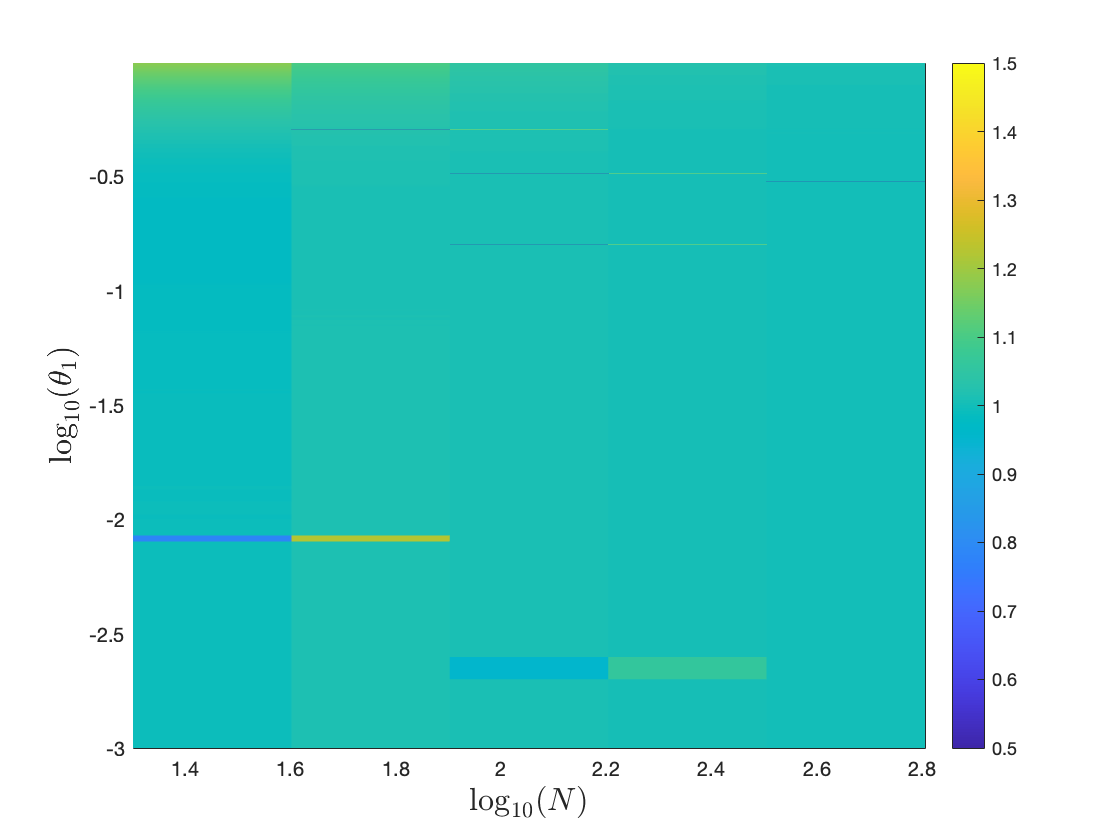}
\put(-0,34){(e)}
\end{overpic}     
\end{minipage}
\begin{minipage}[b]{.31\textwidth}
    \centering
 \begin{overpic}[abs,width=1.1\textwidth,unit=1mm,scale=.25]{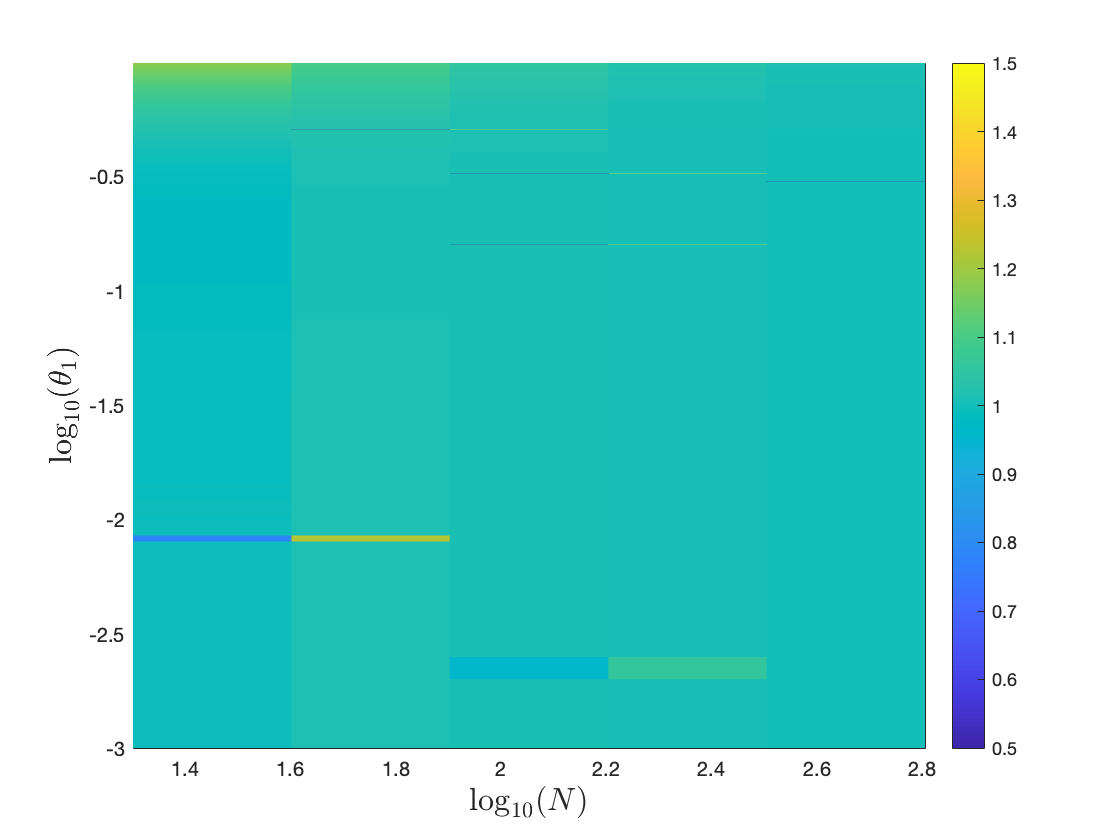}
\put(-0,34){(f)}
\end{overpic}     
\end{minipage}
\begin{minipage}[b]{.31\textwidth}
    \centering
 \begin{overpic}[abs,width=1.1\textwidth,unit=1mm,scale=.25]{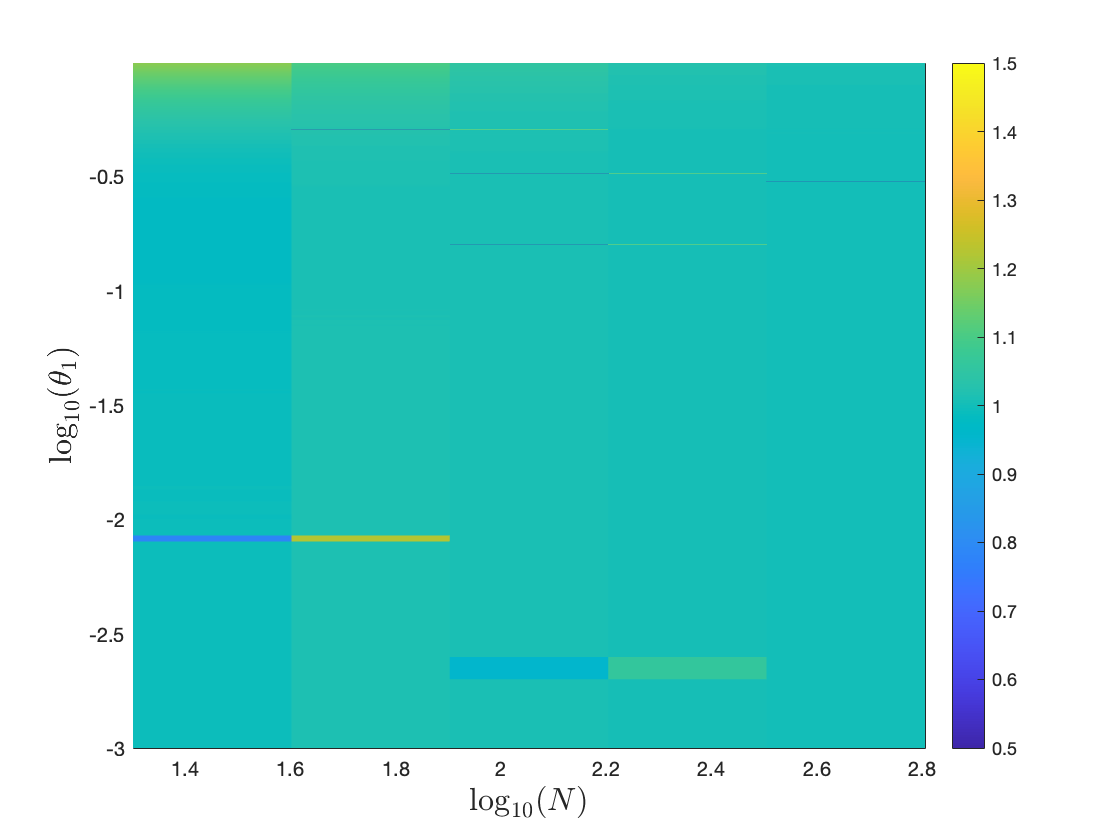}
\put(-0,34){(g)}
\end{overpic}     
\end{minipage}
\begin{minipage}[b]{.31\textwidth}
    \centering
\begin{overpic}[abs,width=1.1\textwidth,unit=1mm,scale=.25]{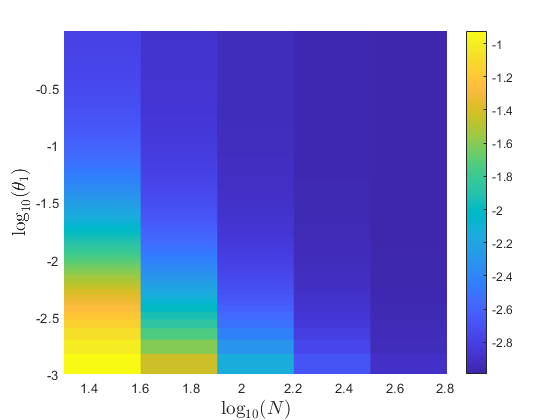}
\put(-0,34){(h)}
\end{overpic}      
\end{minipage}
\begin{minipage}[b]{.31\textwidth}
    \centering
\begin{overpic}[abs,width=1.1\textwidth,unit=1mm,scale=.25]{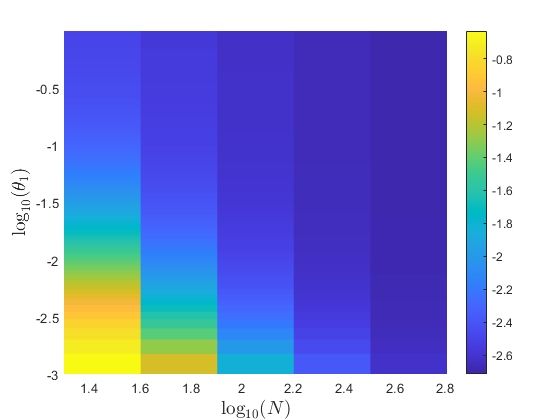}
\put(-0,34){(i)}
\end{overpic}      
\end{minipage}
\begin{minipage}[b]{.31\textwidth}
    \centering
\begin{overpic}[abs,width=1.1\textwidth,unit=1mm,scale=.25]{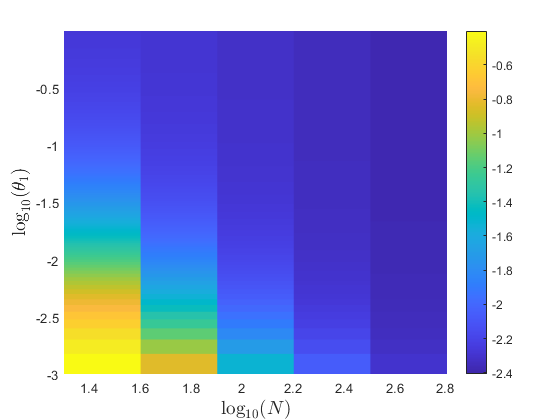}
\put(-0,34){(l)}
\end{overpic}      
\end{minipage}
\caption{\textit{(a): Numerical solution of the mixed problem \eqref{eq:system1M} (asterisks) and comparison
with the exact solution \eqref{uexa}, with $k_0 = 5.0$ and $k_1 = 1.0$. $N=20$ intervals have
been adopted. Convergence order of the solution (b),(c),(d), convergence order of the divergence of the solution (e),(f), (g) and conditioning slope (h),(i), (l) for a fixed value of $\theta_2 = 10^{-3}$, $\theta_1 \in \{0.0010, 0.0015, 0.0020, 0.0025, \cdots, 0.990\}$ and $N \in \{20,40,80,160,320,640 \}$. In this test $\lambda = h^{-\alpha}$, where, in panels (b),(e) $\alpha = 2$, in (c),(f) $\alpha = 1.75$ and in (d),(g) $\alpha = 1.5$.}} 
\label{fig:compare_mix}
\end{figure}

We choose the following exact solution \begin{equation}
	u_{\rm exa} = \sin(k_0 x + k_1), \qquad {\rm in }\, [0,1]
	\label{uexa}
\end{equation}
with $k_0= 5$, $k_1 = 1$, to ensure that both Dirichlet and Neumann conditions are non-zero at the boundaries. In Fig.~\ref{fig:compare_mix}, we show the comparison between the exact solution (blue line) and the numerical solution (red asterisks) for the mixed problem in Eq.~\eqref{eq:system1M}, with Dirichlet boundary conditions in $a$, Neumann boundary conditions in $b$, and number of elements $N = 20$. In panels (b), (c), (d) we show the convergence order for the solution, in (e), (f), (g) the convergence order for the divergence of the solution and in (h), (i), (l) the slope of the conditioning number of the stiff matrix for a fixed value of $\theta_2 = 10^{-3}$, $\theta_1 \in \{0.0010, 0.0015, 0.0020, 0.0025, \cdots, 0.990\}$ and $N \in \{20,40,80,160,320,640 \}$. In this test the penalisation coefficient is $\lambda = h^{-\alpha}$, where, in panels (b),(e) $\alpha = 2$, in (c),(f) $\alpha = 1.75$ and in (d),(g) $\alpha = 1.5$.

\subsection{Numerical results in 2D}
In this section, we provide some numerical tests showing second-order accuracy in two dimensions. We test different smooth geometries: circle, flower and leaf-shaped domain, and, at the end, a saddle point domain. The rectangular region $R \subset \mathbb R^2$, and the different domains $\Omega\subset R$, are defined by level-set functions that we define below. In Figs.~\ref{fig_err_cA_circle}-\ref{fig_err_cA_hourglass}, we show, in panel (a), the shape of the domain that we are considering, in panel (b) the accuracy order of the solution for different values of penalisation coefficient $\lambda = h^{-\alpha}$ and snapping parameters $1/\lambda$ (defined in Algorithm~\ref{alg_snap}), with $\alpha \in \{1.5,1.75,2\}$. Furthermore, in panel (c), we show the slope of the conditioning number of the stiff matrix, and in panel (d), the accuracy order of the divergence of the solution. We observe that for $\alpha = 1.5$, the expected asymptotic convergence rate is achieved once $N$ becomes sufficiently large. The exact solution we choose is $u_{\rm exa} = \cos(2\pi x)\cos(2\pi y)$.

\subsubsection*{Circular domain} 
\begin{figure}[H]
	\centering
	\hfill
\begin{minipage}[b]
		{.49\textwidth}
		\centering
\begin{overpic}[abs,width=\textwidth,unit=1mm,scale=.25]{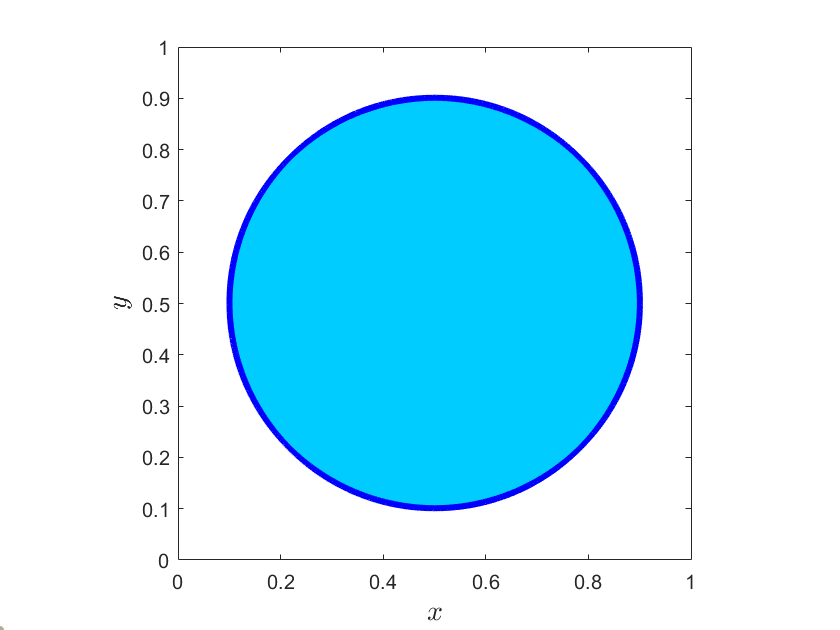}
    \put(-3,48){(a)}
\end{overpic}  
	\end{minipage}\hfill
 \begin{minipage}[b]
		{.49\textwidth}
		\centering
\begin{overpic}[abs,width=\textwidth,unit=1mm,scale=.25]{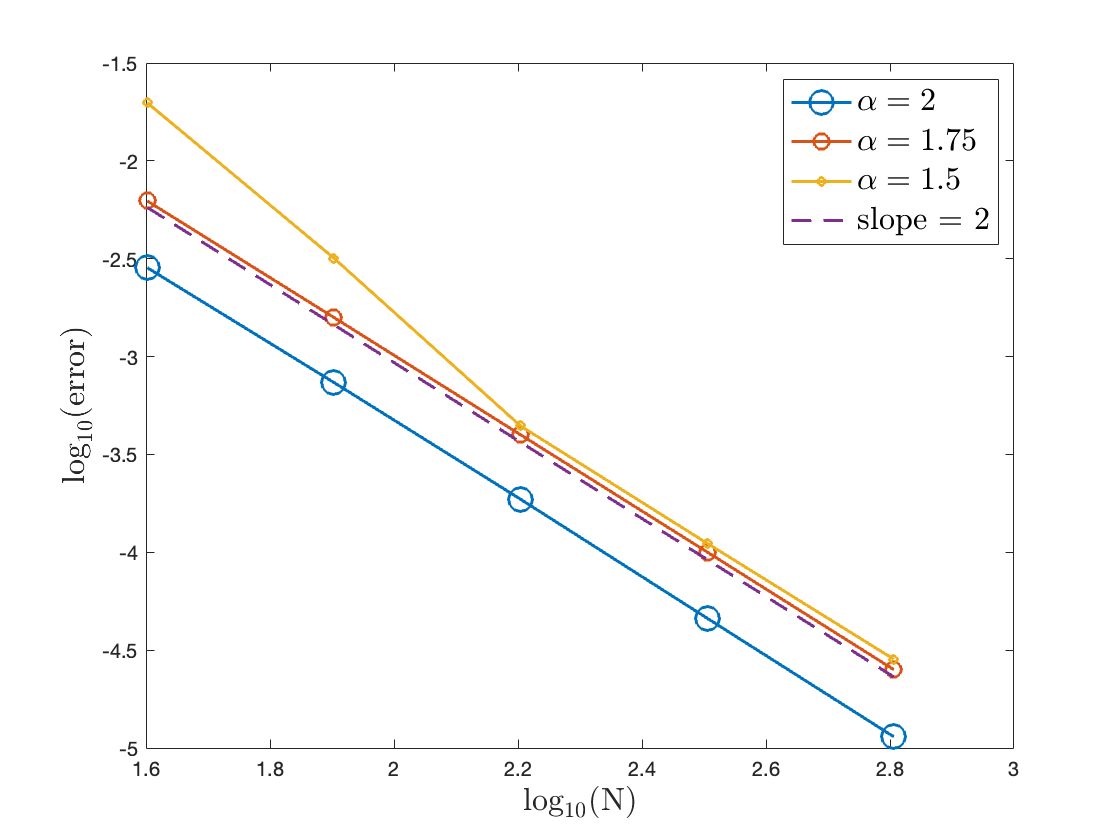}
\put(-2,48){(b)}\end{overpic}  
	\end{minipage}\hfill
	\begin{minipage}[b]
		{.49\textwidth}
		\centering
\begin{overpic}[abs,width=\textwidth,unit=1mm,scale=.25]{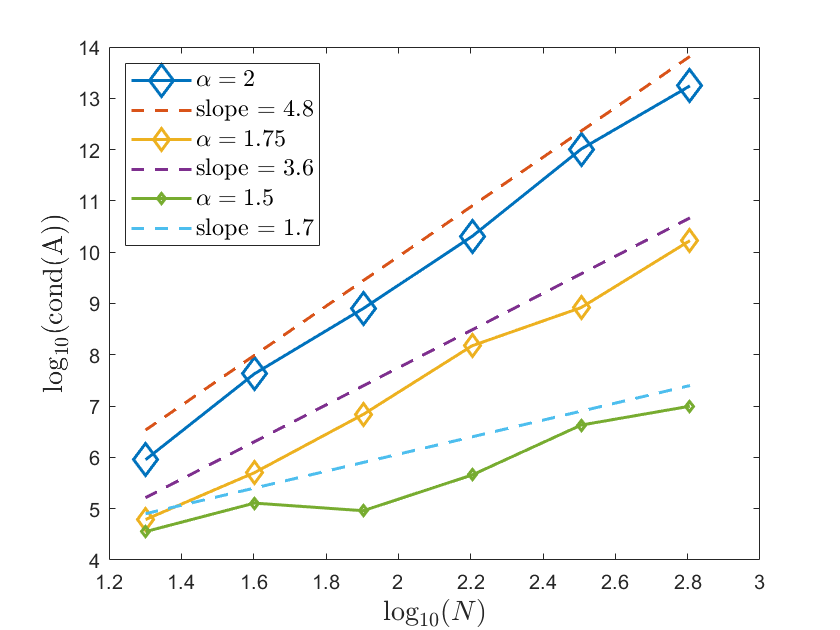}
\put(-2,48){(c)}\end{overpic}  
	\end{minipage}
	\begin{minipage}[b]
		{.49\textwidth}
		\centering
\begin{overpic}[abs,width=\textwidth,unit=1mm,scale=.25]{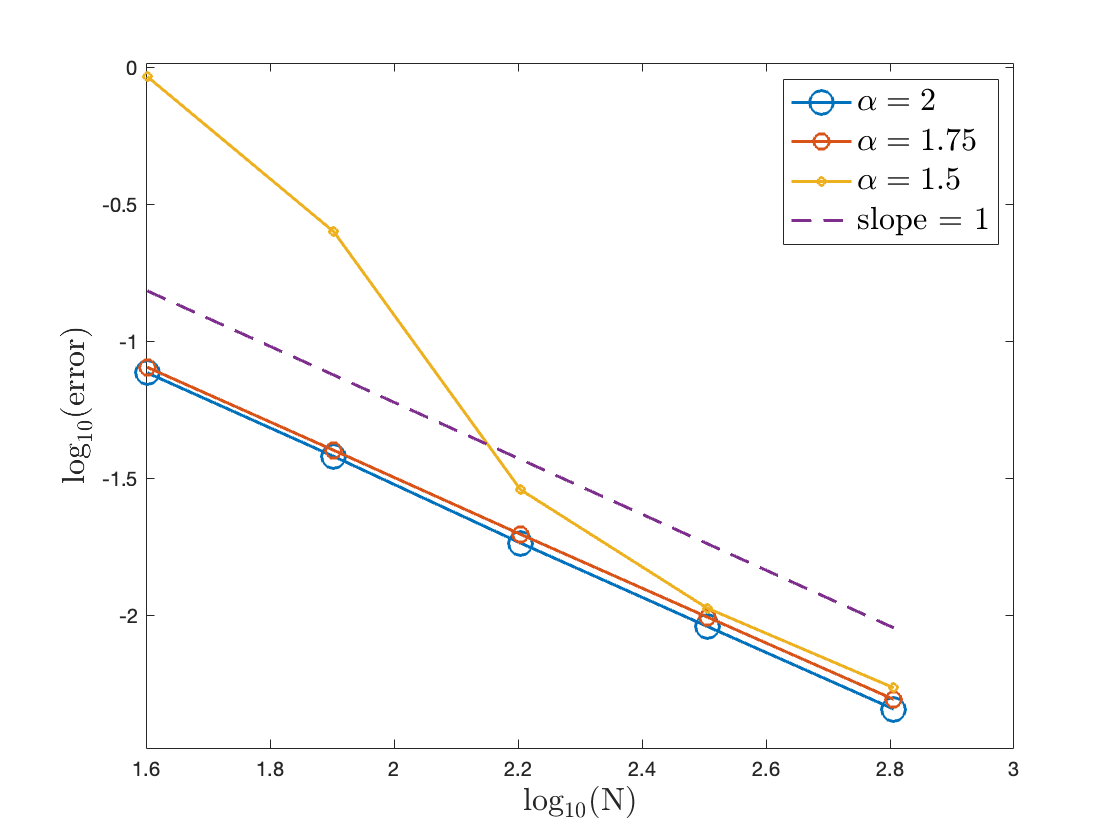}
\put(-2,48){(d)}\end{overpic}  
	\end{minipage} 
	\hspace*{\fill}
	\caption{\textit{Numerical results for the circular domain and Dirichlet boundary conditions: (a) shape of the domain; (b) error plot for different values of the snapping back to grid and penalisation term; (c) conditioning number for different values of the snapping back to grid and penalisation term; (d) error plot for the divergence of the solution for different values of the snapping back to grid and penalisation term.}}
	\label{fig_err_cA_circle}
\end{figure}

\begin{figure}
	\centering
	\hfill
\begin{minipage}[b]
		{.49\textwidth}
		\centering
\begin{overpic}[abs,width=\textwidth,unit=1mm,scale=.25]{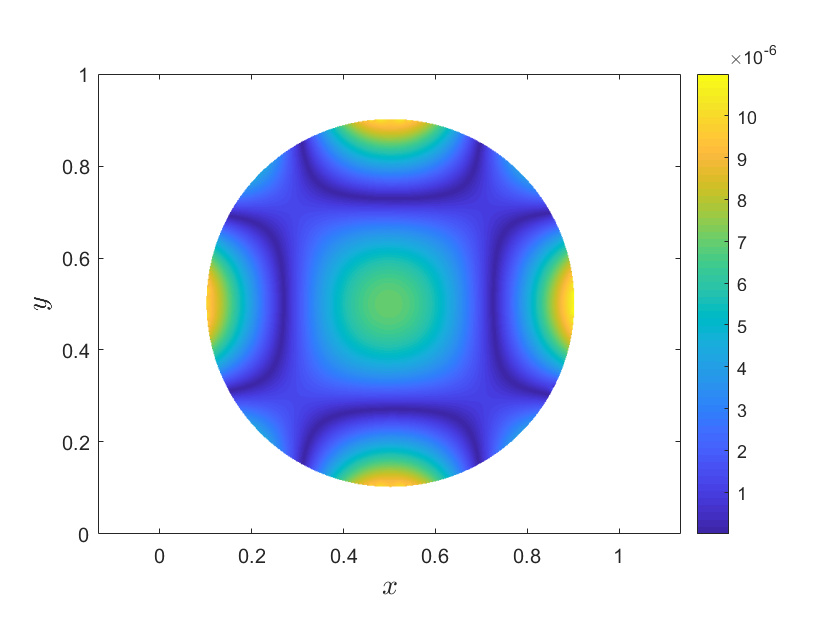}
    \put(-3,48){(a)}
\end{overpic}  
	\end{minipage}\hfill
 \begin{minipage}[b]
		{.49\textwidth}
		\centering
\begin{overpic}[abs,width=\textwidth,unit=1mm,scale=.25]{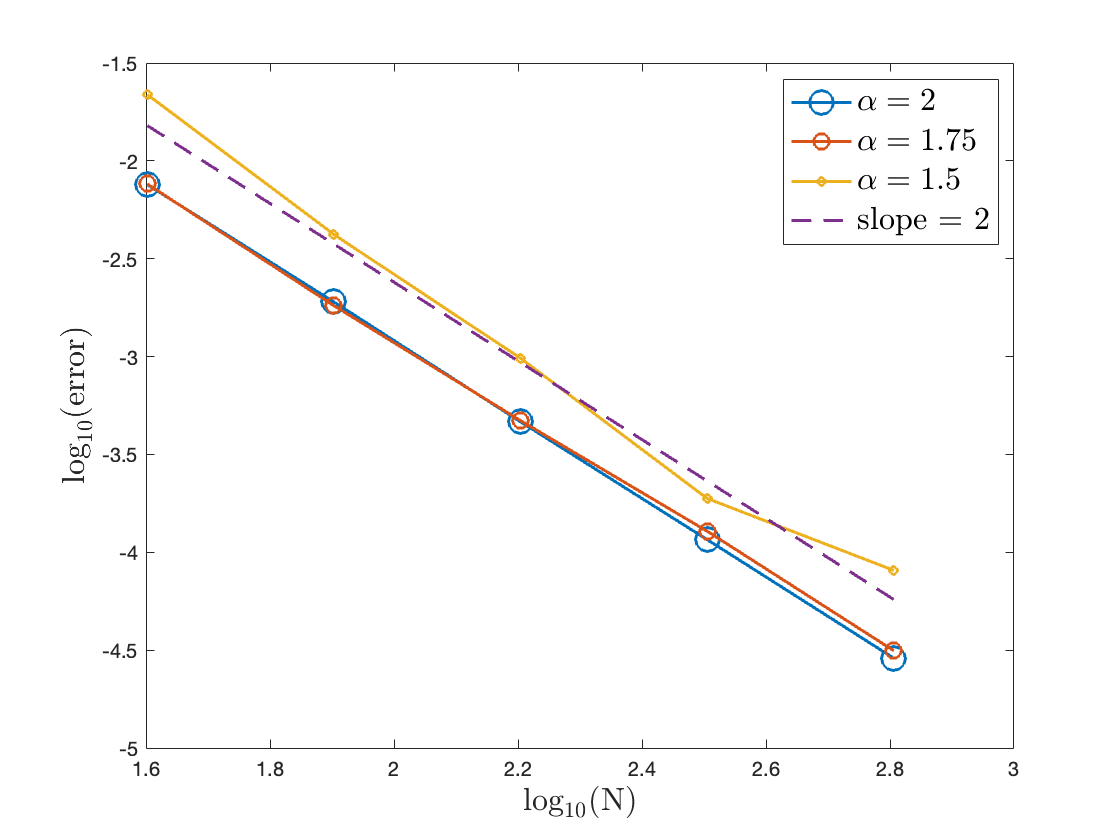}
\put(-2,48){(b)}\end{overpic}  
	\end{minipage}\hfill
	\begin{minipage}[b]
		{.49\textwidth}
		\centering
\begin{overpic}[abs,width=\textwidth,unit=1mm,scale=.25]{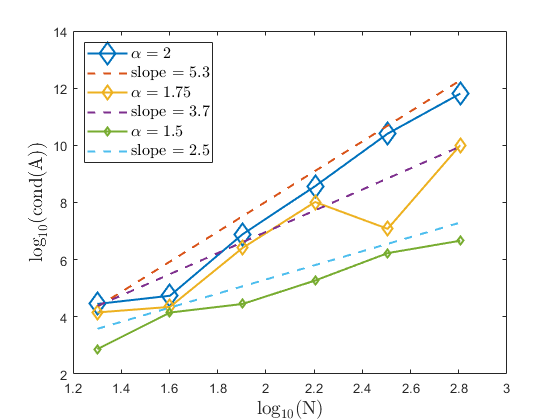}
\put(-2,48){(c)}\end{overpic}  
	\end{minipage}
	\begin{minipage}[b]
		{.49\textwidth}
		\centering
\begin{overpic}[abs,width=\textwidth,unit=1mm,scale=.25]{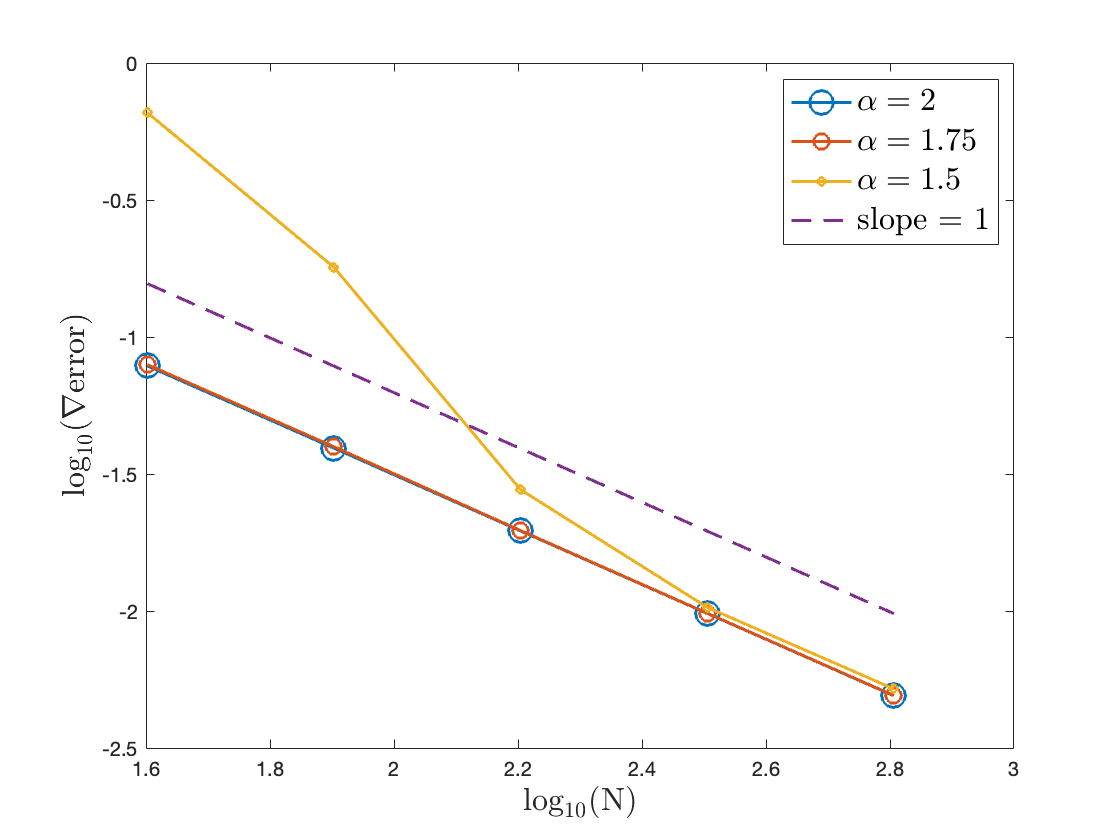}
\put(-2,48){(d)}\end{overpic}  
	\end{minipage} 
	\hspace*{\fill}
	\caption{\textit{Numerical results for the circular domain and mixed boundary conditions: (a) difference $|u_h - u_{\rm exa}|$; (b) error plot for different values of the snapping back to grid and penalisation term; (c) conditioning number for different values of the snapping back to grid and penalisation term; (d) error plot for the divergence of the solution for different values of the snapping back to grid and penalisation term.}}
	\label{fig_err_cA_circle_mixed}
\end{figure}

Let us start with a circular domain $\Omega\subset [0,1]^2$ centered at $(x_c,y_c) \in \Omega$. To prove that the order of accuracy does not depend on the symmetry of the domain, we calculate 10 different numerical solutions coming from different level-set functions. For this reason, we choose the following expression for the centre of the domain:
\[(x_c,y_c) = \left(0.5 + \widehat \epsilon_i\Delta x,0.5 + \widetilde \epsilon_j\Delta y\right), \quad i \in 1,\cdots,10,\] 
where $\widehat \epsilon_i, \widetilde \epsilon_j \in [0,1]$, are randomly chosen. In this case, the level-set function is $\phi=r-\sqrt{(x-x_c)^2+(y-y_c)^2}$, where $r = 0.4$ is the radius of the circle. In the end, we have 10 different variables for the quantities defined in Eqs.~(\ref{eq:error}-\ref{eq:Derror}), ${\rm error},\nabla{\rm error}$, one for each pair of $(x_c,y_c)$. 
In  Figs.~\ref{fig_err_cA_circle}-\ref{fig_err_cA_circle_mixed} 
we show the average of these 10 quantities.

We choose Dirichlet boundary conditions (i.e., $\Gamma_D = \Gamma$) in Fig.~\ref{fig_err_cA_circle}, and mixed boundary conditions in Fig.~\ref{fig_err_cA_circle_mixed}, i.e., $\Gamma_D = \Gamma \cap \{x\leq 0.5\}$ and $\Gamma_N = \Gamma \cap \{x > 0.5\}$.

\subsubsection*{Flower-shaped domain} 
\begin{figure}[H]
	\centering
	\hfill
\begin{minipage}[b]
		{.49\textwidth}
		\centering
\begin{overpic}[abs,width=\textwidth,unit=1mm,scale=.25]{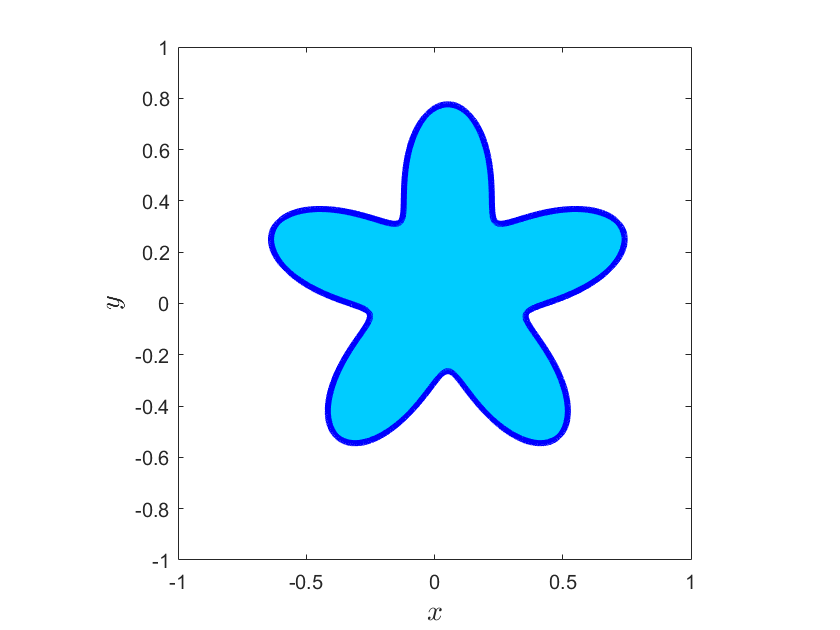}
\put(-3,48){(a)}\end{overpic}  
	\end{minipage}\hfill
 \begin{minipage}[b]
		{.49\textwidth}
		\centering
\begin{overpic}[abs,width=\textwidth,unit=1mm,scale=.25]{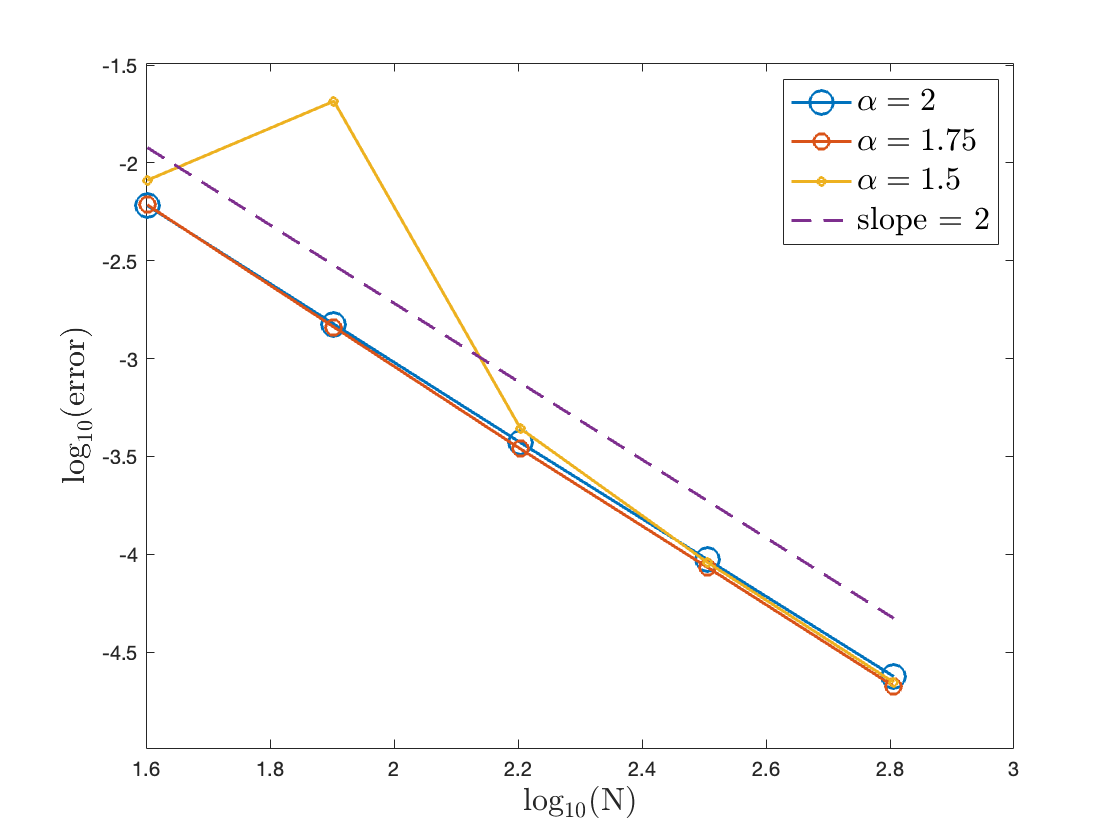}
\put(-2,48){(b)}\end{overpic}  
	\end{minipage}\hfill
	\begin{minipage}[b]
		{.49\textwidth}
		\centering
\begin{overpic}[abs,width=\textwidth,unit=1mm,scale=.25]{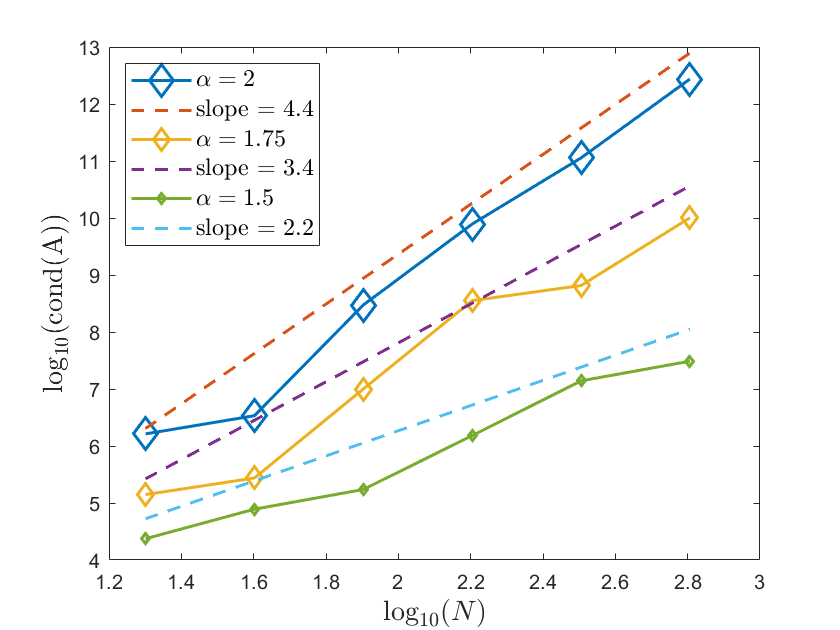}
\put(-2,48){(c)}\end{overpic}  
	\end{minipage}
	\begin{minipage}[b]
		{.49\textwidth}
		\centering
\begin{overpic}[abs,width=\textwidth,unit=1mm,scale=.25]{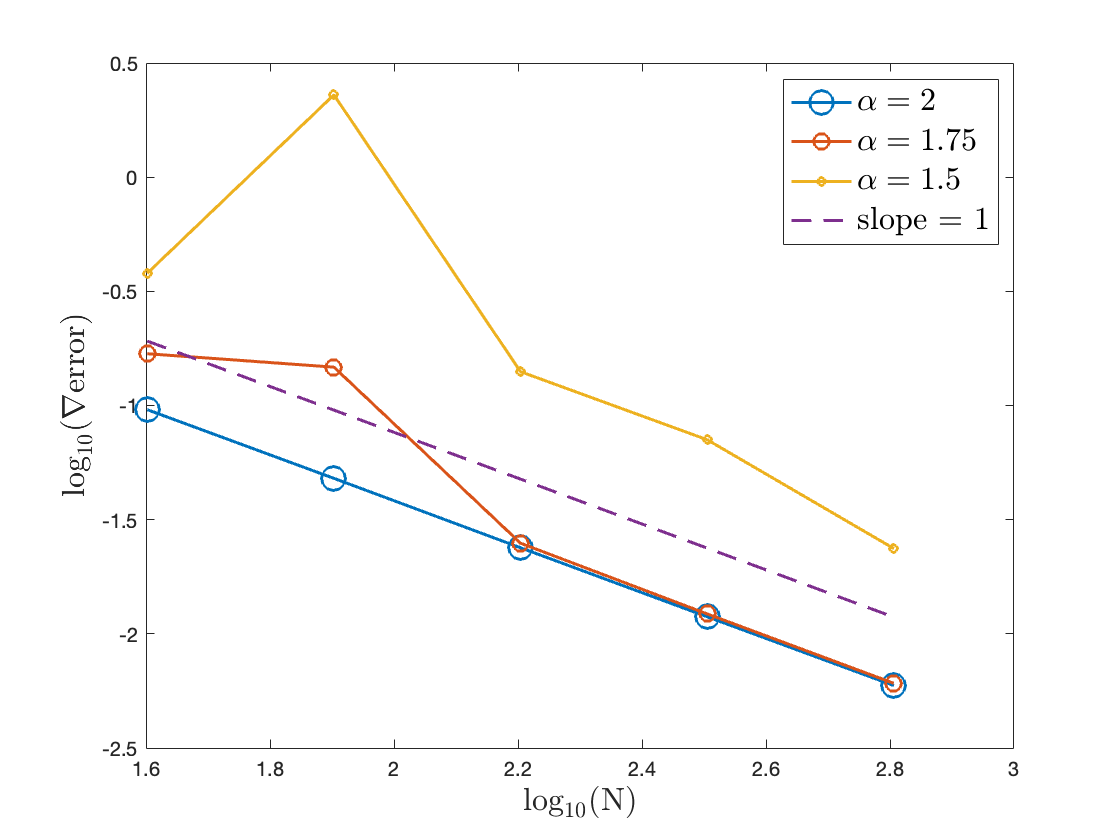}
\put(-2,48){(d)}\end{overpic}  
	\end{minipage} 
	\hspace*{\fill}
	\caption{\textit{Numerical results for the flower domain: (a) shape of the domain; (b) error plot for different values of the snapping back to grid and penalisation term; (c) conditioning number for different values of the snapping back to grid and penalisation term; (d) error plot for the divergence of the solution for different values of the snapping back to grid and penalisation term.}}
	\label{fig_err_cA_flower}
\end{figure}

Here we show a flower-shaped domain $\Omega \subset R = [-1,1]^2$, and Dirichlet boundary conditions, i.e., $\Gamma_D = \Gamma$. In this case, the level-set function is 
\begin{align*}
&X = x-0.03\sqrt 3, \quad  Y = y-0.04\sqrt 2, \quad R = \sqrt{X^2+Y^2} \\
&\phi = \frac{R - 0.52 - (Y^5 +5X^4Y-10X^2Y^3)}{5R^5},   
\end{align*}
and the results are shown in Fig.~\ref{fig_err_cA_flower}.

\subsubsection*{Leaf-shaped domain} 
\begin{figure}
	\centering
	\hfill
\begin{minipage}[b]
		{.49\textwidth}
		\centering
\begin{overpic}[abs,width=\textwidth,unit=1mm,scale=.25]{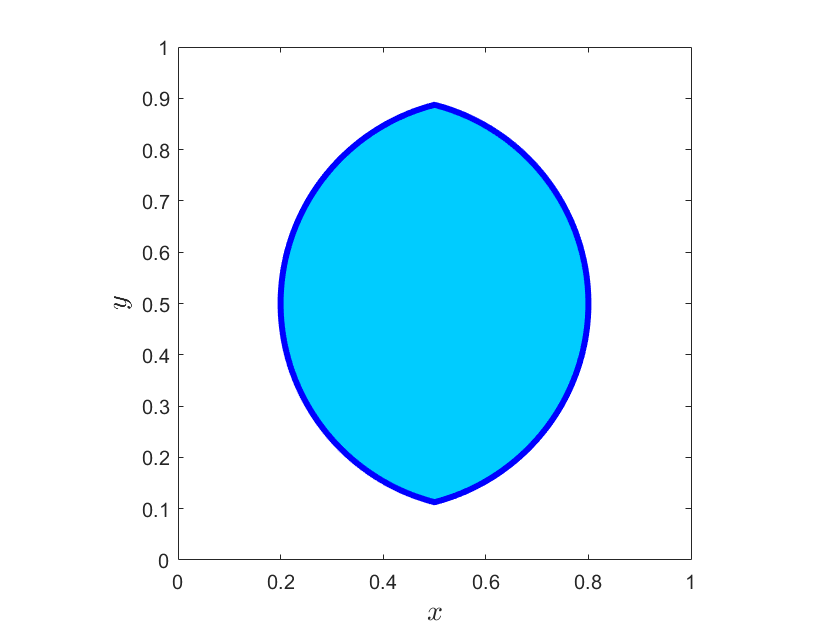}
\put(-3,48){(a)}\end{overpic}  
	\end{minipage}\hfill
 \begin{minipage}[b]
		{.49\textwidth}
		\centering
\begin{overpic}[abs,width=\textwidth,unit=1mm,scale=.25]{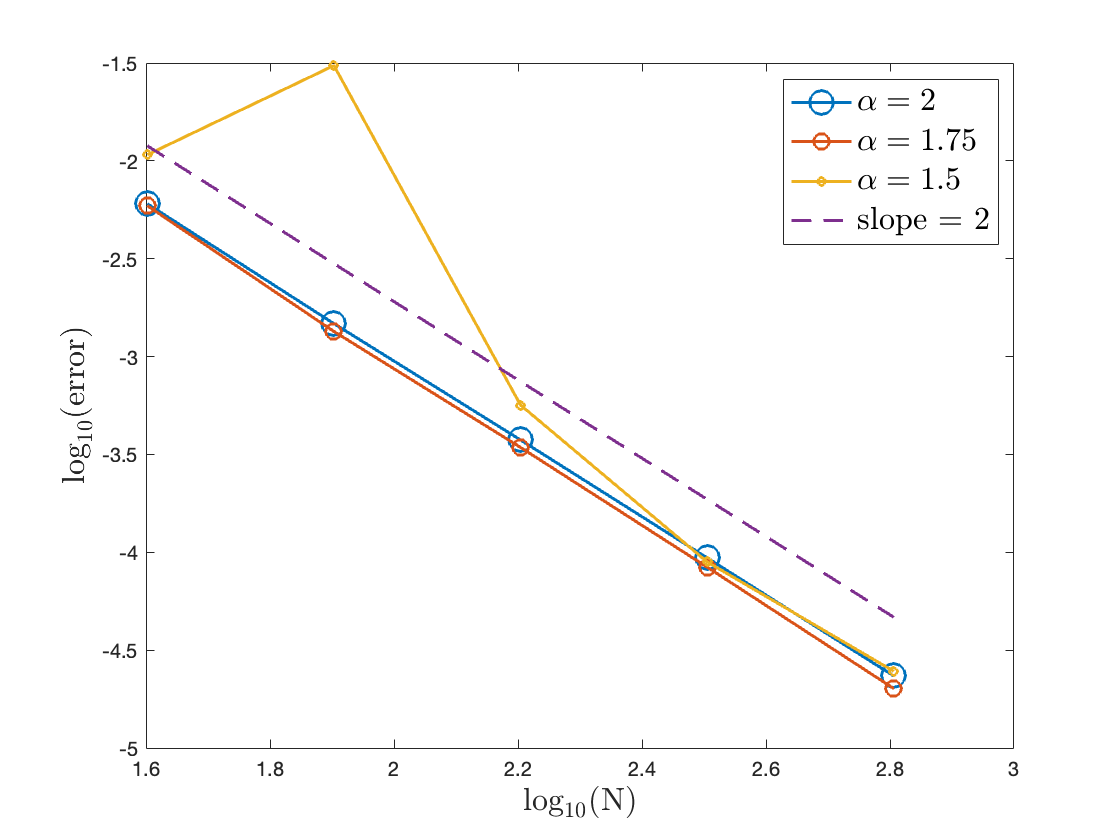}
\put(-2,48){(b)}\end{overpic}  
	\end{minipage}\hfill
	\begin{minipage}[b]
		{.49\textwidth}
		\centering
\begin{overpic}[abs,width=\textwidth,unit=1mm,scale=.25]{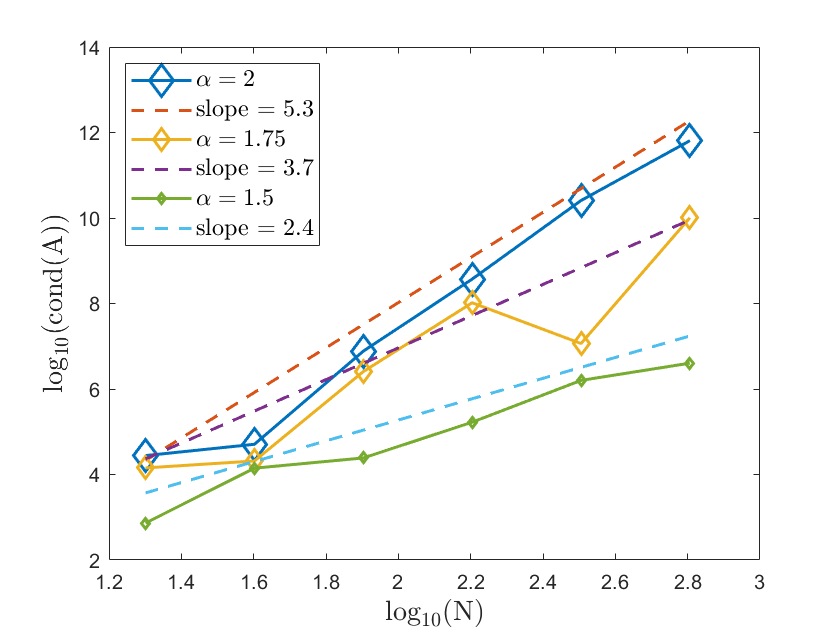}
\put(-2,48){(c)}\end{overpic}  
	\end{minipage}
	\begin{minipage}[b]
		{.49\textwidth}
		\centering
\begin{overpic}[abs,width=\textwidth,unit=1mm,scale=.25]{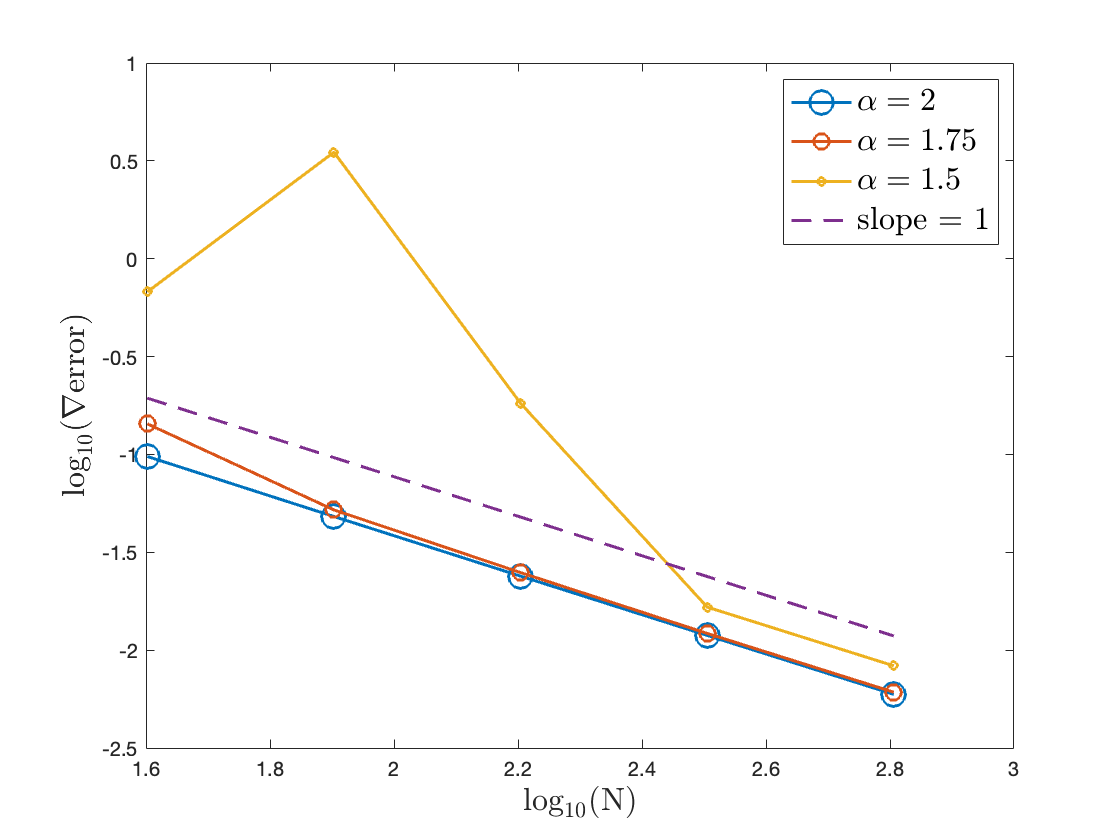}
\put(-2,48){(d)}\end{overpic}  
	\end{minipage} 
	\hspace*{\fill}
	\caption{\textit{Numerical results for the leaf domain: (a) shape of the domain; (b) error plot for different values of the snapping back to grid and penalisation term; (c) conditioning number for different values of the snapping back to grid and penalisation term; (d) error plot for the divergence of the solution for different values of the snapping back to grid and penalisation term.}}
	\label{fig_err_cA_leaf}
\end{figure}

Here we show a leaf domain $\Omega \subset R = [0,1]^2$ with mixed boundary conditions, such that $\Gamma_D = \partial \Omega \cap \{x<0.5\}$ and $\Gamma_N = \partial \Omega \cap \{x\geq 0.5\}$. In this case, the level-set function is 
\begin{align*}
&x_1 = 0.4,\quad y_1 = 0.5, \quad R_1 = \sqrt{({ x}-x_1)^2+(y-y_1)^2} \\
&x_2 = 0.6,\quad y_2 = 0.5, \quad R_2 = \sqrt{(x-x_2)^2+(y-y_2)^2} \\         
&r_0 = 0.4,\quad  \phi_1 = R_1-r_0, \quad \phi_2 = R_2-r_0, \quad \phi = \max\{\phi_1,\phi_2\},
\end{align*}
and the results are in Fig.~\ref{fig_err_cA_leaf}.

\subsubsection*{Hourglass domain} 
\begin{figure}
	\centering
	\hfill
\begin{minipage}[b]
		{.49\textwidth}
		\centering
\begin{overpic}[abs,width=\textwidth,unit=1mm,scale=.25]{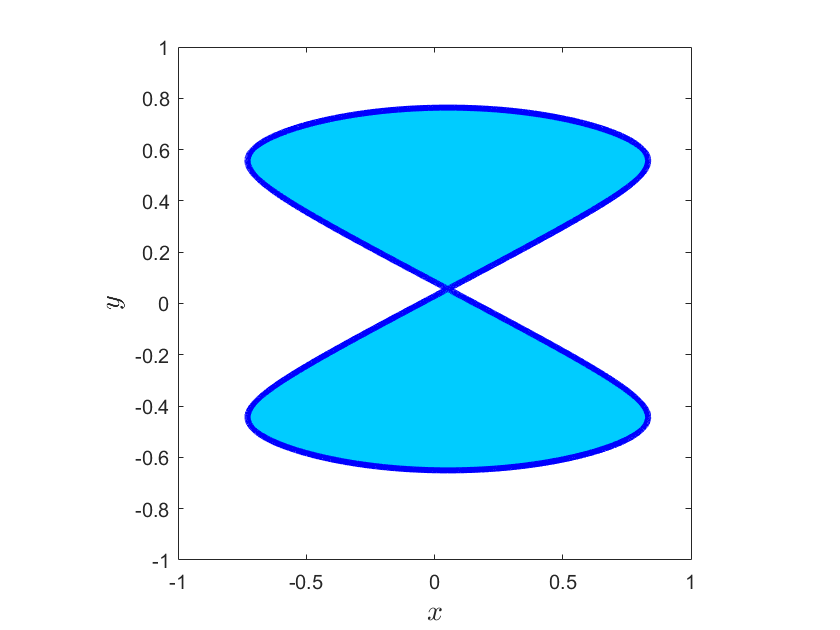}
\put(-3,48){(a)}\end{overpic}  
	\end{minipage}\hfill
 \begin{minipage}[b]
		{.49\textwidth}
		\centering
\begin{overpic}[abs,width=\textwidth,unit=1mm,scale=.25]{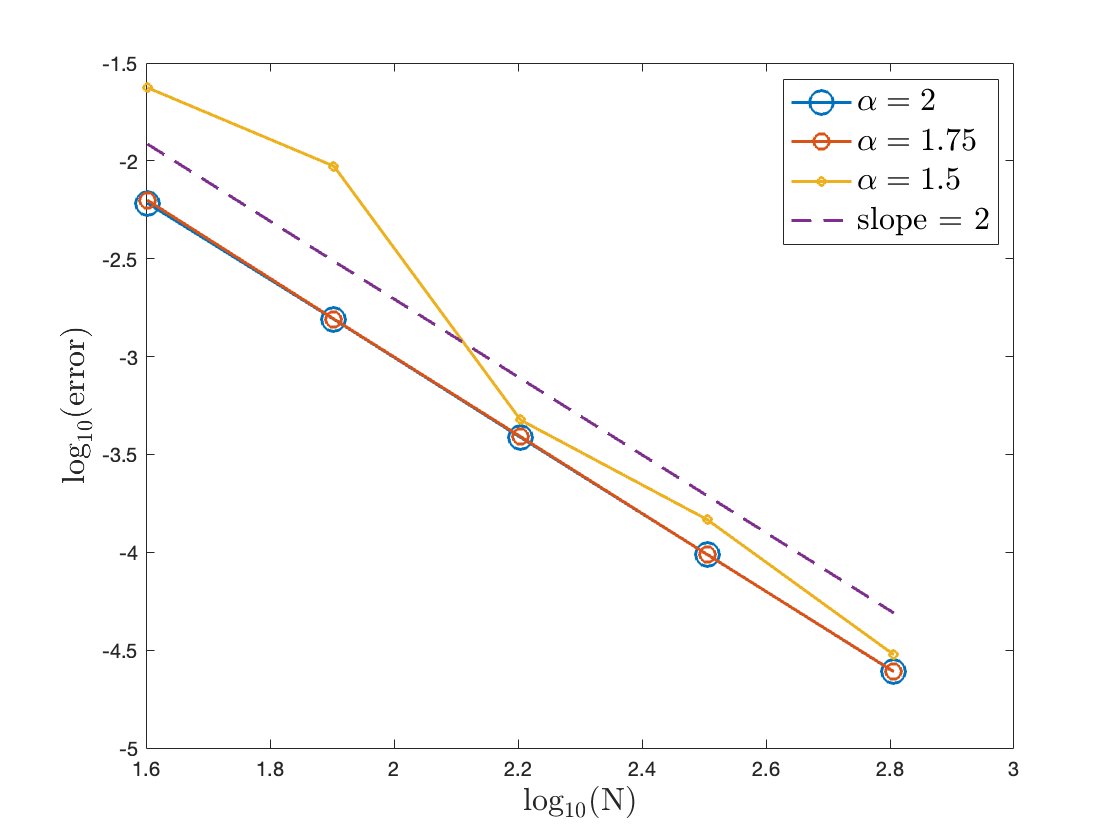}
\put(-2,48){(b)}\end{overpic}  
	\end{minipage}\hfill
	\begin{minipage}[b]
		{.49\textwidth}
		\centering
\begin{overpic}[abs,width=\textwidth,unit=1mm,scale=.25]{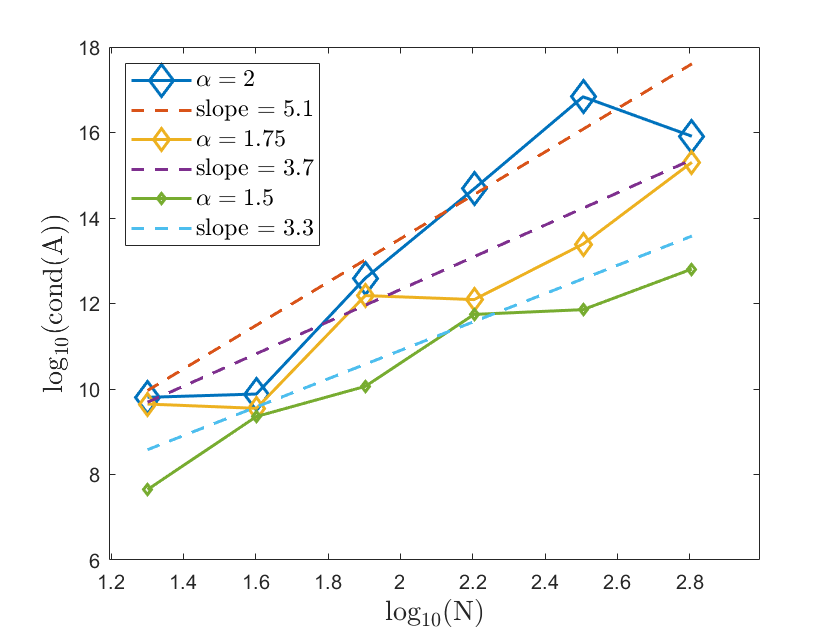}
\put(-2,48){(c)}\end{overpic}  
	\end{minipage}
	\begin{minipage}[b]
		{.49\textwidth}
		\centering
\begin{overpic}[abs,width=\textwidth,unit=1mm,scale=.25]{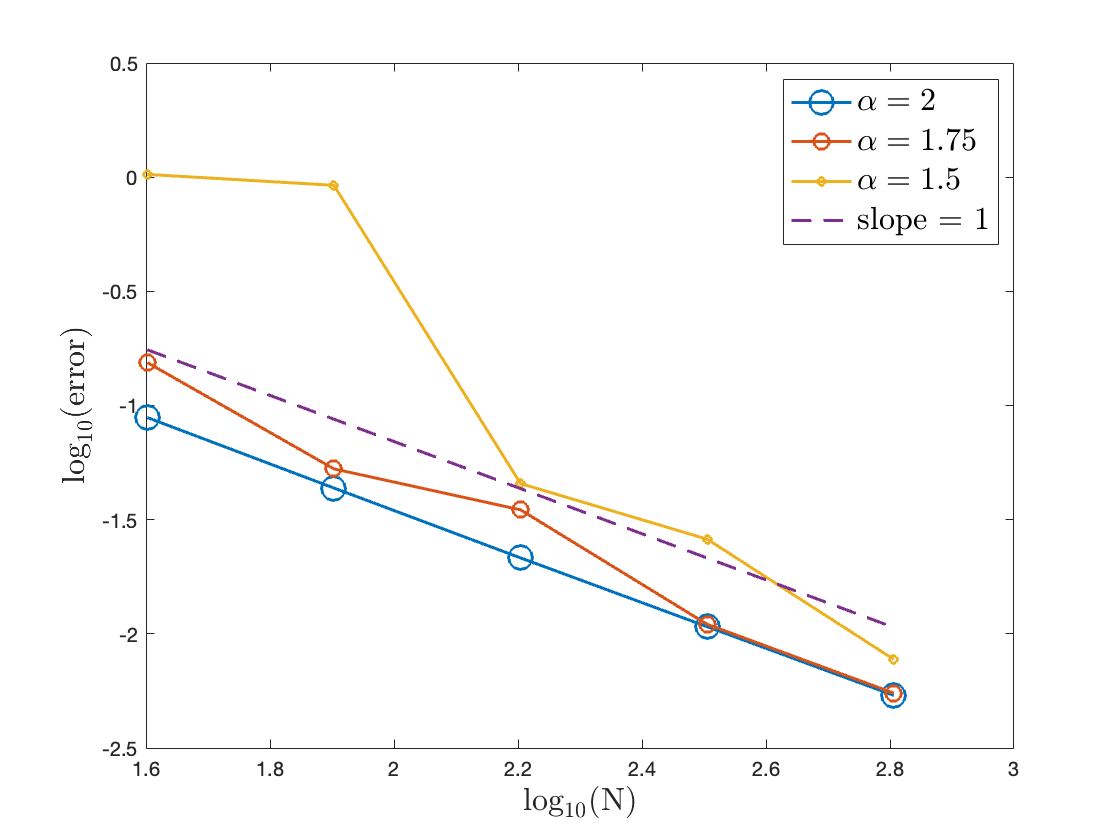}
\put(-2,48){(d)}\end{overpic}  
	\end{minipage} 
	\hspace*{\fill}
	\caption{\textit{Numerical results for the hourglass domain: (a) shape of the domain; (b) error plot for different values of the snapping back to grid and penalisation term; (c) conditioning number for different values of the snapping back to grid and penalisation term; (d) error plot for the divergence of the solution for different values of the snapping back to grid and penalisation term.}}
	\label{fig_err_cA_hourglass}
\end{figure}

\begin{figure}
	\centering
	\hfill
\begin{minipage}[b]
		{.49\textwidth}
		\centering
\begin{overpic}[abs,width=\textwidth,unit=1mm,scale=.25]{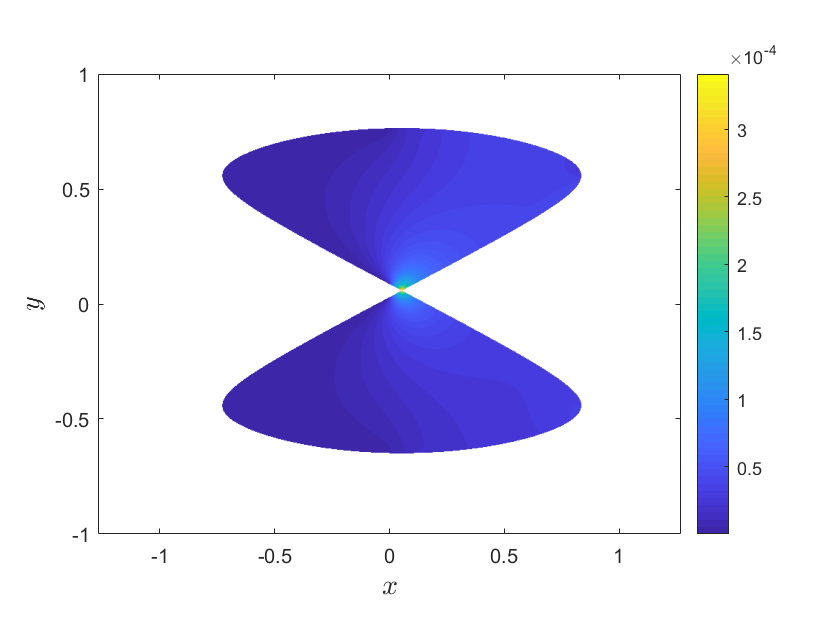}
    \put(-3,48){(a)}
\end{overpic}  
	\end{minipage}\hfill
 \begin{minipage}[b]
		{.49\textwidth}
		\centering
\begin{overpic}[abs,width=\textwidth,unit=1mm,scale=.25]{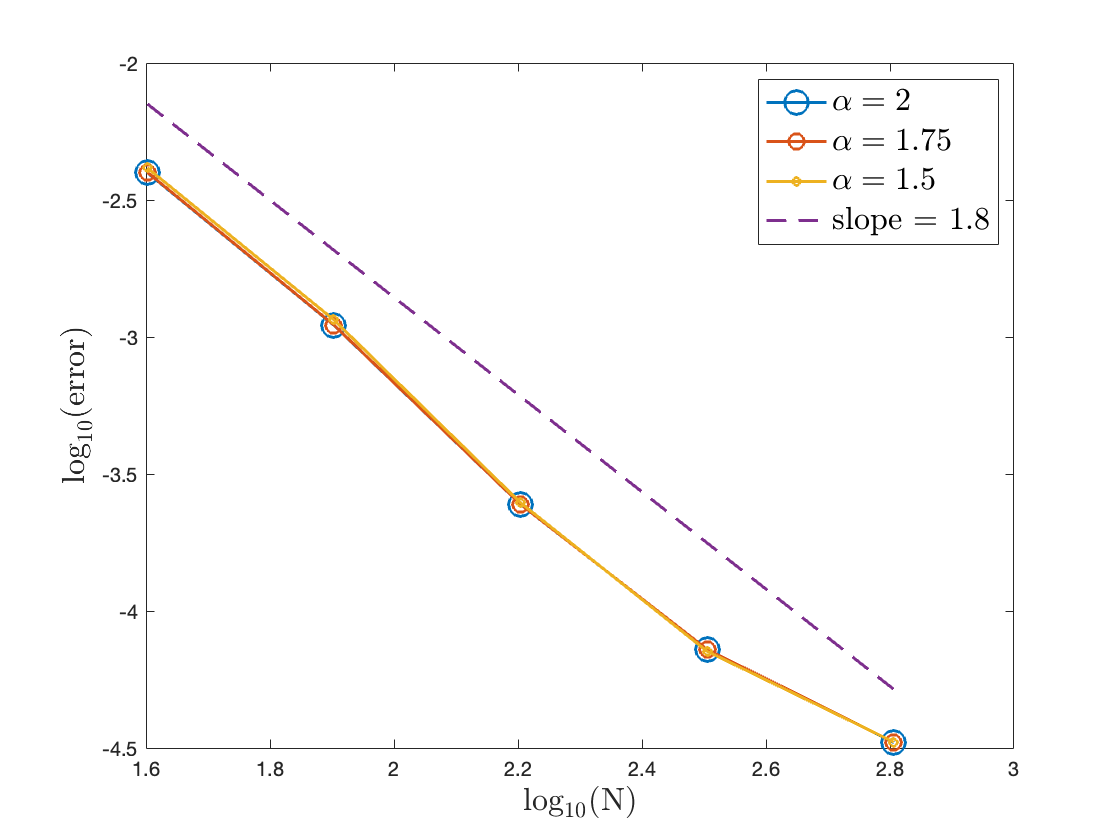}
\put(-2,48){(b)}\end{overpic}  
	\end{minipage}\hfill
	\begin{minipage}[b]
		{.49\textwidth}
		\centering
\begin{overpic}[abs,width=\textwidth,unit=1mm,scale=.25]{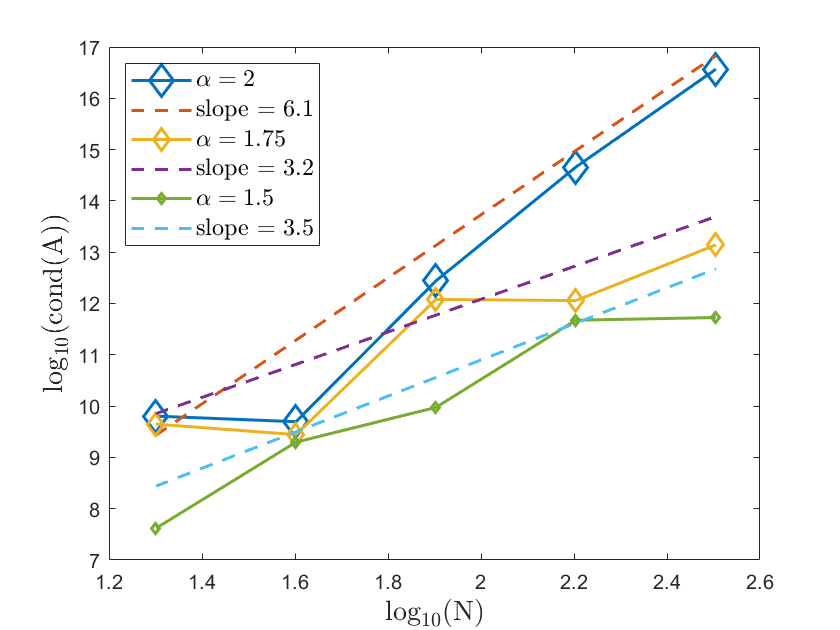}
\put(-2,48){(c)}\end{overpic}  
	\end{minipage}
	\begin{minipage}[b]
		{.49\textwidth}
		\centering
\begin{overpic}[abs,width=\textwidth,unit=1mm,scale=.25]{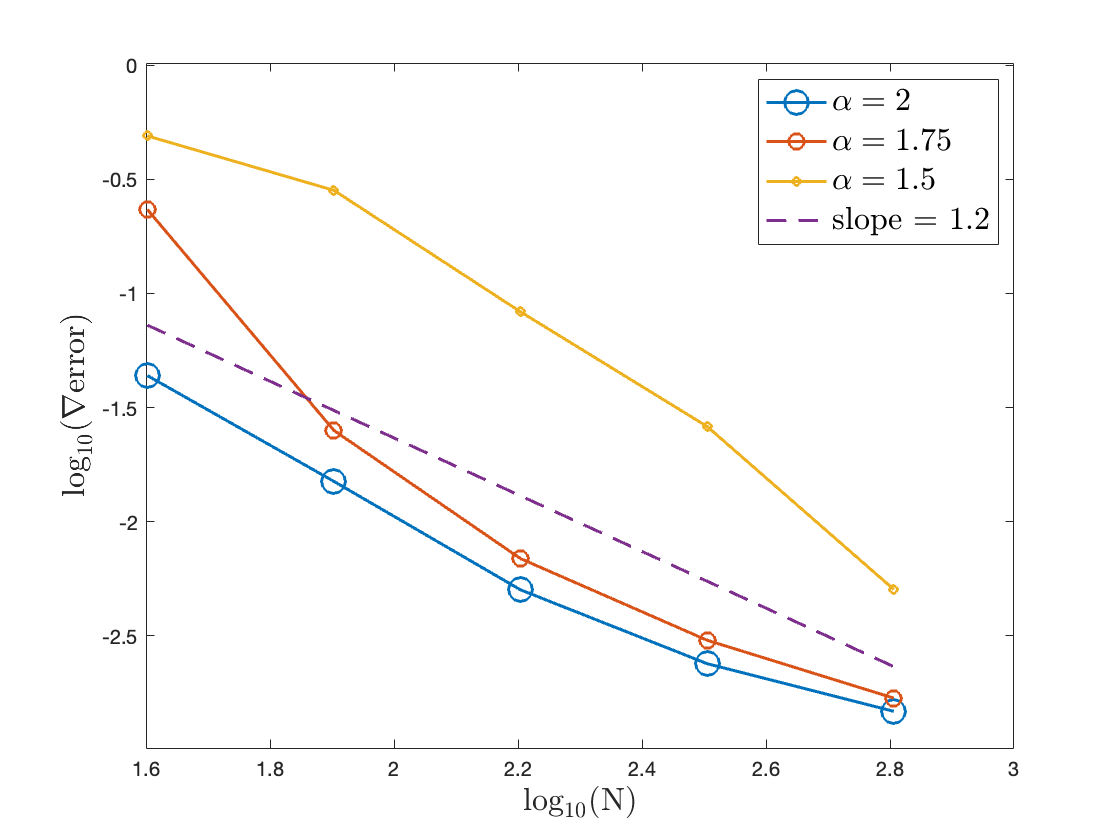}
\put(-2,48){(d)}\end{overpic}  
	\end{minipage} 
	\hspace*{\fill}
	\caption{\textit{Numerical results for the hourglass domain and mixed boundary conditions: (a) difference $|u_h - u_{\rm exa}|$; (b) error plot for different values of the snapping back to grid and penalisation term; (c) conditioning number for different values of the snapping back to grid and penalisation term; (d) error plot for the divergence of the solution for different values of the snapping back to grid and penalisation term.}}
	\label{fig_err_cA_hourglass_mixed}
\end{figure}

Here we consider a domain $\Omega \subset R = [-1,1]^2$, with a saddle point, whose level-set function is
\begin{align*}
    X = x-0.03\sqrt{3}, \qquad Y = y-0.04\sqrt{2}, \qquad \phi(x,y) = 256\,Y^4-16\,X^4-128\,Y^2 + 36\,X^2.
\end{align*}
We choose Dirichlet boundary conditions (i.e., $\Gamma_D = \Gamma$) in Fig.~\ref{fig_err_cA_hourglass}, and mixed boundary conditions in Fig.~\ref{fig_err_cA_hourglass_mixed}, i.e., $\Gamma_D = \Gamma \cap \{x\leq 0\}$ and $\Gamma_N = \Gamma \cap \{x > 0\}$. 

\section{Conclusions}

In the present work, we have introduced a novel nodal ghost discretisation posed on regular square meshes. In particular, our derivation of the numerical scheme here presented found its root in the variational formulation of the Poisson equation. Such formulation, often known also as weak formulation, is at the heart of the finite element method. 

Resorting to a variational formulation allows us to preserve the symmetry and the coercivity of the continuous problem in the discrete setting. The conservation of such properties at a discrete level represents a substantial difference and a great improvement with respect to other nodal ghost numerical schemes such as the Coco-Russo method \cite{Coco2013, Coco2018}.

There are many differences between the method here proposed and a standard finite element discretisation. One particular difference, between our scheme and a standard body-fitted FEM, is that we described the geometry of the domain $\Omega$ via a level-set function. Such an approach is also adopted in CutFEM, \cite{HansboHansbo, Burman2015, Burman2012, Burman2022}, which differs from our scheme for the following reasons.
First of all, in the CutFEM in order to compute the boundary integrals, they introduce a sub-triangulation of cut elements. Such a process requires additional computational work. We avoid the construction of a sub-triangulation via the quadrature technique presented in Section~\ref{sec:comp}.
Secondly, we escape the so-called \textit{small cut} issue, which might result in losing the coercivity of the bilinear form, by perturbing mesh geometry.
In particular, the perturbation procedure we adopted is well known to the engineering community by the name of snapping back to grid.

Adopting standard techniques, as in the FEM and the CutFEM, we can prove optimal converge rate for our scheme with respect to the $L^2$ norm.
To the best of the authors' knowledge up until now, there has not been a rigorous study of a priori convergence rates when the snapping back to grid technique is adopted.
For this reason, we have presented a novel approach to proving a priori convergence rates for our scheme even when the snapping back to grid technique is used to address the \textit{small cuts} issue.

Lastly, we have presented many numerical experiments that validate the $L^2$ convergence rate theoretically proven. Furthermore, we would like to point out that the numerical experiments shown in Section \ref{sec:num} suggest that the error in energy norm is super-converging with respect to the analysis here presented.

In future, we plan to further investigate the scheme here proposed. In particular, we would like to study in greater detail the eigenvalue problem \eqref{eq:eig} and the role it plays when determining the optimal value of the penalisation constant. We are also interested in developing fast solvers for our discretization based on a multigrid approach, similar to \cite{Coco2013}.

\section*{Acknowledgments}
{The work of C.A. and G.R. has been supported by the Spoke 10 Future AI Research (FAIR) of the Italian Research Center funded by the Ministry of University and Research as part of the National Recovery and Resilience Plan (PNRR). \\ 
C.A. and G.R. have also been supported also by Italian Ministerial grant PRIN 2022 “Efficient numerical schemes and optimal control methods for time-dependent partial differential equations”, No. 2022N9BM3N - Finanziato dall’Unione europea - Next Generation EU – CUP: E53D23005830006.} \\
{The work of C.A. and G.R. has also been supported by the Italian Ministerial grant PRIN 2022 PNRR “FIN4GEO: Forward and Inverse Numerical Modeling of hydrothermal systems in volcanic regions with application to geothermal energy exploitation”, No. P2022BNB97 - Finanziato dall’Unione europea - Next Generation EU – CUP: E53D23017960001. } \\
{The authors are members of the Gruppo Nazionale Calcolo Scientifico-Istituto Nazionale di Alta Matematica (GNCS-INdAM).}

\appendix
\section{Discretization of the two-dimensional problem}
\label{app:sec_2D}
In this section we derive in detail the computation to obtain a space discretization in two-dimensions. We start considering the Poisson equation with Dirichlet boundary conditions in $\Omega\subseteq R\subset \RR^2$, where $R$ is a rectangular region:
\begin{equation}
	-\Delta u = f \quad {\rm in}\>  \Omega, \quad u = g_D \> {\rm on}\>  \partial \Omega
	\label{eq2D0_A}	
\end{equation}
Multiplying by a test function $v\in H^1(\Omega)$, integrating in $\Omega$ and then integrating by parts, one obtains: 
\begin{equation}
	\int_{\Omega}\nabla u\cdot \nabla v \, d\Omega = \int_\Omega f v\, d\Omega + \int_{\partial
\Omega}v\pad{u}{n}\, ds.
\end{equation}
Following the same strategy that we applied in the one-dimensional case, we penalise the distance from the prescribed value of the Dirichlet boundary condition by minimizing the following functional 
\begin{equation}
	J_{D\lambda}[u] = \int_\Omega \left[\frac12 (\nabla u)^2 - f u\right]\, d\Omega + 
	\frac{\lambda}2\int_{\partial \Omega} (u-g_D)^2\, ds
	\label{JD2L[u]_A}
\end{equation}
By imposing the first variation of the functional equal to zero we obtain the following penalised variational formulation of the problem: 
\begin{equation}
	\int_\Omega\nabla u\cdot \nabla v\, d\Omega +  \int_{\partial \Omega} \left[\lambda\cdot (u-g_D)-\pad{u}{n}\right]v \, ds = \int_\Omega f v \, d\Omega
	\label{eq:weak5_A}
\end{equation}

The discretization of the problem is obtained by discretizing the rectangular region $R\subset \RR^2$, containing the domain $\Omega$, by a regular rectangular grid (see Figure~\ref{fig:Domain2D}).

The discrete space $V_h$ defined in \eqref{eq:Vh} is given by the piecewise bilinear functions which are continuous in $R$.
As a basis, we choose the following functions:
\begin{equation}
    {\varphi_{i}(x,y) = \max\left\{
        \left(1-\frac{|x-x_i|}{\Delta x}\right),0
    \right\}\cdot \max\left\{
        \left(1-\frac{|y-y_i|}{\Delta x}\right),0
    \right\}}
    \label{eq:V_h2_A}
\end{equation}
here we denote by $i = (i_1,i_2)$ an index that identifies a node on the grid. 

The generic element $u_h\in V_h$ will have the following representation:
\begin{equation}
    u_h(x,y) = \sum_{i\in\nodes}u_i\vphi_i(x,y).
    \label{eq:u_h2_A}
\end{equation}
The discretization of the problem \eqref{eq:weak5_A} is obtained by replacing 
$u$ and $v$ by $u_h$ and $v_h$, both in $V_h$, and evaluating the integral instead of over $\Omega$ over its polygonal approximation $\Omega_h$. Since we are approximating $\Omega$ with $\Omega_h$, we are also approximating $\Gamma$ by its polygonal approximation $\Gamma_{h}$. This leads to the following discrete variational problem
\begin{equation}
	\int_{\Omega_h}\nabla u_h\cdot \nabla v_h\, d{\Omega_h} +  \int_{\Gamma_h} \left[\lambda\cdot (u_h-g_D)-\pad{u_h}{n}\right]v_h \, ds = 
	\int_{\Omega_h} f v_h \, d{\Omega_h}, \quad \forall v_h\in V_h
	\label{eq:weak2h_A}
\end{equation}

When applied to the subspace $V_h$, and using the standard basis, the above equations become
\begin{align}
    &\sum_{j\in \mathcal N} u_j(\nabla\varphi_i,\nabla\varphi_j)_{L^2(\Omega_h)}  - \sum_{j\in \mathcal N} u_j\left(\pad{\varphi_i}{n},\varphi_j\right)_{L^2(\Gamma_{D,h})} + \lambda\sum_{j\in \mathcal N} u_j(\varphi_j,\varphi_i)_{L^2(\Gamma_{D,h})} \nonumber \\ = & \sum_{j\in \mathcal N} f_j (\varphi_j,\varphi_i)_{L^2(\Omega_h)} + \lambda\sum_{j\in \mathcal N} {g_D}_j(\varphi_j,\varphi_i)_{L^2(\Gamma_{D,h})}.
\label{eq:weak_penal_discrete_A}    
\end{align}
where the quantities $f_i$, $g^D_i$ denote the node values of the source data and of the Dirichlet boundary conditions.

Here we rewrite \eqref{eq:weak2h_A}, applying the Nitsche technique, to make the equation symmetric
\begin{align} \nonumber
	\int_{\Omega_h}\nabla u_h\cdot \nabla v_h\, d{\Omega_h} +  \int_{\Gamma_h} \lambda\cdot (u_h-g_D)  v_h \, ds - \int_{\Gamma_h} \pad{u_h}{n} v_h \, ds - \int_{\Gamma_h} u_h \pad{v_h}{n} \, ds \\ = 
	\int_{\Omega_h} f v_h \, d{\Omega_h} - \int_{\Gamma_h} g_D \pad{v_h}{n} \, ds  , \quad \forall v_h\in V_h
     \label{eq:weakdiscrte2h_A}
\end{align}
as follows
\begin{align} \nonumber
    &\sum_{j\in \mathcal N} u_j\left(\underbrace{(\nabla\varphi_j,\nabla\varphi_i)_{L^2(\Omega_h)}}_{\rm S} - \underbrace{\left(\left(\pad{\varphi_j}{n},\varphi_i\right)_{L^2(\Gamma_{D,h})} + \left(\varphi_j,\pad{\varphi_i}{n}\right)_{L^2(\Gamma_{D,h})} \right)}_{\rm S_T} + \lambda\underbrace{(\varphi_j,\varphi_i)_{L^2(\Gamma_{D,h})}}_{\rm P_{\Gamma}} \right) \\ = & \sum_{j\in \mathcal N} f_j \underbrace{(\varphi_j,\varphi_i)_{L^2(\Omega_h)}}_{\rm M} + \sum_{j\in \mathcal N} {g_D}_j\left(-\underbrace{\left(\varphi_j,\pad{\varphi_i}{n}\right)_{L^2(\Gamma_{D,h})}}_{\rm D_{\Gamma}} + \underbrace{\lambda(\varphi_j,\varphi_i)_{L^2(\Gamma_{D,h})}}_{\rm P_{\Gamma}} \right)
\end{align}

Making use of representation \eqref{eq:u_h2_A} and using the functions of the basis as test functions we obtain the following system of equations for $u_i$, $i\in\nodes$:
\begin{equation}
	\sum_{j\in\nodes} a^D_{ij} u_j = f^D_i, \quad i\in\nodes
	\label{eq:system2D_A}
\end{equation}
where 
\begin{equation}
\begin{array}{c}
	\displaystyle	a^D_{ij}  = (\nabla\varphi_i,\nabla\varphi_j)_{L^2(\Omega_h)} + \lambda (\varphi_i,\varphi_j)_{L^2(\Gamma_{D,h})} 
	- (n\cdot\nabla\vphi_j,\vphi_i)_{L^2(\Gamma_{D,h})} \\ - (\vphi_j,n\cdot\nabla\vphi_i)_{L^2(\Gamma_{D,h})},\\[4mm]
	\displaystyle	f^D_i  = \lambda (g_D,\vphi_i)_{L^2(\Gamma_{D,h})}  - (g_D,n\cdot\nabla\vphi_i)_{L^2(\Gamma_{D,h})}  + (f,\vphi_i)_{L^2(\Omega_h)}.
\end{array}
\label{coeffD2_A}
\end{equation}

Now we focus on the Poisson equation with Neumann boundary conditions
\begin{equation}
	-\Delta u = f \quad {\rm in}\>  \Omega, \qquad \pad{u}{n} = g_N \> {\rm on}\>  \Gamma_N,
	\label{eq2N_A}	
\end{equation}
with $\Gamma_N = \partial \Omega$. 

It is sufficient to multiply the equation \eqref{eq2N_A} by a test function $v\in H^1(\Omega)$, integrate over $\Omega$ and then integrate by parts to obtain
\begin{align}
    \int_\Omega\nabla u\cdot \nabla v\, d\Omega  = \int_\Omega f v \, d\Omega 
	+  \int_{\Gamma_N} g_N v \, ds 
	\label{eq:weak2N}
\end{align}
Here we need to be careful that a compatibility condition has to be verified by $f\in L^2(\Omega)$ and  we need to choose a unique solution in the kernel of the gradient, imposing a null-average condition, i.e.,
\begin{equation}
    \int_\Omega f \,d\Omega = \int_{\Gamma_N} g_N\,ds ,\qquad \int_\Omega u \,d\Omega= 0.
\end{equation}

Here we discretize \eqref{eq:weak2N} using the standard basis of the discrete space $V_h$, that becomes
\begin{align}
    \sum_{j\in \mathcal N} u_j(\nabla\varphi_i,\nabla\varphi_j)_{L^2(\Omega_h)} 
    = \sum_{j\in \mathcal N} f_j (\varphi_j,\varphi_i)_{L^2(\Omega_h)} + \sum_{j\in \mathcal N} {g_N}_j(\varphi_j,\varphi_i)_{L^2(\Gamma_{N,h})} 
\label{eq:weak2N_Nitsche_discrete}    
\end{align}
where the quantities $f_i$, $g^N_i$ denote the node values of the source data and of the Neumann boundary conditions, respectively. Here we rewrite \eqref{eq:weak2N_Nitsche_discrete}, as follows
\begin{align} \nonumber
    &\sum_{j\in \mathcal N} u_j\underbrace{(\nabla\varphi_j,\nabla\varphi_i)_{L^2(\Omega_h)}}_{\rm S}  = \sum_{j\in \mathcal N} f_j \underbrace{(\varphi_j,\varphi_i)_{L^2(\Omega_h)}}_{\rm M} + \sum_{j\in \mathcal N} {g_N}_j\underbrace{(\varphi_j,\varphi_i)_{L^2(\Gamma_{N,h})} }_{\rm N_{\Gamma}} 
\end{align}
This results in the discrete system 
\begin{equation}
	\sum_{j\in\nodes} a^N_{ij} u_j = f^N_i, \quad i\in\nodes
	\label{eq:system2N_A}
\end{equation}
where $A^N = \{a^N_{ij} \}$ and $F^N = \{f_i \}$, and
\begin{align}
\label{coeffN_A}
A^N = {\rm S},\\
F^N = {\rm M} f + {\rm N_{\Gamma}}\, g_N,
\end{align}
where the entries of the matrices are 
\begin{align}
    \label{eq_coeff_aN_A}
	\displaystyle	a^N_{ij}  & = (\nabla\varphi_i,\nabla\varphi_j)_{L^2(\Omega_h)},\\ \nonumber
	\displaystyle	f^N_i & = (f,\vphi_i)_{L^2(\Omega_h)} + (g_N,\vphi_i)_{L^2(\Gamma_{N,h})}.
\end{align}

We now consider the mixed problem in two dimensions:
\begin{equation}
	-\Delta u = f \quad {\rm in}\>  \Omega, \quad u = g_D \> {\rm on}\>  \Gamma_D, \> 
	\pad{u}{n} = g_N \> {\rm on}\>  \Gamma_N
	\label{eq2DM_A}	
\end{equation}
where the boundary $\partial\Omega$ is partitioned into a Dirichlet boundary $\Gamma_D$ and 
a Neumann boundary $\Gamma_N$, $\Gamma_D\cup \Gamma_N = \partial \Omega$, $\Gamma_D\cap\Gamma_N=\emptyset$.

The penalised method is obtained by minimizing the following functional, which differs from the original one because of the addition of a penalisation term on the Dirichlet boundary $\Gamma_D$:
\begin{equation}
	J_{D\lambda}[u] = \int_\Omega \left[\frac12 (\nabla u)^2 - f u\right]\, d\Omega + 
	\frac{\lambda}2\int_{\Gamma_D} (u-g_D)^2\, ds
	\label{JD2LM[u]}
\end{equation}
Setting to zero the first variation of the functional leads to the following penalised variational formulation of the problem: 
\begin{align}
	\int_\Omega\nabla u\cdot \nabla v\, d\Omega 
	 +  \int_{\Gamma_D} \Big [\lambda\cdot (u-g_D) & -\pad{u}{n}\Big ]v \, ds \nonumber \\
	=  
	& \int_\Omega f v \, d\Omega 
	+  \int_{\Gamma_N} g_N v \, ds 
	\label{eq:weak2M}
\end{align}
where we made use of Neumann boundary conditions by replacing $\displaystyle \frac{\partial u}{\partial n}$ with $g_N$ on $\Gamma_N$. 
Also in two space dimensions, it is possible to make the problem symmetric by adopting the Nitsche technique, adding and subtracting the same quantity, as follows
\begin{align}
	\int_\Omega\nabla u\cdot \nabla v\, d\Omega 
	 + &  \int_{\Gamma_D} \Big [\lambda\cdot (u-g_D) -\pad{u}{n}\Big ]v \, ds - \int_{\Gamma_D}  u \pad{v}{n} \, ds \nonumber \\
	=  
	& \int_\Omega f v \, d\Omega 
	+  \int_{\Gamma_N} g_N v \, ds - \int_{\Gamma_D} g_D \pad{v}{n} \, ds.
	\label{eq:weak2M_Nitsche}
\end{align}

Notice that for purely Neumann problems there is no penalisation term, since $\Gamma_N = \Gamma$ and $\Gamma_D = \emptyset$.
For such problems, we assume the standard compatibility condition is satisfied, i.e.
\begin{equation}
    \int_\Omega f\, d\Omega + \int_\Gamma g_N\, d\Gamma = 0.
\end{equation}
The discretization of the problem is obtained in the usual way, by replacing the domain $\Omega$ with its polygonal approximation $\Omega_h$, whose boundary is $\Gamma_h$.

We now proceed to consider the variational formulation \eqref{eq:weak2M_Nitsche} over the discrete space $V_h$, and using the standard basis, the above equations become
\begin{align}
    &\sum_{j\in \mathcal N} u_j(\nabla\varphi_i,\nabla\varphi_j)_{L^2(\Omega_h)}  - \sum_{j\in \mathcal N} u_j\left(\pad{\varphi_i}{n},\varphi_j\right)_{L^2(\Gamma_{D,h})}  - \sum_{j\in \mathcal N} u_j\left(\varphi_j,\pad{\varphi_i}{n}\right)_{L^2(\Gamma_{D,h})} \nonumber \\ + &\lambda\sum_{j\in \mathcal N} u_j(\varphi_j,\varphi_i)_{L^2(\Gamma_{D,h})} = \sum_{j\in \mathcal N} f_j (\varphi_j,\varphi_i)_{L^2(\Omega_h)} + \sum_{j\in \mathcal N} {g_N}_j(\varphi_j,\varphi_i)_{L^2(\Gamma_{N,h})} \\ \nonumber - & \sum_{j\in \mathcal N} {g_D}_j\left(\varphi_j,\pad{\varphi_i}{n}\right)_{L^2(\Gamma_{D,h})} + \lambda\sum_{j\in \mathcal N} {g_D}_j(\varphi_j,\varphi_i)_{L^2(\Gamma_{D,h})}
\label{eq:weak2M_Nitsche_discrete}    
\end{align}
where the quantities $f_i$, $g^D_i$, $g^N_i$ denote the node values of the source data and of the mixed boundary conditions, respectively. Here we rewrite \eqref{eq:weak2M_Nitsche_discrete}, as follows
\begin{align} \nonumber
    &\sum_{j\in \mathcal N} u_j\left(\underbrace{(\nabla\varphi_j,\nabla\varphi_i)_{L^2(\Omega_h)}}_{\rm S} - \underbrace{\left(\left(\pad{\varphi_j}{n},\varphi_i\right)_{L^2(\Gamma_{D,h})}  + \left(\varphi_j,\pad{\varphi_i}{n}\right)_{L^2(\Gamma_{D,h})} \right)}_{\rm S_T} \right. \\ + & \left. \lambda\underbrace{(\varphi_j,\varphi_i)_{L^2(\Gamma_{D,h})}}_{\rm P_{\Gamma_D}} \right)  = \sum_{j\in \mathcal N} f_j \underbrace{(\varphi_j,\varphi_i)_{L^2(\Omega_h)}}_{\rm M} + \sum_{j\in \mathcal N} {g_N}_j\underbrace{(\varphi_j,\varphi_i)_{L^2(\Gamma_{N,h})} }_{\rm N_{\Gamma_N}} \\ &+ \sum_{j\in \mathcal N} {g_D}_j\left(-\underbrace{\left(\varphi_j,\pad{\varphi_i}{n}\right)_{L^2(\Gamma_{D,h})}}_{\rm D_{\Gamma_D}} \!\!\!\!\!\!\!+ \underbrace{\lambda(\varphi_j,\varphi_i)_{L^2(\Gamma_{D,h})}}_{\rm P_{\Gamma_D}} \right)
\end{align}
This results in the discrete system 
\begin{equation}
	\sum_{j\in\nodes} a^M_{ij} u_j = f^M_i, \quad i\in\nodes
	\label{eq:system2M_A}
\end{equation}
where $A^M = \{a^M_{ij} \}$ and $F^M = \{f_i \}$, and
\begin{align}
\label{coeffM2_A}
A^M = {\rm S - S_T} +\lambda {\rm P_{\Gamma_D}},\\ \nonumber
F^M = {\rm M} f + (\lambda{\rm P_{\Gamma_D}- D_{\Gamma_D}}) g_D + {\rm N_{\Gamma_N}}\, g_N,
\end{align}
where the entries of the matrices are 
\begin{equation}
    \begin{array}{c}
    \label{eq_coeff_aM_A}
	\displaystyle	a^M_{ij}  = (\nabla\varphi_i,\nabla\varphi_j)_{L^2(\Omega_h)} 
	+ \lambda (\varphi_i,\varphi_j)_{L^2(\Gamma_{D,h})} 
	- (n\cdot\nabla\vphi_j,\vphi_i)_{L^2(\Gamma_{D,h})} \\ - (\vphi_j,n\cdot\nabla\vphi_i)_{L^2(\Gamma_{D,h})},\\[4mm]
	\displaystyle	f^M_i  = \lambda (g_D,\varphi_i)_{L^2(\Gamma_{D,h})} - (g_D,n\cdot\nabla\vphi_i)_{L^2(\Gamma_{D,h})}
	+ (g_N,\vphi_i)_{L^2(\Gamma_{N,h})} \\ + (f,\vphi_i)_{L^2(\Omega_h)}.
\end{array}
\end{equation}

\section{Computational details of the one-dimensional problem}
\label{Appendix:B}
Referring to Fig.~\ref{fig:1D_setup}, this section is dedicated to the explicit calculation of the entries for the matrices introduced in Section~\ref{sec:one_dim_discr}. Our primary objective is to determine all the non-zero entries of the matrices 
$A^D$ and $A^M$. While the computation of these entries is trivial in the interior region of the domain, the complexity arises in the proximity of the boundaries. Therefore, we focus our attention on the scalar product of following basis functions: 
\begin{equation*}
\begin{array}{c}
\displaystyle	(\varphi_0',\varphi_0')_{L^2([a,b])} = \frac{\theta_1}{h},\quad 
	(\varphi_0',\varphi_1')_{L^2([a,b])} = -\frac{\theta_1}{h},\quad 
	(\varphi_1',\varphi_1')_{L^2([a,b])} = \frac{1+\theta_1}{h}, \\[4mm] 
\displaystyle	(\varphi_{N-1}',\varphi_{N-1}')_{L^2([a,b])} = \frac{1+\theta_2}{h},\quad 
	(\varphi_{N-1}',\varphi_N')_{L^2([a,b])} = -\frac{\theta_2}{h}, \quad	
	(\varphi_N',\varphi_N')_{L^2([a,b])} = \frac{\theta_2}{h}, \\[4mm]
\displaystyle	(\varphi_{i-1}',\varphi_i')_{L^2([a,b])} = -\frac{1}{h}, \> i=1,\ldots,N \quad
	(\varphi_i',\varphi_i')_{L^2([a,b])} = \frac{2}{h}, i=2,\ldots,N-2,\\[4mm]
\displaystyle	\varphi_0(a) = \theta_1,\quad \varphi_1(a) = 1-\theta_1, \quad
	\varphi_N(b) = \theta_2,\quad \varphi_{N-1}(b) = 1-\theta_2, \\[4mm]
\displaystyle	\varphi'_0(a) = -1/h,\quad \varphi'_1(a) = 1/h, \quad
	\varphi'_N(b) = 1/h, \quad \varphi'_{N-1}(b) = -1/h.
\end{array}
\end{equation*}

The rest of the source terms can be computed as follows. The term related to the mass matrix is
 \begin{equation*}
 \begin{array}{c}
\displaystyle	(f,\varphi_i)_{L^2([a,b])} \approx \sum_{j = 0}^N f_j(\varphi_j,\varphi_i)_{L^2([a,b])} = \frac{h}{6}\left(f_{i-1} + 4f_i + f_{i+1} \right), \quad i=2,\ldots,n-2. \\

\displaystyle	(f,\varphi_0)_{L^2([a,b])} \approx \sum_{j = 0}^N f_j(\varphi_j,\varphi_0)_{L^2([a,b])} = h\left(\frac{\theta_1^3}{3}f_0 +  \left(\frac{\theta_1^2}{2}-\frac{\theta_1^3}{3}\right)f_1\right), \\[5mm]

\displaystyle	(f,\varphi_1)_{L^2([a,b])} \approx \sum_{j = 0}^N f_j(\varphi_j,\varphi_1)_{L^2([a,b])} = h\left(\left(\frac{\theta_1^2}{2}-\frac{\theta_1^3}{3}\right)f_0 + \frac{1+\theta_1^3}{3}f_1\right),    \\[5mm]

\displaystyle	(f,\varphi_{N-1})_{L^2([a,b])} \approx \sum_{j = 0}^N f_j(\varphi_j,\varphi_{N-1})_{L^2([a,b])} = h\left( \left(\frac{1+\theta_2^3}{3}-\theta_2^2+\theta_2\right) f_{N-1} + \left(\frac{\theta_2^2}{2}-\frac{\theta_2^3}{3}\right) f_N\right)  \\[5mm]

\displaystyle	(f,\varphi_{N})_{L^2([a,b])} \approx \sum_{j = 0}^N f_j(\varphi_j,\varphi_{N-1})_{L^2([a,b])} = h\left( 
\left(\frac{\theta_2^2}{2}-\frac{\theta_2^3}{3}\right) f_{N-1} + \frac{\theta_2^3}{3} f_N\right).
\end{array}
\end{equation*}

Regarding the penalisation term
\begin{equation*}
    \begin{array}{c}
\displaystyle u_a\varphi_0(a) \approx \sum_{j = 0}^N{u_a} \varphi_j(a)\varphi_0(a) = \theta_1^2 {u_a}_0 +  \theta_1(1-\theta_1) {u_a}_1 \\
\displaystyle u_a\varphi_1(a) \approx \sum_{j = 0}^N{u_a} \varphi_j(a)\varphi_1(a)= \theta_1(1-\theta_1) {u_a}_0  + (1-\theta_1^2)  {u_a}_1,
    \end{array}
\end{equation*}
analogously, we find the coefficients for the penalisation term in $x = b$.




\end{document}